\documentclass{amsart}

\usepackage{amssymb}
\usepackage{eepic}
\usepackage{color}

\allowdisplaybreaks

\numberwithin{equation}{section}

\newtheorem{theorem}{Theorem}[section]
\newtheorem{lemma}[theorem]{Lemma}

\newtheorem{corollary}[theorem]{Corollary}
\newtheorem{remark}[theorem]{Remark}
\newtheorem{problem}[theorem]{Problem}

\newtheorem{Atheorem}{Theorem}

\newtheorem{Aremark}[Atheorem]{Remark}

\newtheorem{TheoA}{Theorem A}
\newtheorem{TheoB}{Theorem B}
\newtheorem{TheoC}{Theorem C}
\newtheorem{TheoD}{Theorem D}

\newcommand{\Z}{\mathbb{Z}}
\newcommand{\R}{\mathbb{R}}
\newcommand{\C}{\mathbb{C}}

\newcommand{\summ}{\sum\nolimits}

\def\ten{\otimes}
\def\G{\mathrm{G}}
\def\1{\mathbf{1}}
\def\H{\mathcal{H}}
\def\Q{\mathcal{Q}}
\def\M{\mathcal{M}}
\def\A{\mathcal{A}}
\def\RR{\mathcal{R}}

\def\T{\mathcal{S}}
\def\S{\mathcal{S}}

\def\V{\mathrm{\mathcal{L}(G)}}
\def\Rd{\R_\mathrm{disc}}
\def\hatRd{\widehat{\R}_\mathrm{disc}}

\newcommand{\dem}{\noindent {\bf Proof. }}
\newcommand{\demA}{\noindent {\bf Proof of Theorem A. }}
\newcommand{\demB}{\noindent {\bf Proof of Theorem B. }}

\newcommand{\demD}{\noindent {\bf Proof of Theorem D. }}

\newcommand{\fin}{\hspace*{\fill} $\square$ \vskip0.2cm}

\def\esssup{\mathop{\mathrm{ess \, sup}}}

\addtocounter{tocdepth}{0}

\def\mean{- \hskip-10.5pt \int}

\def\meann{- \hskip-10.7pt \int}

\addtocounter{tocdepth}{0}

\begin{document}


\title[Smooth multipliers on group von Neumann algebras]{Smooth Fourier multipliers \\ on group von Neumann algebras}

\author[Junge, Mei, Parcet]
{Marius Junge, Tao Mei, Javier Parcet}

\maketitle

\begin{abstract}
\hskip-3pt We investigate Fourier multipliers on the compact dual of arbitrary discrete groups. Our main result is a H\"ormander-Mihlin multiplier theorem for finite-dimensional cocycles with optimal smoothness condition. We also find Littlewood-Paley type inequalities in group von Neumann algebras, prove $L_p$ estimates for noncommutative Riesz transforms and characterize $L_\infty \to \mathrm{BMO}$ boundedness for radial Fourier multipliers. The key novelties of our approach are to exploit group cocycles and cross products in Fourier multiplier theory in conjunction with BMO spaces associated to semigroups of operators and a noncommutative generalization of Calder\'on-Zygmund theory. 
\end{abstract}

\addtolength{\parskip}{+1ex}

\section*{{\bf Introduction and main results}}

Convergence of Fourier series and norm estimates for Fourier multipliers are central in harmonic analysis. As far as Calder\'on-Zygmund methods are involved very few results have been successfully transferred to other noneuclidean topological groups. In this paper we study smooth Fourier multipliers frequency supported by an arbitrary discrete group. For instance, the frequency group associated to the $n$-dimensional torus is the integer lattice $\Z^n$ and ---by de Leeuw's compactification theorem \cite{dL}--- we may impose the discrete topology on the frequency group of $\R^n$ and still obtain the same family of $L_p$-multipliers. What can we say about $L_p$-boundedness of Fourier multipliers for arbitrary discrete groups? What do we mean by smoothness of the multiplier in that case? Basic unexplored examples include duals of Cantor cubes, other discrete abelian groups, finite groups of large cardinality, the discretized Heisenberg group, free groups... For nonabelian discrete groups, the compact dual is a quantum group whose underlying space is a group von Neumann algebra. These algebras are widely accepted and very well understood in noncommutative geometry \cite{Con} and operator algebra \cite{KR}. In harmonic analysis this approach was first considered in the ground-breaking results of Haagerup \cite{H} and Cowling/Haagerup \cite{CH} on the approximation property of group von Neumann algebras. Up to isolated contributions \cite{Ha,NR}, the $L_p$-theory for Fourier multipliers on these algebras is very much unexplored. 

Let $\G$ be a discrete group with left regular representation $\lambda: \G \to \mathcal{B}(\ell_2(\G))$ given by $\lambda(g) \delta_h = \delta_{gh}$, where the $\delta_g$'s form the unit vector basis of $\ell_2(\G)$. Write $\mathcal{L}(\G)$ for its group von Neumann algebra, the weak operator closure of the linear span of $\lambda(\G)$. Given $f \in \mathcal{L}(\G)$, we consider the standard normalized trace $\tau_\G(f) = \langle \delta_e, f \delta_e \rangle$ where $e$ denotes the identity of $\G$. Any such $f$ has a Fourier series $$\sum_{g \in \G} \widehat{f}(g) \lambda(g) \quad \mbox{with} \quad \widehat{f}(g) = \tau_\G(f \lambda(g^{-1})) \quad \mbox{so that} \quad \tau_\G(f) = \widehat{f}(e).$$ Let $L_p(\widehat{\mathbb{G}}) = L_p(\mathcal{L}(\mathrm{G}), \tau_\G)$ denote the $L_p$ space over the noncommutative measure space $(\mathcal{L}(\mathrm{G}), \tau_\G)$ ---so called noncommutative $L_p$ spaces--- equipped with the norm $$\|f\|_p = \Big\| \summ_g \widehat{f}(g) \lambda(g) \Big\|_p = \Big( \tau_\G \Big[ \hskip1pt \Big| \summ_g \widehat{f}(g) \lambda(g) \Big|^p \hskip1pt \Big] \Big)^\frac1p.$$ We invite the reader to check that $L_p(\widehat{\mathbb{G}}) = L_p(\mathbb{T}^n)$ for $\G = \Z^n$, after identifying $\lambda(k)$ with $e^{2\pi i \langle k, \cdot \rangle}$. In general, the absolute value and the power $p$ are obtained from functional calculus for this (unbounded) operator on $\ell_2(\G)$. A Fourier multiplier is then given by $$T_m: \summ_g \widehat{f}(g) \lambda(g) \mapsto \summ_g m_g \widehat{f}(g) \lambda(g).$$ Going back to $\mathrm{G} = \Z^n$, we find Fourier multipliers on the $n$-torus. We will say that any smooth function $\widetilde{m}: \R^n \to \C$ is a lifting multiplier for $m$ whenever its restriction to $\Z^n$ coincides with $m$. According to de Leeuw's restriction theorem \cite{dL}, the $L_p$ boundedness of $T_m$ follows whenever there exists a well-behaved lifting multiplier $\widetilde{m}$ defining an $L_p$-bounded map in the ambient space $\R^n$. In particular, if $1 < p < \infty$ it suffices to check the H\"ormander-Mihlin smoothness condition \cite{Ho,Mi} $$\widetilde{m} \in \mathcal{C}^{[\frac{n}{2}]+1}(\R^n \setminus \{0\}) \quad \mbox{and} \quad \big| \partial_\xi^\beta \, \widetilde{m}(\xi) \big| \, \lesssim \, |\xi|^{-|\beta|} \quad \mbox{for all} \quad |\beta| \le \left[\frac{n}{2}\right]+1.$$ In the context of Lie groups, the role of $\R^n$ is replaced by the Lie algebra to find similar formulations. Recall the remarkable work by M\"uller, Ricci, Stein and their coauthors on nilpotent groups, see \cite{M1,MRS,RR,RS2} and the references therein. A fundamental goal in this paper is to give sufficient differentiability conditions for the $L_p$ boundedness of multipliers on the compact dual of discrete groups. Unlike for $\Z^n$, there is no standard differential structure to construct smooth lifting multipliers for an arbitrary discrete $\G$. The main novelty in our approach is to identify the right endpoint spaces ---intrinsic BMO's over certain semigroups--- using a broader interpretation of tangent spaces in terms of length functions and cocycles.

An affine representation of $\G$ is an orthogonal representation $\alpha: \G \to O(\H)$ over a real Hilbert space $\H$ together with a mapping $b: \G \to \H$ satisfying the cocycle law $b(gh) = \alpha_g(b(h)) + b(g)$. Given an affine representation $(\H,\alpha,b)$, the function $\psi(g) = \langle b(g), b(g) \rangle_\H$ vanishes at the identity $e$, it is symmetric $\psi(g) = \psi(g^{-1})$ and conditionally negative, i.e. $\sum_g \beta_g = 0 \Rightarrow \sum_{g,h} \overline{\beta}_g \beta_h \psi(g^{-1}h) \le 0$. In this paper any $\psi: \G \to \R_+$ satisfying the properties above will be called a \emph{length}. Conversely, we know from Schoenberg's theorem \cite{Sc} that any length $\psi$ determines a precise affine representation $(\H_\psi, \alpha_\psi, b_\psi)$ satisfying $\psi(g) = \langle b_\psi(g), b_\psi(g) \rangle_{\H_\psi}$, more details will be given in the body of the paper. H\"ormander-Mihlin and de Leeuw classical theorems are formulated in terms of the trivial cocycle coming from the heat semigroup. We propose the Hilbert spaces $\H_\psi$ as cocycle substitutes of the Lie algebra. Here is a fairly simple formulation ---see also Theorems B and \ref{MainHorMult}--- of our cocycle form of H\"ormander-Mihlin theorem, valid for group von Neumann algebras.

\begin{TheoA}
Let $\G$ be a discrete group equipped with any length $\psi: \G \to \R_+$ satisfying $\dim \H_\psi = n < \infty$. Let $T_m: \sum_g \widehat{f}(g) \lambda(g) \mapsto \sum_g m_g \widehat{f}(g) \lambda(g)$ and assume there exists $\varepsilon > 0$ and $\widetilde{m}: \H_\psi \to \C$ satisfying $m_g = \widetilde{m}(b_\psi(g))$ and $$\big| \partial_\xi^\beta \widetilde{m}(\xi) \big| \, \lesssim \, \min \Big\{ |\xi|^{- |\beta|+\varepsilon}, |\xi|^{- |\beta|-\varepsilon} \Big\} \quad \mbox{for} \quad |\beta| \le \Big[\frac{n}{2}\Big]+1.$$ Then, $T_m: L_p(\widehat{\mathbb{G}}) \to L_p(\widehat{\mathbb{G}})$ is a completely bounded multiplier for all $1 < p < \infty$.
\end{TheoA}

Completely bounded (\emph{cb}) means that $T_m \ten id$ is a multiplier on $L_{p}(\widehat{\mathbb{G\times H}})$ for every discrete group $\mathrm{H}$. Theorem A is sharp in terms of the number of derivatives related to the cocycle dimension. On the contrary, the additional $\varepsilon$ is a prize we pay for noncommutativity. It can be removed under assumptions like
\begin{itemize}
\item[i)] $\G$ is abelian,

\item[ii)] $b_\psi(\G)$ is a lattice in $\R^n$,

\item[iii)] $\alpha_\psi(\G)$ is a finite subgroup of $O(n)$,

\item[iv)]The multiplier is $\psi$-radial, i.e. $m_g = h(\psi(g))$.
\end{itemize}
Theorem A is a cocycle extension of the Mihlin multiplier theorem, more than merely a noncommutative form of it. Indeed, it provides new results even for finite or abelian groups. For instance, we may find low dimensional injective cocycles for finite groups of large cardinality, like $\mathbb{Z}_n$ or the symmetric permutation groups $\mathrm{S}_n$ where we find injective cocycles with $\dim \H_\psi = 2 << n$ and $\dim \H_\psi = n << n!$ respectively. To the best of our knowledge, there are no similar results for finite groups. In the Euclidean (non-discrete) context of $\R^n$, we will also analyze the classical H\"ormander-Mihlin and de Leeuw's theorems under our perspective, which naturally gives rise to a large family of exotic Fourier multipliers.

The presence of $\varepsilon > 0$ in Theorem A excludes some central examples, like the $\psi$-\emph{directional Riesz transforms} which are naturally defined for $\eta \in \H_\psi$ ---recalling that $\psi(g) = \langle b_\psi(g), b_\psi(g) \rangle_{\H_\psi}$--- as $$R_\eta \Big( \summ_g \widehat{f}(g) \lambda(g) \Big) = - i \summ_g \frac{\langle b_\psi(g), \eta \rangle_{\H_\psi}}{\sqrt{\psi(g)}} \widehat{f}(g) \lambda(g).$$ These multipliers are covered by Theorem B below, which exploits the semigroup $\T_\psi = (S_{\psi,t})_{t \ge 0}$ with $S_{\psi,t}: \lambda(g) \mapsto e^{-t \psi(g)} \lambda(g)$ and the BMO space constructed with it. Namely, define $$\mathrm{BMO}_{\T_\psi} \, = \, \Big\{ f \in L_2^\circ(\widehat{\mathbb{G}}) \ \big| \ \|f\|_{\mathrm{BMO}_{\T_\psi}} = \max \big\{\|f\|_{\mathrm{BMO}_{\T_\psi}^c}, \|f^*\|_{\mathrm{BMO}_{\T_\psi}^c}\big\} < \infty \Big\},$$ where $f \in L_2^\circ(\widehat{\mathbb{G}}) \subset L_2(\widehat{\mathbb{G}})$ if $\widehat{f}(g)=0$ whenever $\psi(g)=0$ and $$\|f\|_{\mathrm{BMO}_{\T_\psi}^c} \, = \, \sup_{t > 0} \Big\| \Big( S_{\psi,t}|f|^2 - |S_{\psi,t}f|^2 \Big)^\frac12 \Big\|_{\V} \quad \mbox{with} \quad |f|^2 = f^* f.$$ The following result proves the $L_p$-boundedness of $\psi$-directional Riesz transforms.

\begin{TheoB}
Given $(\G, \psi)$ and $$T_m: \summ_g \widehat{f}(g) \lambda(g) \mapsto \summ_g m_g \widehat{f}(g) \lambda(g)$$ as above, let $\widetilde{m}: \H_\psi \to \C$ be a lifting multiplier for $m = \widetilde{m} \circ b_\psi$ such that
\begin{itemize}
\item[i)] $L_2$-row/column condition $$\big\| \widetilde{m} \big\|_{schur}
\, = \, \inf_{\begin{subarray}{c}
\widetilde{m}(\alpha_{\psi,g}(\xi)) = \langle A_\xi, B_g
\rangle_\mathcal{K} \\ (\xi,g) \in \H_\psi \times \G \\ \mathcal{K} \ \mathrm{Hilbert} 
\end{subarray}} \Big( \sup_{\xi \in \H_\psi} \|A_\xi\|_\mathcal{K} \,
\sup_{g \in \G} \|B_g\|_\mathcal{K} \Big) \, < \, \infty.$$

\item[ii)] H\"ormander-Mihlin smoothness $$\qquad \widetilde{m}
\in \mathcal{C}^{n+2}(\R^n \setminus\{0\}) \ \ \mbox{and} \ \
\big| \partial_\xi^\beta \widetilde{m}(\xi) \big| \, \lesssim \, |\xi|^{- |\beta|} \ \ \mbox{for all} \ \ |\beta| \le n+2.$$
\end{itemize}
Then, $T_m:  \V \stackrel{cb}{\longrightarrow} \mathrm{BMO}_{\T_\psi}$ and $T_m: L_p(\widehat{\mathbb{G}}) \stackrel{cb}{\longrightarrow} L_p(\widehat{\mathbb{G}})$ for every $1 < p < \infty$.
\end{TheoB}

As Theorem A, the above result holds for any finite-dimensional cocycle. In the particular case of Riesz transforms, it is meaningful to wonder if Riesz transforms associated to infinite-dimensional cocycles are also bounded. This follows from a dimension free estimate \cite{JMP00} with applications for other multipliers. Both Theorems A and B can be used either to construct Fourier multipliers or to test the $L_p$ boundedness of a fixed multiplier. The real challenge in the latter case is to find the right length/cocycle $b_\psi$ and the lifting $\widetilde{m}$ such that $m_g = \widetilde{m}(b_\psi(g))$. This is exactly the topic of Fefferman's \emph{smooth interpolation of data} \cite{Fe1,Fe2,Fe3} relative to the set $b_\psi(\G)$. If $\Delta_\psi = \inf_{g \neq h} \| b_\psi(g)-b_\psi(h) \|_{\H_\psi}^2 >0$, we call the cocycle $b_\psi$ well-separated. The H\"ormander-Mihlin theorem for the $n$-torus corresponds to the standard cocycle $\Z^n \subset \R^n$ with the trivial action $\alpha$. Up to changes of basis, it is the only finite-dimensional, injective, well-separated cocycle of $\Z^n$. Accordingly we call $b_\psi$ standard if it is injective and well-separated. By a classical theorem of Bieberbach \cite{Bie}, a standard finite-dimensional cocycle of $\G$ can only exist if $\G$ is virtually abelian. This excludes for instance infinite discrete groups with Kazhdan's property $(\mathrm{T})$. The novelty in our approach is to allow for \emph{clustering} in the set $b_\psi(\G)$ and thus go beyond the class of virtually abelian groups. Of course, our hypotheses lead to look for disperse clouds $b_\psi(\G)$ living in low dimensional spaces $\H_\psi$. We refer to Paragraph \ref{Ending} for a description of the intriguing interplay between these \lq\lq competing" requirements.

Our methods also lead to Littlewood-Paley type estimates. Square function inequalities in noncommutative $L_p$ spaces require to combine the so-called row and column square functions appropriately. Given a von Neumann algebra $\M$, the spaces $L_p(\M; \ell_2^r)$ and $L_p(\M; \ell_2^c)$ are the closure of finite sequences $(f_j)$ in $L_p(\M)$ with respect to the norms $$\Big\| \big( \summ_j f_j f_j^* \big) ^\frac12 \Big\|_{L_p(\M)} \quad \mbox{and} \quad \Big\| \big( \summ_j f_j^* f_j \big) ^\frac12 \Big\|_{L_p(\M)}$$ respectively. The right combination $$L_p(\M; \ell_{rc}^2) = \begin{cases} L_p(\M; \ell_2^r) + L_p(\M; \ell_2^c) & (1 \le p \le 2) \\ L_p(\M; \ell_2^r) \cap \hskip1pt L_p(\M; \ell_2^c) & (2 \le p \le \infty) \end{cases}$$ was discovered in the formulation of the noncommutative Khintchine inequalities \cite{Lu,LuP0}. Given a length function $\psi: \G \to \R_+$ arising from an $n$-dimensional cocycle, consider a sequence of functions $(h_j)_{j \ge 1}$ in $\mathcal{C}^{[\frac{n}{2}]+1}(\R_+ \setminus \{0\})$ such that $\sum_j | \frac{d^k}{d\xi^k} h_j(\xi) |^2 \le c_n |\xi|^{-2k}$ for all $k \le [\frac{n}{2}] +1$ and define $$L_{\psi, j} f = \summ_g h_j(\psi(g)) \widehat{f}(g) \lambda(g) \quad \mbox{and} \quad \Lambda_\psi(f) = \summ_j L_{\psi,j} f \otimes \delta_j.$$ 

\vskip-2pt

\begin{TheoC}
If $1 < p < \infty$, we have 
\begin{itemize}
\item[i)] $\Lambda_\psi: L_p(\widehat{\mathbb{G}}) \stackrel{cb}{\longrightarrow} L_p(\widehat{\mathbb{G}}; \ell_{rc}^2)$, that is $$\Big\| \summ_j L_{\psi,j} f \otimes \delta_j \Big\|_{L_p(\widehat{\mathbb{G}}; \ell_{rc}^2)} \, \le_{cb} \, c_p \, \|f\|_{L_p(\widehat{\mathbb{G}})}.$$ 

\item[ii)] Additionally, if $\sum_j |h_j(\xi)|^2 = 1$ for all $\xi \neq 0$ $$\|f\|_{L_p(\widehat{\mathbb{G}})} \, \le_{cb} \, c_p \, \Big\| \summ_j L_{\psi,j}
 f \otimes \delta_j \Big\|_{L_p(\widehat{\mathbb{G}}; \ell_{rc}^2)},$$ so that $\Lambda_\psi: L_p(\widehat{\mathbb{G}}) \to L_p(\widehat{\mathbb{G}}; \ell_{rc}^2)$ becomes a complete embedding.
\end{itemize}
\end{TheoC}

Theorem C is a crucial tool for the main results in \cite{JMP00,PRo}, see Theorem \ref{ThNCLP} for a formulation including $L_\infty \to \mathrm{BMO}$ type estimates. A noncommutative form of Calder\'on-Zygmund theory requires to find suitable substitutes for the interplay metric/measure in commutative spaces. Theorems A, B and C emerge from CZO's on von Neumann algebras $\RR \rtimes \G$, where $\G$ acts  on a von Neumann algebra $\RR$ which factors as a tensor product of the Bohr compactification of $\R^n$ and some other noncommutative measure space $(\M, \tau)$. The key link with our main results comes from the \emph{intertwining identities} $$\pi_\psi \circ S_{\psi,t} = \big( S_t \rtimes id_\G \big) \circ \pi_\psi \quad \mbox{and} \quad \pi_\psi \circ T_m = \big( T_{\widetilde{m}} \rtimes id_{\V} \big) \circ \pi_\psi,$$ where $(S_t)_{t \ge 0}$ denotes the heat semigroup on $\H_\psi$ ---equipped with the discrete topology--- and $\pi_\psi: \V \to \mathcal{L}(\H_\psi) \rtimes_{\alpha_\psi} \G$ is the $*$-homomorphism determined by $\lambda(g) \mapsto \exp(2\pi i \langle b_\psi(g), \, \cdot \, \rangle_{\H_\psi}) \rtimes \lambda(g)$. The first intertwining identity yields an embedding $\mathrm{BMO}_{\T_\psi} \to \mathrm{BMO}_{\T_\rtimes}$ with $S_{\rtimes,t} = S_t \rtimes id_\G$. This explains our interest on BMO spaces over semigroups of cp maps. In the classical theory we find BMO spaces associated to a metric or a martingale filtration. Duong and Yan \cite{DY1,DY2} extended it to certain semigroups ---see also \cite{ADM,HM}--- but still imposing the existence of a metric in the underlying space, something that a priori we do not have at our disposal. Interpolation results with $L_p$ spaces  \cite{JM2} are deduced from the theory of noncommutative martingales with continuous index set \cite{JPe} and a theory of Markov dilations with continuous path \cite{JRS}.

Let $\widehat{\R}_{\mathrm{disc}}^n$ be the Bohr compactification of $\R^n$, the Pontryagin dual of $\R^n$ equipped with the discrete topology. Since $\psi(g) = \| b_\psi(g) \|_{\H_\psi}^2$ a multiplier $m: \G \to \C$ is called $\psi$-radial if $m_g = h(\psi(g))$ for some $h: \R_+ \to \C$, so that we find a lifting $\widetilde{m}$ which is radial on $\H_\psi$. We will use the little Grothendieck theorem \cite{Gro,PBook,PBull} for the following characterization of $L_\infty \to \mathrm{BMO}$ boundedness for $\psi$-radial multipliers.

\begin{TheoD}
If $h: \R_+ \to \C$, $\mathrm{TFAE}$
\begin{itemize}
\item[i)] $T_{h \circ | \, |^2}: L_\infty(\R^n) \to
\mathrm{BMO}_{\R^n}$ bounded,

\vskip1pt

\item[ii)] $T_{h \circ | \, |^2}: L_\infty(\hatRd^n) \to \mathrm{BMO}_{\T'}$ bounded,

\vskip3pt

\item[iii)] $T_{h \circ \psi} \! : \mathcal{L}(\G)
\! \to \! \mathrm{BMO}_{\T_\psi}$ \hskip-1pt bounded
\hskip-1pt for \hskip-1pt all $\G$ \hskip-1pt discrete \hskip-1pt
with \hskip-1pt $\dim \H_\psi = n$,
\end{itemize}
with $\T'$ the heat semigroup on $\hatRd^n$. Moreover, \emph{ii)} $\Leftrightarrow$ \emph{iii)} still holds for $n=\infty$.
\end{TheoD}

Our results in this paper show some impact of cohomology in Fourier multiplier theory. Tools from classical harmonic analysis impose that $\dim \H_\psi < \infty$, but  many interesting cocycles are constructed on infinite-dimensional Hilbert spaces. In fact certain exotic groups like the Tarski monster or some Burnside groups, do not admit finite-dimensional cocycles at all. Fortunately, Theorem D for radial multipliers and the dimension free estimates for Riesz transforms in \cite{JMP00} indicate that this is not the end of Fourier multiplier theory. Since $A_\psi(\lambda(g)) = \psi(g) \lambda(g)$ generates $\T_\psi$, radial multipliers are of the form $h(A_\psi)$, already considered by McIntosh's $H_\infty$-calculus for analytic $h$. Theorem A imposes considerably weaker conditions and Theorem D provides new Fourier multipliers even for infinite-dimensional cocycles. The imaginary powers $\psi(g)^{is}$ and other examples of Laplace transforms are included as shown in \cite{JM2}. As an illustration in $\G=\R^n$, the lengths $$\psi(\xi) = 1 - \int_{\R^n} f_0(x) f_0(x-\xi) \, dx$$ with $\|f_0\|_2 = 1$ or $\psi(\xi) = \int_\Omega | \sum_j \xi_j f_j | \, d\mu$ with $f_j \in L_1(\Omega,\mu)$ come from infinite dimensional cocycles and taking $f_0 = \chi_\Sigma$, we obtain highly irregular $\psi$'s. Other examples will be given in \cite{JMP00}. These $m_\xi = \psi(\xi)^{is}$ are just exotic forms of Stein's imaginary powers and $L_p$-boundedness for $1 < p < \infty$ is guaranteed. The main novelty from our method is that we may identify the endpoint estimates for $T_m$ associated to $m_\xi = \psi(\xi)^{is}$, so that $T_m: L_\infty(\R^n) \to \mathrm{BMO}_{\T_\psi}$. 

The paper essentially follows the order established in this Introduction. At the end of the paper, we review some classical multiplier theorems in $\R^n$ under our cocycle formulation. We illustrate our main results on noncommutative tori and the free group algebra. As an application, we also construct in Corollary \ref{EXQMS} new examples of Rieffel's quantum metric spaces for the compact dual of virtually abelian groups. We conclude with a geometric analysis of our results. We shall assume some background on von Neumann algebras, noncommutative $L_p$ and Hardy spaces, as well as some operator space terminology. Standard references on operator algebra are \cite{KR,Ta}. We refer to \cite[Section 1]{Pa1} for a brief review of the results from noncommutative integration needed for this paper. A more in depth discussion is given in Pisier/Xu's survey \cite{PX2}. The $p$-norms of row/column square functions and corresponding Hardy spaces appear in \cite{JLX,P2,PX1}. Two excellent books on operator space theory are Effros/Ruan and Pisier monographs \cite{ER,P3}. 

\vskip5pt

\noindent \textbf{Acknowledgement.} Over the last years, we have discussed our results with many colleagues. We thank the interesting comments and bibliographic references from A. Carbery, G. Garrig\'os, D. M\"uller, N. Ozawa, J. Peterson, \'E. Ricard, A. Seeger, A. Thom and J. Wright. We are also indebted to the referee for his comments, which led to a significantly more transparent presentation. Junge is partially supported by the NSF DMS-0901457 and DMS-1201886, Mei by the NSF DMS-1266042 and Parcet by the ERC Grant StG-256997-CZOSQP and MTM2010-16518. Junge and Parcet are also supported in part by ICMAT Severo Ochoa Grant SEV-2011-0087.

\section{{\bf BMO spaces}}

We begin by introducing BMO spaces on finite von Neumann algebras associated to a Markov semigroup of operators. We will relate this construction with the standard definition of BMO in the Euclidean spaces. In the context of Fourier multipliers, we will prove an $L_\infty \to \mathrm{BMO}$ form of de Leeuw's compactification theorem \cite{dL} and construct a normal extension of Fourier multipliers defined on finite von Neumann algebras.  

\subsection{Semigroup $\mathrm{BMO}$ spaces} \label{BMOPrelim}

Let us briefly review the theory developed in \cite{JM2} of BMO spaces constructed over Markov semigroups. A semigroup of operators $\S = (S_t)_{t \ge 0}$ acting on a noncommutative probability space $(\M,\tau)$ ---a finite von Neumann algebra $\M$ equipped with a normal finite faithful trace $\tau$--- is called a noncommutative Markov semigroup when:
\begin{itemize}
\item[i)] $S_t(\1_\M) = \1_\M$ for all $t \ge 0$,

\vskip2pt

\item[ii)] Each $S_t$ is weak-$*$ continuous and completely positive on $\M$,

\vskip2pt

\item[iii)] Each $S_t$ is self-adjoint, i.e. $\tau ((S_tf)^*g) = \tau (f^*(S_tg))$ for $f,g \in \M \cap L_2(\M)$,

\vskip2pt

\item[iv)] $S_tf \rightarrow f$ as $t \rightarrow 0^+$ in the weak-$*$ topology of $\M$.
\end{itemize}
These conditions are reminiscent of Stein's notion \cite{St}. They imply that $S_t$ is trace preserving and extends to a semigroup of complete contractions on $L_p(\M)$ for any $1 \le p \le \infty$. $\S$ admits an infinitesimal generator $$-A = \lim_{t \rightarrow 0} \frac{S_t-id_\M}{t} \ \quad \mbox{with} \quad \ S_t = \exp(-tA).$$ In the $L_2$ setting, $A$ is a densely defined positive unbounded operator. In general we write $A_p$ for the generator of the realization of $\S$ on $L_p(\M)$. Note that $\mathrm{ker}A_p$ is a complemented subspace of $L_p(\M)$. Let $E_p$ denote the corresponding projection and $J_p = id_{L_p(\M)} - E_p$. Consider the complemented spaces $$L_p^\circ(\M) = J_p(L_p(\M)) = \Big\{ f \in L_p(\M) \, \big| \ \lim_{t \to \infty} S_t f = 0 \Big\}.$$ The BMO space on $(\M,\tau)$ associated to $\S = (S_t)_{t \ge 0}$ is defined as $$\mathrm{BMO}_\S(\M) \, = \, \Big\{ f \in L_2^\circ (\M) \ \big| \ \max \big\{ \|f\|_{\mathrm{BMO}_\S^c}, \|f^*\|_{\mathrm{BMO}_\S^c} \big\} < \infty \Big\},$$ where the column BMO seminorm is given by $$\|f\|_{\mathrm{BMO}_\S^c} \, = \, \sup_{t > 0} \Big\| \Big( S_t (f^*f) - S_t(f)^*S_t(f) \Big)^\frac12 \Big\|_\M.$$ This definition makes sense since we know from the Kadison-Schwarz inequality that $S_t(f)^* S_t(f) \le S_t(f^*f)$. The null space of this seminorm is the fixed-point subspace $\mathrm{ker} A_\infty$ of our semigroup. Indeed, if $\|f\|_{\mathrm{BMO}_\T^c} = 0$ we know from \cite{Ch} that $f$ belongs to the multiplicative domain of $S_t$, so that $$\tau(gf) = \tau(S_{t/2}(gf)) = \tau(S_{t/2}(g) S_{t/2}(f)) = \tau(g S_t(f)).$$ This proves that $f$ is fixed by the semigroup. Reciprocally, $\mathrm{ker} A_\infty$ is a $*$-subalgebra of $\M$ by \cite{JX2}. Thus, the seminorm vanishes on $\mathrm{ker} A_\infty$. In particular, we obtain a norm after quotienting out $\mathrm{ker} A_\infty$, which justifies our definition of $\mathrm{BMO}_\S(\M)$ as a subspace of $L_2^\circ(\M) = J_2(L_2(\M))$. The definition of $\mathrm{BMO}_\S(\M)$ for semifinite pairs $(\M,\tau)$ is more subtle and it demands some terminology from the theory of Hilbert modules. We will avoid that degree of generality by specifying in Paragraph \ref{EuclideanBMO} the relation of these BMO spaces with the classical one in $\R^n$. 

\begin{remark}
\emph{We impose the operator space structure given by $$M_k(\mathrm{BMO}_\T(\M)) = \mathrm{BMO}_{\widehat{\T}_{(k)}}(M_k(\M)) \quad \mbox{with} \quad \widehat{S}_{(k),t} = S_t \otimes id_{M_k}.$$ If $S_{\otimes,t} = S_t \otimes id_{\mathcal{B}(\ell_2)}$ we say $f \in \mathrm{BMO}_{\S_\otimes}(\M \bar\otimes \mathcal{B}(\ell_2))$ when the norm of the $M_k$ truncations $(id \otimes \Pi_k) f (id \otimes \Pi_k)$ in $M_k(\mathrm{BMO}_\T(\M))$ is uniformly bounded in $k \ge 1$.}
\end{remark}

It will be essential for us to provide interpolation results between semigroup type $\mathrm{BMO}$ spaces and the corresponding noncommutative $L_p$ spaces. It is a hard problem to identify the minimal regularity on the semigroup $\T = (S_t)_{t \ge 0}$ which suffices for this purpose. The first substantial progress was announced in a preliminary version of \cite{JM}, where the gradient form $2 \Gamma(f_1,f_2) = A(f_1^*)f_2 + f_1^*A(f_2) - A(f_1^*f_2)$ with $A$ the infinitesimal generator of the semigroup, was a key 
tool in finding sufficient regularity conditions in terms of nice Markov dilations. Concretely, we will say that a Markov semigroup admitting a reversed Markov dilation with almost uniformly (a.u.) continuous path ---see \cite{JM} for precise definitions--- is regular. As it follows from \cite{JRS}, all the semigroups that we handle in this paper are regular. The following result will be crucial in what follows, we refer the reader to \cite{JM2} for the proof. 

\begin{theorem} \label{Interpolation}
If $\T = (S_t)_{t \ge 0}$ is regular on $(\M,\tau)$ \vskip-7pt $$\big[ \mathrm{BMO}_\T(\M), L_p^\circ(\M) \big]_{p/q} \, \simeq_{cb} \, L_q^\circ(\M) \quad \mbox{for all} \quad 1 \le p < q < \infty.$$ 
\end{theorem}

\subsection{Relation with Euclidean BMO}
\label{EuclideanBMO}

Let us briefly recall the construction of operator-valued BMO spaces from \cite{Me1,NPTV} and relate it with the semigroup formulas above. Given a noncommutative measure space $(\M,\tau)$, we will write $(\RR,\varphi)$ for the von Neumann algebra generated by essentially bounded functions $f: \R^n \to \M$ which comes equipped with the trace $\varphi(f) = \int_{\R^n} \tau(f(y)) \, dy$. In other words, we set $\RR = L_\infty(\R^n) \bar\otimes \M$. On the other hand, recall that given a measure space $(\Omega,\mu)$ the norm in $L_\infty(\M; L_2^c(\Omega,\mu))$ is $$\|f\|_{L_\infty(\M;L_2^c(\Omega,\mu))} \, = \, \Big\| \Big( \int_{\Omega} |f(\omega)|^2 d\mu(\omega) \Big)^\frac12 \Big\|_\M.$$ A more precise definition of these spaces can be found for instance in \cite{JLX}.  If we set $(\Omega, \mu) = (\R^n, \mu_n)$ equipped with the measure $d\mu_n(x) = (1 + |x|^{n+1})^{-1} dx$, we define the corresponding column BMO spaces as follows $$\mathrm{BMO}_\RR^c \, = \, \Big\{ f \in L_\infty \big( \M; L_2^c(\R^n,\mu_n) \big) \ \big| \ \sup_{Q \in \Q} \Big\| \Big( \meann_Q |f(x) - f_Q|^2 \, dx \Big)^\frac12 \Big\|_\M < \infty \Big\},$$ where $\Q$ is the set of Euclidean balls in $\R^n$ and $f_Q$ denotes the average of $f$ over $Q$. The measure $\mu_n$ is chosen so that $\mathrm{BMO}_\RR^c$ defined as above is complete. Also, we will use that $f_Q \in \M$ for all $f \in L_\infty(\M; L_2^c(\R^n,\mu_n))$, see \cite{Me1} for further details. The row BMO space is defined similarly, the norm being $\|f\|_{\mathrm{BMO}_\RR^r} = \|f^*\|_{\mathrm{BMO}_\RR^c}$. As it is customary, the space $\mathrm{BMO}_\RR$ is defined as the intersection of row and column BMO and it comes equipped with the maximum of both norms, see \cite{Me1} for further details. As above, we impose the operator space structure $M_k(\mathrm{BMO}_\RR) = \mathrm{BMO}_{M_k(\RR)}$.

Let us now explain the relation between these spaces and semigroup BMO type norms. Note that we have deliberately omitted the definition of the latter spaces for non-finite von Neumann algebras, although the row and column BMO seminorms for bounded elements may be defined as in the finite case. The heat semigroup on $\R^n$ has an integral representation $$S_t(f)(x) = \int_{\R^n} k_t(x,y)f(y)dy \quad \mbox{with} \quad k_t(x,y)=\frac1{(4\pi t)^\frac n2} \exp \Big( \frac{-|x-y|^2}{4t} \Big).$$ Let $\mathrm{B}_r(x)$ be the Euclidean ball of radius $r$ centered at $x$, and let
$$\mathsf{E}_{\mathrm{B}_r(x)} f \, = \, \frac1{|\mathrm{B}_r(x)|}\int_{\mathrm{B}_r(x)} f(y) dy \, = \, \frac {\Gamma(\frac n2+1)}{\pi^\frac n2r^n}\int_{\mathrm{B}_r(x)} f(y) dy.$$ Then $S_t$ is a `global' mean value operator as a convex combination of $\mathsf{E}_{\mathrm{B}_r(x)}$'s 
\begin{eqnarray*}
S_tf(x) & = & \frac1{(4\pi t)^{\frac{n}{2}}} \int_{\R^n}  \Big( \int_{\frac {|x-y|^2}{4t}}^\infty e^{-u} du \Big) f(y)dy \\
&=&\frac1{(4\pi t)^{\frac{n}{2}}} \int_{\R_+} \Big( \int_{\mathrm{B}_{\sqrt{4ut}}(x)} f(y)dy \Big) e^{-u} du \\
&=&\frac{1}{\Gamma(\frac n2+1)}\int_{\R_+} e^{-u} \, u^{\frac n2} \, \big( \mathsf{E}_{\mathrm{B}_{\sqrt{4ut}}(x)} f \big) \, du.
\end{eqnarray*}
Of course, the same identity holds for operator-valued functions $f: \R^n \to \M$.

\begin{lemma} \label{HeatBMO}
Given $k \ge 1$ and $f \in M_k(\RR)$ $$\|f\|_{M_k(\mathrm{BMO}_\RR^c)} \, \sim_{c_n} \, \sup_{t > 0} \Big\| \Big( S_{\otimes,t}^{(k)} |f|^2 - |S_{\otimes,t}^{(k)} f|^2 \Big)^\frac12 \Big\|_{M_k(\RR)},$$ where $S_{\otimes,t}^{(k)} = S_t \otimes id_{M_k(\M)}$ and the constant $c_n$ only depends on the dimension $n$. 
\end{lemma}

\dem It suffices to assume $k=1$. The identity $$\Big( S_{\otimes,t}|f|^2 - |S_{\otimes,t}f|^2 \Big)(x) \, = \, \Big( S_{\otimes,t} \big| f - S_{\otimes,t}f(x) \big|^2 \Big)(x)$$ with $S_{\otimes,t} = S_t \otimes id_\M$ follows from our expression of $S_t$ as a global mean. Thus 
\begin{eqnarray*}
\lefteqn{\Big( S_{\otimes,t}|f|^2 - |S_{\otimes,t}f|^2 \Big)(x)} \\ & = & \frac{1}{\Gamma(\frac n2+1)}\int_{\R_+} e^{-u} \, u^{\frac n2} \, \mathsf{E}_{\mathrm{B}_{\sqrt{4ut}}(x)} \big| f - S_{\otimes,t}f(x) \big|^2 \, du \\ & = & \frac{1}{\Gamma(\frac n2+1)}\int_{\R_+} e^{-u} \, u^{\frac n2} \, \mathsf{E}_{\mathrm{B}_{\sqrt{4ut}}(x)} \Big| \frac{1}{\Gamma(\frac n2+1)} \int_{\R_+} e^{-v} v^{\frac n2} \big( f - \mathsf{E}_{\mathrm{B}_{\sqrt{4vt}}(x)} f \big) \, dv \Big|^2 \, du \\ & \le & \frac{1}{\Gamma(\frac n2+1)} \int_{\R_+} e^{-u} \, u^{\frac n2} \, \frac{1}{\Gamma(\frac n2+1)} \int_{\R_+} e^{-v} \, v^{\frac n2} \, \mathsf{E}_{\mathrm{B}_{\sqrt{4ut}}(x)} \big| f - \mathsf{E}_{\mathrm{B}_{\sqrt{4vt}}(x)} f \big|^2 \, dv \, du \\ & \le & \underbrace{\frac{1}{\Gamma(\frac n2+1)^2} \Big( \int_{\R_+} \int_{\R_+} e^{-u}\, u^{\frac n2} \, e^{-v} \, v^{\frac n2} \, \Big( \frac{u^2 + v^2}{uv} \Big)^{\frac{n}{2}} \, du \, dv \Big)}_{c_n} \, \|f\|_{\mathrm{BMO}_\RR^c}^2, 
\end{eqnarray*}  
where we have used the operator convexity of the function $| \ |^2$. For the upper estimate, fix a ball $\mathrm{B}_r(x)$ and fix $t$ so that $|\mathrm{B}_r(x)| = (4\pi t)^\frac n2$. With this choice it is very simple to observe that $\mathsf{E}_{\mathrm{B}_r(x)}h \le c_n' \, S_{\otimes,t}h(x)$ for any $h \ge 0$. Moreover, using again the operator convexity of $| \ |^2$ and Kadison-Schwarz inequality we deduce the following estimate  
\begin{eqnarray*}
\lefteqn{\mathsf{E}_{\mathrm{B}_r(x)} \big| f - \mathsf{E}_{\mathrm{B}_r(x)}f \big|^2} \\ & \lesssim & \mathsf{E}_{\mathrm{B}_r(x)} \Big( \big| f - S_{\otimes,t}f(x) \big|^2 + \big| \mathsf{E}_{\mathrm{B}_r(x)}f - S_{\otimes,t}f(x) \big|^2 \Big) \ \lesssim \ \mathsf{E}_{\mathrm{B}_r(x)} \big| f - S_{\otimes,t}f(x) \big|^2.
\end{eqnarray*}
Combining the estimates above we finally obtain 
\begin{eqnarray*}
\Big( \mathsf{E}_{\mathrm{B}_r(x)} \big| f - \mathsf{E}_{\mathrm{B}_r(x)}f \big|^2 \Big)^{\frac12} & \lesssim & \Big( \mathsf{E}_{\mathrm{B}_r(x)} \big| f-S_{\otimes,t}f(x) \big|^2 \Big)^{\frac12} \\ & \le & c_n' \Big( S_{\otimes,t} \big| f - S_{\otimes,t}f(x) \big|^2 \Big)^{\frac12} (x) \\ & \le & c_n' \, \sup_{t> 0} \Big\| \Big( S_{\otimes,t}|f|^2 - |S_{\otimes,t}f|^2 \Big)^{\frac12} \Big\|_\RR
\end{eqnarray*}
and the assertion follows taking norms in $\M$ and suprema in the balls $\mathrm{B}_r(x)$. \fin

\begin{remark}
\emph{When $\RR = L_\infty(\R^n)$ and $\M = \C$, we will write $\mathrm{BMO}_{\R^n}$ for $\mathrm{BMO}_\RR$.}
\end{remark}

\subsection{An $L_\infty \to \mathrm{BMO}$ de Leeuw theorem}

Let us now present a variation of de Leeuw's compactification theorem \cite{dL}. Instead of focusing on $L_p$-boundedness, we will be concerned instead with $L_\infty \to \mathrm{BMO}$ type estimates. Recall that we write $\Rd^n$ for $\R^n$ equipped with the discrete topology. Since $\R^n$ and $\Rd^n$ coincide as sets, we let $\Xi$ denote it. According to Pontryagin duality, the continuous characters on $\R^n$ and its compactification are both indexed by $\Xi$. Write $\chi_\xi^{\null}$ and $\chi'_\xi$ for the continuous characters on $\R^n$ and its Bohr compactification respectively. According to the construction of the Bohr compactification, we find a universal inclusion map $\Psi: \R^n \to \hatRd^n$ with dense image and such that $$\chi'_\xi(\Psi(x)) = \chi_\xi^{\null}(x)$$ for all $(\xi, x) \in \Xi \times \R^n$. The key point here is that the $L_\infty$ norms coincide on trigonometric polynomials. More concretely, let $\Lambda$ be a finite subset of $\Xi$ and consider the polynomials $f = \sum_{\xi \in \Lambda} a_\xi \chi_\xi$ and $f' = \sum_{\xi \in \Lambda} a_\xi \chi_\xi'$. Then continuity of $f, f'$ and density of $\Psi(\R^n)$ give $$\|f'\|_{L_\infty(\hatRd^n)} \ = \ \sup_{x \in \R^n} |f' \circ \Psi(x)| \ = \ \sup_{x \in \R^n} |f(x)| \ = \ \|f\|_{L_\infty(\R^n)}.$$ On the other hand, let $\T = (S_t)_{t \ge 0}$ and $\T' = (S_t')_{t \ge 0}$ denote respectively the heat semigroups on $\R^n$ and its Bohr compactification. Let us also consider a symbol $\widetilde{m}: \Xi \to \C$ which yields a Fourier multiplier $$\widehat{T_{\widetilde{m}}f} (\xi) \, = \, \widetilde{m}(\xi) \widehat{f}(\xi)$$ in $\R^n$ and $\hatRd^n$, considering of course the corresponding Fourier transform for each case. Then, since the algebra of trigonometric polynomials is preserved in both cases by the heat semigroup and the Fourier multiplier, we may use the above isometry together with Lemma \ref{HeatBMO} to obtain 
\begin{eqnarray*}
\|T_{\widetilde{m}}f'\|_{\mathrm{BMO}_{\T'}(\hatRd^n)} \!\! & = & \!\! \sup_{t > 0} \Big\| \Big( S_{t}'|T_{\widetilde{m}}f'|^2 - |S_t' T_{\widetilde{m}}f'|^2 \Big)^{\frac12} \Big\|_{L_\infty(\hatRd^n)} \\ \!\! & = & \!\! \sup_{t > 0} \Big\| \Big( S_{t}|T_{\widetilde{m}}f|^2 - |S_{t} T_{\widetilde{m}}f|^2 \Big)^{\frac12} \Big\|_{L_\infty(\R^n)} \sim_{c_n} \|T_{\widetilde{m}} f\|_{\mathrm{BMO}_{\R^n}^{\null}}.
\end{eqnarray*}
Let $\A_{\R^n}$ and $\A_{\hatRd^n}$ be the algebras of trigonometric polynomials in $\R^n$ and $\hatRd^n$.

\begin{lemma} \label{deLeeuwBMO}
We have $$\big\| T_{\widetilde{m}}: \A_{\R^n} \to \mathrm{BMO}_{\R^n} \big\| \, \sim_{c_n} \, \big\| T_{\widetilde{m}}: \A_{\hatRd^n} \to \mathrm{BMO}_{\T'}(\hatRd^n) \big\|.$$ Moreover, the same holds for cb-norms under the given operator space structures.
\end{lemma}

\dem \hskip-2pt Even \hskip-2pt in the operator space level it follows from the considerations above. \hskip-3pt \fin

\subsection{Normal extension of Fourier multipliers}

Lemma \ref{deLeeuwBMO} has been formulated for simplicity in the (weak-$*$ dense) algebra of trigonometric polynomials. Other auxiliary results below will be also formulated in weak-$*$ dense algebras. This is possible since we can always find a unique weak-$*$ continuous extension of Fourier multipliers which are defined over trigonometric polynomials. Before proving such statement, we need a duality result. Let $\psi: \G \to \R_+$ be any length on some discrete group $\G$ and denote by $\S_\psi = (S_{\psi,t})_{t \ge 0}$ the associated Markov semigroup $\lambda(g) \mapsto \exp(-t \psi(g)) \lambda(g)$ on $(\V,\tau_\G)$. Let $H_1(\S_\psi)$ denote the closure of $L_2^\circ(\widehat{\mathbb{G}})$ with respect to $$\|f\|_{H_1(\S_\psi)} = \sup_{\|h\|_{\mathrm{BMO}_{\S_\psi}} \le 1} |\tau_\G(fh^*)|$$ with the operator space structure $$\|f\|_{M_k(H_1(\S_\psi))}  \, = \, \|f\|_{\mathcal{CB}(\mathrm{BMO}_{\S_\psi},M_k)}.$$ Then, we shall prove in the Appendix that the map $$f \in \mathrm{BMO}_{\S_\psi} \mapsto \phi_f \in H_1(\S_\psi)^*$$ with $\phi_f$ densely defined by $\phi_f(h) = \tau_\G(f^*h)$ for $h \in L_2^\circ(\widehat{\mathbb{G}})$ is a complete isomorphism.

\begin{lemma} \label{NormalExt}
Let $\psi: \G \to \R_+$ be any length on some discrete group $\G$ and denote by $\S_\psi = (S_{\psi,t})_{t \ge 0}$ the associated Markov semigroup $\lambda(g) \mapsto \exp(-t \psi(g)) \lambda(g)$ on $(\V,\tau_\G)$. Let $\A_\G$ be the algebra of trigonometric polynomials in $\V$ and assume that $T_m: \A_\G \to \mathrm{BMO}_{\S_\psi}$ is a bounded Fourier multiplier for some bounded symbol $m: \G \to \C$. Then, $T_m$ extends uniquely to a normal $($i.e. weak-$*$ continuous$)$ bounded operator $$\widetilde{T}_m: \V \to \mathrm{BMO}_{\S_\psi}.$$ Moreover, if $T_m: \A_\G \stackrel{cb}{\longrightarrow} \mathrm{BMO}_{\S_\psi}$ its normal extension is also completely bounded. 
\end{lemma}

\dem Let $T_m^*$ be the adjoint of $T_m$ as an operator on $L_2(\widehat{\mathbb{G}})$. It suffices to show that $$\|T_m^*f\|_{L_1(\widehat{\mathbb{G}})} \, \lesssim \, \|f\|_{H_1(\S_\psi)}$$ for all $f \in L_2^\circ(\widehat{\mathbb{G}})$. Indeed, by the density of  $L_2^\circ(\widehat{\mathbb{G}})$ in $H_1(\S_\psi)$ this would imply that $T_m^*: H_1(\S_\psi) \to L_1(\widehat{\mathbb{G}})$. Taking adjoints and applying the duality theorem proved in the Appendix, we see that $$\widetilde{T}_m = T_m^{**}: \V \to \mathrm{BMO}_{\S_\psi}$$ is the weak-$*$ continuous extension of $T_m$ we were looking for. The uniqueness trivially follows from the weak-$*$ density of $\A_\G$. In order to prove the inequality at the beginning of this proof, we note that $T_m^*f \in L_2^\circ(\widehat{\mathbb{G}}) \subset L_1^\circ(\widehat{\mathbb{G}})$ for all $f \in L_2^\circ(\widehat{\mathbb{G}})$ since $m: \G \to \C$ is bounded. Moreover, by Kaplansky density theorem the unit ball of $\A_\G$ is weak-$*$ dense in the unit ball of $\V$. Therefore, we obtain the following inequality $$\|T_m^*f\|_{L_1(\widehat{\mathbb{G}})} \, = \sup_{\begin{subarray}{c} a \in \A_\G \\ \|a\|_{\V} \le 1 \end{subarray}} \big| \langle a, T_m^*f \rangle \big| \, \le \, \big\| T_m: \A_\G \to \mathrm{BMO}_{\S_\psi} \big\| \, \|f\|_{H_1(\S_\psi)}$$ for any $f \in L_2^\circ(\widehat{\mathbb{G}})$. A similar argument applies for completely bounded norms. \fin

\section{{\bf Crossed product extensions}}
\label{SectSemi}

We now study the $L_\infty \to \mathrm{BMO}$ boundedness of semidirect product extensions of semicommutative CZO's. Our results are of independent interest, regarded as a first step through Calder\'on-Zygmund theory on fully noncommutative von Neumann algebras, see \cite{JMP} for related results. 

\subsection{Crossed products}

Given a discrete group $\G$ with left regular representation $\lambda: \G \to \mathcal{B}(\ell_2(\G))$, let $\V$ denote its group von Neumann algebra and $L_p(\widehat{\mathbb{G}})$ the associated noncommutative $L_p$ space, as defined in the Introduction. Note that for $\G$ abelian we get the $L_p$ space on the dual group equipped with its normalized Haar measure. We will sometimes keep the terminology $\V$ for $p=\infty$. Given another noncommutative probability space $(\M, \tau)$ with $\M \subset \mathcal{B(\H)}$, assume that there exists a trace preserving action $\alpha: \G \to \mathrm{Aut}(\M)$. Define the crossed product algebra $\M \rtimes_\alpha \G$ as the weak operator closure in $\M \bar\otimes \mathcal{B}(\ell_2(\G))$ of the $*$-algebra generated by $\rho(\M)$ and $\Lambda(\G)$, where the $*$-representations $\rho: \M \to \M \bar\otimes \mathcal{B}(\ell_2(\G))$ and $\Lambda: \G \to \M \bar\otimes \mathcal{B}(\ell_2(\G))$ are given by $$\rho(f) = \sum_{h \in \G} \alpha_{h^{-1}}(f) \otimes e_{h,h} \quad \mbox{and} \quad \Lambda(g) = \sum_{h \in \G} \1_\M \otimes e_{gh,h},$$ with $e_{g,h}$ the matrix units for $\mathcal{B}(\ell_2(\G))$. A generic element of $\M \rtimes_\alpha \G$ can be formally written as $\sum_{g \in \G} f_g \rtimes_\alpha \lambda(g)$ (understanding the infinite sum as a limit in the weak-$*$ topology), where each $f_g \in \M$. With this convention, we may embed the crossed product algebra $\M \rtimes_\alpha \G$ into $\M \bar\otimes \mathcal{B}(\ell_2(\G))$ by means of the map $j = \rho \rtimes \Lambda$. Indeed, we have
\begin{eqnarray*}
\lefteqn{\hskip-10pt j \Big( \summ_g f_g \rtimes_\alpha \lambda(g) \Big) \ = \ \summ_g \rho(f_g) \Lambda(g)} \\
& = & \summ_g \Big( \summ_{h,h'} \big( \alpha_{h^{-1}}(f_g) \otimes e_{h,h} \big) \big( \1_\M \otimes  e_{gh',h'} \big) \Big)
\\ & = & \summ_{g} \Big( \summ_h \alpha_{h^{-1}}(f_g) \otimes e_{h,g^{-1}h} \Big)  \
= \  \summ_{g} \Big( \summ_h \alpha_{(gh)^{-1}}(f_g) \otimes e_{gh,h} \Big).
\end{eqnarray*}
Similar computations lead to
\begin{itemize}
\item $(f \rtimes_\alpha \lambda(g))^* = \alpha_{g^{-1}}(f^*) \rtimes_\alpha \lambda(g^{-1})$,

\item $(f \rtimes_\alpha \lambda(g)) (f' \rtimes_\alpha \lambda(g')) = f \alpha_g(f') \rtimes_\alpha \lambda(gg')$,

\item $\tau \rtimes \tau_\G(f \rtimes_\alpha \lambda(g)) = \tau \otimes \tau_\G(f \otimes \lambda(g)) = \delta_{g=e} \tau(f)$.
\end{itemize}
Since $\alpha$ will be fixed, we relax the terminology and write $\sum_g f_g \lambda(g) \in \M \rtimes \G$ for generic elements in the crossed product, instead of $\sum_g f_g \rtimes_\alpha \lambda(g)$. We say that a semigroup $\T = (S_t)_{t \ge 0}$ on $(\M,\tau)$ is $\G$-equivariant if $$\alpha_g S_t = S_t \alpha_g \quad \mbox{for} \quad (t,g) \in \R_+ \times \G.$$

If $\S$ is a $\G$-equivariant Markov semigroup on $\M$, let $\T_\rtimes = (S_t \rtimes id_\G)_{t \ge 0}$ and $\T_\otimes = (S_t \otimes id_{\mathcal{B}(\ell_2(\G))})_{t \ge 0}$ denote the cross/tensor product amplification of our semigroup on $\M \rtimes \G$ and $\M \bar\otimes \mathcal{B}(\ell_2(\G))$ respectively. Note that the Markovianity of $\S_\rtimes$ relies on the $\G$-equivariance of $\S$. Let $\A_{\M \rtimes \G}$ be the subalgebra of finite sums $\sum_g f_g \lambda(g) \in \M \rtimes \G$. In the following result, we compute the $\mathrm{BMO}_{\S_\rtimes}$-norm of elements in $\A_{\M \rtimes \G}$ in terms of the map $j = \rho \rtimes \Lambda$ defined above. Recall that we write $S_{\otimes,t}^{(k)}$ for the $k$-th matrix amplification $S_{\otimes,t} \otimes id_{M_k}$.

\begin{lemma} \label{BMOisometryBMO}
If $\T$ is $\G$-equivariant, $k \ge 1$ and $f \in M_k(\A_{\M \rtimes \G})$ $$\|f\|_{M_k(\mathrm{BMO}^c_{\T_\rtimes}(\M \rtimes \G))} \, = \, \sup_{t > 0} \Big\| \Big( S_{\otimes,t}^{(k)}|j(f)|^2 - |S_{\otimes,t}^{(k)}j(f)|^2 \Big)^{\frac12} \Big\|_{M_k(\M \bar\otimes \mathcal{B}(\ell_2(\G)))}.$$ 
\end{lemma}

\dem If $f \in \A_{\M \rtimes \G}$, we have $$S_{\rtimes, t} |f|^2 - |S_{\rtimes, t} f|^2 = \sum_{g,h \in \G}^{\null} \alpha_{g^{-1}} \big( S_t(f_g^*f_h) - S_t(f_g^*)S_t(f_h) \big) \lambda(g^{-1}h)$$ since $\T$ is $\G$-equivariant. Then, simple algebraic calculations yield $$j \Big( S_{\rtimes, t} |f|^2 - |S_{\rtimes, t} f|^2 \Big) \, = \, S_{\otimes, t}|j(f)|^2 - |S_{\otimes, t} j(f)|^2.$$ Moreover, the exact same calculations show that the identity above remains valid after matrix amplification. Therefore, the assertion follows from the fact that the map $j = \rho \rtimes \Lambda$ is a complete isometry $\M \rtimes \G \to \M \bar\otimes \mathcal{B}(\ell_2(\G))$. \fin

\subsection{Equivariant extension}
\label{SectEquivCZO}

Given $(\M_j, \tau_j)_{j=1,2}$ noncommutative probability spaces, assume $\G \curvearrowright \M_j$ by trace preserving actions $\alpha_j$. Let $\S_2 = (S_{2,t})_{t \ge 0}$ denote a $\G$-equivariant Markov semigroup on $(\M_2,\tau_2)$. Given a weak-$*$ dense $\alpha_1$-invariant subalgebra $\A_1$ of $\M_1$, consider a map $$T: \A_1 \to \mathrm{BMO}_{\T_2}(\M_2)$$ satisfying $T(\A_1) \subset \M_2$. Then we say that $T$ is $\G$-equivariant if $\alpha_{2,g} T f = T \alpha_{1,g} f$ for all $g \in \G$ and all $f \in \A_1$. Note that our assumption $T(\mathcal{A}_1) \subset \M_2$ is used in the definition. Let $\S_{2 \rtimes} = (S_{2,t} \rtimes id_\G)_{t \ge 0}$ denote the cross product (Markovian) extension of $\S_2$ acting on $\M_2 \rtimes \G$. Let $$\A_{\A_1 \rtimes \G} = \Big\{ f = \sum_{g \in \Lambda}^{\null} f_g \lambda(g) \ \big| \ f_g \in \A_1, \ \Lambda \subset \G \ \mbox{finite} \Big\} \subset \M_1 \rtimes \G.$$ Note that $\A_{\A_1 \rtimes \G}$ is a subalgebra by $\alpha_1$-invariance of $\A_1$. Moreover, Kaplansky density theorem implies that the unit ball of $\A_{\A_1 \rtimes \G}$ ---the space of trigonometric polynomials in $\A_1 \rtimes \G$--- is weak-$*$ dense in the unit ball of $\M_1 \rtimes \G$. We now construct a completely bounded extension $T \rtimes id_\G: \A_{\A_1 \rtimes \G} \to \mathrm{BMO}_{\S_{2\rtimes}}(\M_2 \rtimes \G)$.       

\begin{lemma} \label{LemEquivCZO}
If $T: \A_1 \stackrel{cb}{\longrightarrow} \mathrm{BMO}_{\T_2}(\M_2)$ is $\G$-equivariant $$T \rtimes id_\G: \A_{\A_1 \rtimes \G} \to \mathrm{BMO}_{\T_{2\rtimes}}(\M_2 \rtimes \G) \quad \mbox{is also completely bounded}.$$ The same conclusion holds when $T$ is bounded and the $\M_j$'s are commutative.
\end{lemma}

\dem Let $$j_1: \M_1 \rtimes \G \to \M_1 \bar\otimes \mathcal{B}(\ell_2(\G)),$$ $$j_2: \M_2 \rtimes \G \to \M_2 \bar\otimes \mathcal{B}(\ell_2(\G)),$$ stand for the natural injections. Given $g \in \G$ and $f_g \in \A_1$, we have $$j_2 \big( T \rtimes id_\G (f_g \lambda(g)) \big) = \summ_h T(\alpha_{1,(gh)^{-1}}(f_{g})) \otimes e_{gh,h} = T \otimes id_{\mathcal{B}(\ell_2(\G))} \big( j_1 (f_g \lambda(g)) \big)$$ by $\G$-equivariance of $T$. By linearity, this identity trivially extends to arbitrary elements in $\A_{\A_1 \rtimes \G}$. On the other hand, the $\alpha_1$-invariance of $\A_1$ implies that $j_1(\A_{\A_1 \rtimes \G}) \subset \A_1 \otimes \mathcal{B}(\ell_2(\G))$. In particular, if we set
\begin{eqnarray*}
T_{\rtimes\G}^{(k)} & = & T \rtimes id_\G \otimes id_{M_k}, \\ T_{\otimes\G}^{(k)} & = & T \otimes id_{\mathcal{B}(\ell_2(\G))} \otimes id_{M_k},
\end{eqnarray*}
we obtain from Lemma \ref{BMOisometryBMO} the following estimate for $f \in M_k(\A_{\A_1 \rtimes \G})$
\begin{eqnarray*}
\lefteqn{\big\| T_{\rtimes\G}^{(k)}f \big\|_{M_k(\mathrm{BMO}_{\S_{2\rtimes}}^c)}} \\ & = & \sup_{t > 0} \Big\| \Big( S_{2\otimes,t}^{(k)} \big| j_2^{(k)}(T_{\rtimes\G}^{(k)}f) \big|^2 - \big| S_{2\otimes,t}^{(k)} j_2^{(k)}(T_{\rtimes\G}^{(k)}f) \big|^2 \Big)^{\frac12} \Big\|_{M_k(\M_2 \bar\otimes \mathcal{B}(\ell_2(\G)))} \\ & = & \sup_{t > 0} \Big\| \Big( S_{2\otimes,t}^{(k)} \big| T_{\otimes\G}^{(k)}(j_1^{(k)}f) \big|^2 - \big| S_{2\otimes,t}^{(k)} T_{\otimes\G}^{(k)}(j_1^{(k)}f) \big|^2 \Big)^{\frac12} \Big\|_{M_k(\M_2 \bar\otimes \mathcal{B}(\ell_2(\G)))}. 
\end{eqnarray*}
Now, since the semigroup $S_{2\otimes,t}^{(k)}$ is given by $S_t \otimes id_{\mathcal{B}(\ell_2(\G))} \otimes id_{M_k}$ and $$T_{\otimes \G}^{(k)}(j_1^{(k)}(f)) \in T(\A_1) \bar\otimes \mathcal{B}(\ell_2(\G)) \otimes M_k \subset \M_2 \bar\otimes \mathcal{B}(\ell_2(\G)) \otimes M_k,$$ we may use the operator space structure of $\mathrm{BMO}_{\S_2}^c(\M_2)$ to deduce that $$\big\| T_{\rtimes\G}^{(k)}f \big\|_{M_k(\mathrm{BMO}_{\S_{2\rtimes}}^c)} \, \le \, \big\| T: \A_1 \to \mathrm{BMO}_{\S_2}^c \big\|_{cb} \, \|f\|_{M_k(\A_{\A_1 \rtimes \G})}.$$ Here we have also used that $j_1$ is a complete isometry $\M_1 \rtimes \G \to \M_1 \bar\otimes \mathcal{B}(\ell_2(\G))$. Note that our argument above is also valid for row BMO norms. Hence, this completes the proof of the first assertion in the statement. For the second one, we may assume that $(\M_j, \tau_j)$ is of the form $L_\infty(\Omega_j, \mu_j)$. According to the first part of the statement it suffices to see that any bounded map $$T: L_\infty(\Omega_1) \to \mathrm{BMO}_{\T_2}(\Omega_2)$$ is indeed cb bounded. Our argument is row/column symmetric and we just consider the column case. Given a matrix-valued function $f = (f_{ij}): \Omega_1 \to M_k$, let us write $T^{(k)}$ for $T \otimes id_{M_k}$. Then we have 
\begin{eqnarray*}
\lefteqn{\hskip-15pt \big\| T^{(k)}f \big\|_{M_k(\mathrm{BMO}_{\T_2}^c(\Omega_2))}} \\ & = & \sup_{t > 0} \Big\| \Big( S_{2,t}^{(k)} \big| T^{(k)}f \big|^2 - \big| S_{2,t}^{(k)}T^{(k)}f \big|^2 \Big)^\frac12 \Big\|_{L_\infty(\Omega_2; M_k)} \\ & = & \sup_{\begin{subarray}{c} t > 0 \\ \|\xi\|_{\ell_2(k)} \le 1 \end{subarray}} \esssup_{w \in \Omega_2} \Big\langle \xi, \Big[ S_{2,t}^{(k)} \big| T^{(k)}f \big|^2 - \big| S_{2,t}^{(k)}T^{(k)}f \big|^2 \Big](w) \xi \Big\rangle_{\ell_2(k)}^\frac12.
\end{eqnarray*}
This quantity may be approximated by fixing some $(t_0, \xi_0) \in \R_+ \times \ell_2(k)$ and the essential supremum can also be replaced by an average over a set $\Sigma \subset \Omega_2$ of finite positive measure where the given function is close enough to its esssup. In other words, taking $$\mu_\Sigma(A) = \frac{\mu_2(A \cap \Sigma)}{\mu_2(\Sigma)}$$ we may approximate as follows 
\begin{eqnarray*}
\lefteqn{\hskip-15pt \big\| T^{(k)}f \big\|_{M_k(\mathrm{BMO}_{\T_2}^c(\Omega_2))}} \\ & \sim & \Big\langle \xi_0, \Big( \int_{\Omega_2} \Big[ S_{2,t_0}^{(k)} \big| T^{(k)}f \big|^2 - \big| S_{2,t_0}^{(k)}T^{(k)}f \big|^2 \Big](w) d\mu_\Sigma(w) \Big) \xi_0 \Big\rangle_{\ell_2(k)}^\frac12.
\end{eqnarray*}
It is useful to write the expression above in terms of Hilbert modules. Set $$\Big( z_{ij} \Big)_{i,j \le k} = \Big( z(f_{ij}) \Big)_{i,j \le k} = \Big( Tf_{ij} \otimes \1_{\Omega_2} - \1_{\Omega_2} \otimes S_{2,t_0} Tf_{ij} \Big)_{i,j \le k}$$ and consider the bracket $\langle f_1 \otimes f_2, f_1' \otimes f_2' \rangle = f_2^* S_{2,t_0}(f_1^*f_1') f_2'$. Then we find the identity $$S_{2,t_0}^{(k)} \big| T^{(k)}f \big|^2 - \big| S_{2,t_0}^{(k)} T^{(k)} f \big|^2 = \sum_{i,j=1}^k \Big( \sum_{r=1}^k \langle z_{ri},z_{rj} \rangle \Big) \otimes e_{ij}.$$ Consider the Hilbert module $\Lambda_{t_0}(\Omega_2,\mu_\Sigma)$ defined as the closure of the algebraic tensor product $L_\infty(\Omega_2, \mu_\Sigma) \otimes L_\infty(\Omega_2, \mu_\Sigma)$ in the topology determined by $\xi_\alpha \to \xi$ iff $$\int_{\Omega_2} \big\langle \xi - \xi_\alpha, \xi - \xi_\alpha \big\rangle \phi \, d\mu_\Sigma \to 0 \quad \mbox{for all} \quad \phi \in L_1(\Omega_2, \mu_\Sigma).$$ Combining the GNS construction for completely positive unital maps \cite{Lance} with the Hilbert space associated to W$^*$-modules studied by  Paschke [52], we see that there exists a map $u: \Lambda_{t_0}(\Omega_2, \mu_\Sigma) \to L_\infty(\Omega_2,\mu_\Sigma; \ell_2^c(\mathcal{I}))$ ---for some index set $\mathcal{I}$--- which satisfies $$\langle f_1,f_2 \rangle = u(f_1)^*u(f_2).$$ Note that the definition of the target space of $u$ can be found in the Introduction before the statement of Theorem C. Now, if we consider the map $v = u \circ z$ we deduce 
\begin{eqnarray*}
\sum_{i,j=1}^k \Big( \sum_{r=1}^k \langle z_{ri},z_{rj} \rangle \Big) \otimes e_{ij} & = & \sum_{i,j=1}^k \Big( \sum_{r=1}^k u(z_{ri})^*u(z_{rj}) \Big) \otimes e_{ij} \\ & = & \sum_{i,j=1}^k \Big( \sum_{r=1}^k v(f_{ri})^*v(f_{rj}) \Big) \otimes e_{ij} \ = \ \big| v^{(k)} f \big|^2 
\end{eqnarray*}
Combining our results so far, we get
\begin{eqnarray*}
\lefteqn{\hskip-15pt \big\| T^{(k)}f \big\|_{M_k(\mathrm{BMO}_{\T_2}^c(\Omega_2))}} \\ & \sim & \Big\langle \xi_0, \Big( \int_{\Omega_2} \big| v^{(k)}(f)(w) \big|^2 \, d\mu_\Sigma(w) \Big) \xi_0 \Big\rangle_{\ell_2(k)}^\frac12 \\ & \le & \Big\| \Big( \int_{\Omega_2} \big| v^{(k)}(f)(w) \big|^2 \, d\mu_\Sigma(w) \Big)^{\frac12} \Big\|_{M_k} \ = \ \big\| v^{(k)}(f) \big\|_{M_k(L_2^c(\Omega_2,\mu_\Sigma; \ell_2(\mathcal{I})))}.
\end{eqnarray*}
Thus, we have reduced the problem to show that $v: L_\infty(\Omega_1) \to L_2^c(\Omega_2,\mu_\Sigma; \ell_2(\mathcal{I}))$ is a completely bounded map. Assume for a moment that $v$ is bounded when regarded as a Banach space operator. By the little Grothendieck inequality \cite{PBook,PBull}, this means that $v$ is absolutely $2$-summing so that we can find a factorization $v = w \circ j_\xi$ where $j_\xi$ is the map $f \in L_\infty(\Omega_1, \mu_1) \mapsto f\xi \in L_2(\Omega_1, \mu_1)$ for some $\xi$ with $\int_{\Omega_1} |\xi|^2 \, d\mu_1 = 1$ and we have $$\|w\| \le \frac{2}{\sqrt{\pi}} \|v\|.$$ This immediately gives that $$\|v\|_{cb} \le \big\| w: L_2^c(\Omega_1,\mu_1) \to L_2^c(\Omega_2, \mu_\Sigma; \H) \big\|_{cb} \big\| j_\xi: L_\infty(\Omega_1,\mu_1) \to L_2^c(\Omega_1,\mu_1) \big\|_{cb}$$ and yields $\|v\|_{cb} \le \|w\| \le 2/\sqrt{\pi} \|v\|$, because $j_\xi$ is a complete contraction and column Hilbert spaces are homogeneous operator spaces, see e.g. \cite{P3}. Thus, we just need to compute the Banach space norm of $v$. However, applying again the properties of the right module map $u$, we obtain for $f \in L_\infty(\Omega_1)$ and $$z(f) = Tf \otimes \1_{\Omega_2} - \1_{\Omega_2} \otimes S_{2,t_0} Tf$$ the following estimate
\begin{eqnarray*}
\|v(f)\|_{L_2(\Omega_2,\mu_\Sigma; \ell_2(\mathcal{I}))} & \le & \|v(f)\|_{L_\infty(\Omega_2,\mu_\Sigma; \ell_2(\mathcal{I}))} \\ [4pt] & = & \|u(z(f))^*u(z(f))\|_{L_\infty(\Omega_2,\mu_\Sigma)}^{\frac12} \\ & = & \big\| S_{2,t_0}|Tf|^2 - |S_{2,t_0}Tf|^2 \big\|_{L_\infty(\Omega_2,\mu_\Sigma)}^{\frac12} \ \le \
\|Tf\|_{\mathrm{BMO}_{\T_2}^c(\Omega_2)} 
\end{eqnarray*}
and $\|v: L_\infty(\Omega_1,\mu_1) \to L_2(\Omega_2,\mu_\Sigma; \ell_2(\mathcal{I}))\| \le \|T: L_\infty(\Omega_1,\mu_1) \to \mathrm{BMO}_{\T_2}^c(\Omega_2)\|$. \fin
  
\subsection{Semicommutative CZO's}
\label{SubSemi}

Given noncommutative measure spaces $(\M_j,\tau_j)$ for $j=1,2$, we will write $(\RR_j,\varphi_j)$ to denote the von Neumann algebra generated by essentially bounded functions $f: \R^n \to \M_j$ which comes equipped with the trace $\varphi_j(f) = \int_{\R^n} \tau_j(f(y)) \, dy$. In other words, we have $$\RR_j = L_\infty(\R^n) \bar\otimes \M_j.$$ Let us now consider Calder\'on-Zygmund operators in these algebras associated to operator-valued kernels. Our construction is standard, we refer to \cite{D,RRT} for further details. Let us we write $L_0(\M_j)$ for the $*$-algebra of $\tau_j$-measurable operators affiliated with $\M_j$. Consider kernels $k: \R^{2n} \setminus \Delta \to \mathcal{L}(L_0(\M_1), L_0(\M_2))$ which are defined away from the diagonal $\Delta = \{(x,x) \, | \, x \in \R^n\}$ and take values in linear maps from $\tau_1$-measurable to $\tau_2$-measurable operators. We further assume that  $k(x,y)$ is bounded when regarded as a map $\M_1 \to \M_2$ for $x \neq y$, and that $$\|k(x,y)\|_{\mathcal{B}(\M_1, \M_2)} \lesssim \frac{1}{|x-y|^n}.$$  We will consider linear operators associated to this class of kernels. By that, we only mean that $T$ is well-defined on certain (nice) space of functions ---typically $L_2(\RR_1)$, but we will need to use a different space in Lemma \ref{LEMSCCZ} i) below--- and that for $f$ in such space, we have $$Tf(x) = \int_{\R^n} k(x,y) \, (f(y)) \, dy \quad \mbox{for} \quad x \notin \mathrm{supp}_{\R^n} f.$$ Additional information for these operators require to impose some smoothness in our kernels. In our case, this is given in Lemma \ref{LEMSCCZ} ii). We refer to \cite[Chapter V]{D} for a detailed discussion on to what extend these conditions determine the operator $T$ completely. Note that $L_p(\RR_j) = L_p(\R^n; L_p(\M_j))$, but endpoint estimates do not follow from the vector-valued theory, see \cite{Pa1} for further explanations.

Recalling the definition of $\mathrm{BMO}_{\RR_2}$ from Paragraph \ref{EuclideanBMO}, our goal is to analyze conditions for the $L_\infty(\RR_1) \to \mathrm{BMO}_{\RR_2}$ boundedness of CZO's. To that end, we shall also be working with the spaces $L_\infty(\M_j; L_2^c(\R^n))$, whose rigorous definition can be found in \cite{JLX}. 

\begin{lemma} \label{LEMSCCZ}
We have $T: L_\infty(\RR_1) \to \mathrm{BMO}_{\RR_2}^c$ provided

\begin{itemize}
\item[i)] $L_2$-column condition,
$$\Big\| \Big( \int_{\R^n} |Tf(x)|^2 \, dx \Big)^\frac12 \Big\|_{\M_2} \ \lesssim \ \Big\| \Big( \int_{\R^n} |f(x)|^2 \, dx \Big)^\frac12 \Big\|_{\M_1}.$$ That is, $T: L_\infty(\M_1;L_2^c(\R^n)) \to L_\infty(\M_2;L_2^c(\R^n))$ is bounded.

\vskip3pt

\item[ii)] Smoothness condition for the kernel, $$\esssup_{x_1, x_2} \int_{|x_1 - y| > 2 |x_1 - x_2|} \big\| k(x_1,y) - k(x_2,y) \big\|_{\mathcal{B}(\M_1, \M_2)} \, dy \, < \, \infty.$$
\end{itemize}
\end{lemma}

\dem Given $g \in \mathrm{BMO}_\RR^c$, we first observe that $$\|g\|_{\mathrm{BMO}_{\RR}^c} \ \sim_2 \ \sup_{Q \in \Q} \, \inf_{a_Q \in \M} \Big\| \Big( \frac{1}{|Q|} \int_Q \big| g(x) - a_Q \big|^2 dx \Big)^\frac12 \Big\|_\M.$$ Indeed, we have $$\Big\| \Big( \mean_Q \big| g(x) - g_Q \big|^2 dx \Big)^\frac12 \Big\|_\M \le \ \Big\| \Big( \mean_Q \big| g(x) - a_Q \big|^2 dx \Big)^\frac12 \Big\|_\M + \big\| a_Q - g_Q \big\|_\M$$ and Kadison-Schwarz inequality for the ucp map $u(g) = g_Q \otimes \1_{\R^n}$ gives rise to $u(h)^*u(h) \le u(h^*h)$ for the function $h(x) = g(x) - a_Q$. Therefore, we obtain the estimate $$\big\| a_Q - g_Q \big\|_\M = \| u(h) \|_\RR \le \big\| u(h^*h) \big\|_\RR^\frac12 = \Big\| \Big( \mean_Q \big| g(x) - a_Q \big|^2 dx \Big)^\frac12 \Big\|_\M.$$ This proves the upper estimate, the lower estimate is clear. Now, given $f \in L_\infty(\RR_1)$ and a ball $Q$, we set as usual $f_1 = f \chi_{5Q}$ and $f_2 = f - f_1$ where $5Q$ denotes the ball concentric to $Q$ whose radius is 5 times the radius of $Q$. Then we pick $a_Q = - \hskip-9pt \int_Q Tf_2(x) \, dx \in \M_2$. It therefore suffices to prove $$\mathrm{A} + \mathrm{B} = \Big\| \Big( \mean_Q \big| Tf_1(x) \big|^2 dx \Big)^\frac12 \Big\|_{\M_2} + \ \Big\| \Big( \mean_Q \big| Tf_2(x) - a_Q \big|^2 dx \Big)^\frac12 \Big\|_{\M_2} \ \lesssim \ \|f\|_{L_\infty(\RR_1)}.$$ According to the $L_2$-column condition i), we find $$\mathrm{A} \ \le \ \frac{1}{\sqrt{|Q|}} \Big\| \Big( \int_{5Q} \big| f(x) \big|^2 dx \Big)^\frac12 \Big\|_{\M_1} \ \le \ 5^n \, \|f\|_{L_\infty(\RR_1)}.$$ On the other hand, since $\mathrm{supp}_{\R^n} f_2 \cap Q = \emptyset$ we have for $x \in Q$ $$Tf_2(x) - a_Q = \mean_Q \big( Tf_2(x) - Tf_2(z) \big) \, dz = \mean_Q \int_{\R^n} \big( k(x,y) - k(z,y) \big) (f_2(y)) \, dy \, dz.$$ Using again the Kadison-Schwarz inequality, this gives rise to
\begin{eqnarray*}
\null \hskip15pt \mathrm{B} & = & \Big\| \Big( \mean_Q \big| Tf_2(x) - a_Q \big|^2 dx \Big)^\frac12 \Big\|_{\M_2} \\ & \le & \Big( \mean_Q \mean_Q \Big\| \int_{\R^n} \big( k(x,y) - k(z,y) \big) (f_2(y)) \, dy \Big\|_{\M_2}^2 \, dz \, dx \Big)^\frac12 \\ & \le & \Big( \mean_Q \mean_Q \Big[ \int_{\R^n \setminus 5Q} \big\| k(x,y) - k(z,y) \big\|_{\mathcal{B}(\M_1, \M_2)} \, dy \Big]^2 \, dz \, dx \Big)^\frac12 \|f\|_{L_\infty(\RR_1)} \\ & \le & \Big( \esssup_{x,z \in \R^n} \int_{|x-y| > 2 |x-z|} \big\| k(x,y) - k(z,y) \big\|_{\mathcal{B}(\M_1, \M_2)} \, dy \Big) \, \|f\|_{L_\infty(\RR_1)}. 
\end{eqnarray*}
Then we apply to the last expression our smoothness condition for the kernel. \fin

\vskip5pt

\begin{remark} \label{Op-Vld model}
\emph{The $L_2$-boundedness condition i) reduces to the classical one when $\M_1 = \M_2$ and $k(x,y)$ acts on $f(y)$ by left multiplication. Indeed, if we assume that $T$ is bounded on $L_2(\RR)$ and use $\M \subset \mathcal{B}(\H)$ for $\H = L_2(\M)$, the following inequality holds for $f \in L_\infty(\M;L_2^c(\R^n))$
\begin{eqnarray*}
\lefteqn{\hskip-10pt \Big\| \Big( \int_{\R^n} |Tf(y)|^2 \, dy
\Big)^\frac12 \Big\|_\M} \\ & = & \sup_{\|h\| \le 1} \Big(
\int_{\R^n} \big\langle h, |Tf(y)|^2 h \big\rangle_\H \, dy
\Big)^{\frac12} \\ [7pt] & = & \sup_{\|h\| \le 1}
\big\| T(f \, (\1_{\R^n} \otimes h)) \big\|_{L_2(\RR)} \ \lesssim \
\sup_{\|h\| \le 1} \big\| f \, (\1_{\R^n} \otimes h) \big\|_{L_2(\RR)} \\ & = &
\sup_{\|h\| \le 1} \Big( \int_{\R^n}
\big\langle h, |f(y)|^2 h \big\rangle_\H \, dy \Big)^{\frac12}
\ = \ \Big\| \Big( \int_{\R^n}
|f(y)|^2 \, dy \Big)^{\frac12} \Big\|_\M.
\end{eqnarray*}
This is false for other operator kernels and our $L_2$-condition seems the natural one.}
\end{remark}

\begin{remark} \label{RemBMOCZO}
\emph{Since $M_k(\mathrm{BMO}_\RR) = \mathrm{BMO}_{M_k(\RR)}$, it suffices to replace $\M$ by $M_k(\M)$ everywhere, amplify all the involved maps by tensorizing with $id_{M_k}$ and require that the hypotheses hold with $k$-independent constants to deduce complete boundedness in the statement above.}
\end{remark}

\subsection{Nonequivariant extension}
\label{NEQCZOSect}

Set $\Lambda^\dag f = (\Lambda f^*)^*$ for any mapping $\Lambda$. In the nonequivariant setting, the arguments are not row/column symmetric because the map $(T \rtimes id_\G)^\dag$ is not similar to $T \rtimes id_\G$. This will be specially relevant in the $L_2$-boundedness conditions that we obtain. Indeed, we have for finite sums 
\begin{eqnarray*}
(T \rtimes id_\G)^\dag \big( \summ_g f_g \lambda(g) \big) \!\!\! &
= & \!\!\! \Big[ (T \rtimes id_\G) \big( \summ_g
\alpha_{g^{-1}}(f_g^*) \lambda(g^{-1}) \big) \Big]^* \\ \!\!\! & =
& \!\!\! \summ_g \alpha_g \big( T(\alpha_{g^{-1}}(f_g^*))^* \big)
\lambda(g) =  \summ_g \alpha_g T^\dag \alpha_{g^{-1}} (f_g)
\lambda(g).
\end{eqnarray*}
Thus, $(T \rtimes id_\G)^\dag$ is a map of the form $\summ_g f_g \lambda(g)
\mapsto \summ_g T_g(f_g) \lambda(g)$ and recalling the
embedding $j = \rho \rtimes \Lambda: \M \rtimes \G \to \M \bar\otimes
\mathcal{B}(\ell_2(\G))$, we see after the change of variables $g \mapsto gh^{-1}$ that
\begin{eqnarray*}
j \big( \summ_g T_g(f_g) \lambda(g) \big) & = & \summ_{g,h}
\alpha_{g^{-1}}(T_{gh^{-1}}(f_{gh^{-1}})) \otimes  e_{g,h} \\ & =
& \Big( \alpha_{g^{-1}} T_{gh^{-1}} \alpha_g \Big) \bullet j \big(
\summ_g f_g \lambda(g) \big) \ = \ \Phi \Big( j \big( \summ_g f_g
\lambda(g) \big) \Big),
\end{eqnarray*}
where the $\bullet$ stands for the generalized Schur product of matrices, in the sense that $\alpha_{g^{-1}} T_{gh^{-1}} \alpha_g$ only acts in the $(g,h)$-th entry for each $g,h \in \G$. We will use the following terminology
\begin{itemize}
\item[$\bullet$] $\chi_\xi'$ are the characters of $\hatRd^n$ for all $\xi \in \Rd^n$,

\vskip3pt

\item[$\bullet$] Trigonometric polynomials are denoted by $f' = \sum_\xi a_\xi \chi_\xi'$,

\item[$\bullet$] $\A_{\hatRd^n}$ is the algebra of trigonometric polynomials in $\hatRd^n$,

\vskip1pt

\item[$\bullet$] $\A_{\hatRd^n \rtimes \G}$ is the algebra of trigonometric polynomials in $\A_{\hatRd^n} \rtimes \G$ $$\hskip2pt \A_{\hatRd^n \rtimes \G} \, = \Big\{ \sum_{g \in \Lambda} f_g' \lambda(g) \ \big| \ f_g' \in \A_{\hatRd^n}, \ \Lambda \ \mbox{finite} \Big\} \, \subset \, L_\infty(\hatRd^n) \rtimes \G,$$

\item[$\bullet$] $\S_\rtimes' = (S'_{\rtimes,t})_{t \ge 0}$ with $S_{\rtimes,t}' = S_t' \rtimes_id_\G$ for the heat semigroup $\S'$ on $\hatRd^n$.
\end{itemize} 
Let $\beta: \G \to O(n)$ be an orthogonal action. The action $\alpha: \G \curvearrowright L_\infty(\hatRd^n)$ determined by $\alpha_g(\chi_\xi') = \chi'_{\beta_g(\xi)}$ is clearly trace preserving. In the result below we provide crossed product extensions of Fourier multipliers defined in the Bohr compactification of $\R^n$ under this class of actions.

\begin{lemma} \label{LemNonequiv}
Let $\alpha: \G \curvearrowright L_\infty(\hatRd^n)$ be a trace preserving action implemented by $\beta$ as pointed above. Let us consider a family of Fourier multipliers $T_{\widetilde{m}_g}$ indexed by $g \in \G$ so that $$\widehat{T_{\widetilde{m}_g}f'}(\xi) = \widetilde{m}_g(\xi) \widehat{f'}(\xi) \quad \mbox{for} \quad f' \in \A_{\hatRd^n}.$$ Then, the cross product extension $$\A_{\hatRd^n \rtimes \G} \ni \summ_g f_g'
\lambda(g) \mapsto \summ_g T_{\widetilde{m}_g}(f_g') \lambda(g) \in \mathrm{BMO}_{\S_\rtimes'}^c$$ is a completely bounded map provided the following conditions hold in $\R^n$
\begin{itemize}
\item[i)] $L_2$-column condition, $$\hskip25pt \Big\| \Big( \int_{\R^n} \big| (T_{\widetilde{m}_{gh^{-1}}}) \bullet \rho (x) \big|^2 \, dx \Big)^\frac12 \Big\|_{\mathcal{B}(\ell_2(\G))} \lesssim_{cb} \Big\| \Big( \int_{\R^n} |\rho(x)|^2 \, dx \Big)^\frac12 \Big\|_{\mathcal{B}(\ell_2(\G))},$$ i.e. $\rho \in L_\infty(\mathcal{B}(\ell_2(\G));L_2^c(\R^n)) \stackrel{cb}{\longmapsto} (T_{\widetilde{m}_{gh^{-1}}}) \bullet \rho \in L_\infty(\mathcal{B}(\ell_2(\G));L_2^c(\R^n))$.

\vskip3pt

\item[ii)] Smoothness condition for the kernel, $$\hskip25pt
\esssup_{x_1, x_2} \int_{|x_1-y| > 2|x_1-x_2|} \big\| K(x_1,y) - K(x_2,y) \big\|_{\mathcal{CB}(\mathcal{B}(\ell_2(\G)))} \, dy \, < \, \infty,$$ where $\displaystyle K(x,y) = \sum_{g,h} \widehat{\widetilde{m}}_{gh^{-1}}(\beta_g(x-y)) \otimes e_{g,h}$ acts by the Schur product $\bullet$.
\end{itemize}
\end{lemma}

\dem According to Lemma \ref{BMOisometryBMO} and since $$j \Big( \underbrace{\summ_g T_{\widetilde{m}_g}(f_g') \lambda(g)}_{Tf'} \Big) \ = \ \Phi \Big( j \big( \underbrace{\summ_g f_g' \lambda(g)}_{f'} \big) \Big),$$ it suffices to show that we have for such $f' \in \A_{\hatRd^n \rtimes \G}$ $$\Big\| \Big( S_{\otimes,t}' \big| \Phi j (f') \big|^2 - \big| S_{\otimes,t}' \Phi j (f') \big|^2 \Big)^{\frac12} \Big\|_{L_\infty(\hatRd^n) \bar\otimes \mathcal{B}(\ell_2(\G))} \lesssim \|j(f')\|_{L_\infty(\hatRd^n) \bar\otimes \mathcal{B}(\ell_2(\G))}$$ with constants independent of $t > 0$. Of course, we also need to prove the $k$-th matrix amplification of this inequality for each $k \ge 1$, but the argument in that case is identical and we shall omit it here. In order to prove such inequality our first observation is that we have $$\Big( S_{\otimes,t}' \big| \Phi j (f') \big|^2 - \big| S_{\otimes,t}' \Phi j (f') \big|^2 \Big) \circ \Psi = S_{\otimes,t} \big| \Phi j (f) \big|^2 - \big| S_{\otimes,t} \Phi j (f) \big|^2,$$ where $\Psi: \R^n \to \hatRd^n$ is the universal inclusion map used in the proof of Lemma \ref{deLeeuwBMO}; $f = \sum_g f_g \lambda(g)$ is the trigonometric polynomial in $L_\infty(\R^n) \rtimes \G$ which results after replacing $\chi_\xi'$'s by $\chi_\xi^{\null}$'s in each $f_g'$; and $\S = (S_t)_{t \ge 0}$ is the heat semigroup acting on $\R^n$. We keep the same terminology for $j$ and $\Phi$ understanding that now they act on $\R^n$ instead of its Bohr compactification. Its proof boils down to the identity $\chi_\xi = \chi_\xi' \circ \Psi$ and the fact that all the involved operators respect the structure of trigonometric polynomials. We leave the (easy) details to the reader. Once we have such identity ---together with $j(f') \circ \Psi = j(f)$---  we may argue as in the proof of Lemma \ref{deLeeuwBMO} to show that the desired inequality above is equivalent to the following $$\Big\| \Big( S_{\otimes,t} \big| \Phi j (f) \big|^2 - \big| S_{\otimes,t} \Phi j (f) \big|^2 \Big)^{\frac12} \Big\|_{L_\infty(\R^n) \bar\otimes \mathcal{B}(\ell_2(\G))} \lesssim \|j(f)\|_{L_\infty(\R^n) \bar\otimes \mathcal{B}(\ell_2(\G))}.$$ According to Lemma \ref{HeatBMO}, it suffices to prove that $$\Phi: \A_{\R^n} \bar\otimes \mathcal{B}(\ell_2(\G)) \to \mathrm{BMO}_{\RR}^c$$ is a bounded (in fact completely bounded) map, where $\A_{\R^n}$ is the algebra of trigonometric polynomials in $\R^n$ and $\RR = L_\infty(\R^n) \bar\otimes \mathcal{B}(\ell_2(\G))$. Recall that for $a = \sum_\xi a_\xi \chi_\xi$ in $\A_{\R^n}$, we have $\alpha_{g^{-1}}(a)(x) = \sum_\xi a_\xi \chi_{\beta_{g^{-1}}(\xi)}(x) = a(\beta_g(x))$. Thus, considering $$\rho = \sum_{g,h} a_{g,h} \otimes e_{g,h} \in \A_{\R^n} \bar\otimes \mathcal{B}(\ell_2(\G)),$$ we find that for $x \notin \mbox{supp}_{\R^n} \rho$
\begin{eqnarray*}
\Phi(\rho)(x) \!\!\!\!\! & = & \!\!\!\!\! \summ_{g,h} \alpha_{g^{-1}} \int_{\R^n} \widehat{\widetilde{m}}_{gh^{-1}}(x-y)
(a_{g,h}(\beta_{g^{-1}}(y))) \hskip1pt dy \otimes e_{g,h} \\
\!\!\!\!\! & = & \!\!\!\!\! \summ_{g,h} \int_{\R^n}
\widehat{\widetilde{m}}_{gh^{-1}}(\beta_g(x-y)) (a_{g,h}(y)) \hskip1pt dy \otimes
e_{g,h} \! = \! \int_{\R^n} K(x,y) (\rho(y)) \hskip1pt dy.
\end{eqnarray*}
Therefore, we may regard $\Phi$ as a semicommutative CZO and apply Lemma \ref{LEMSCCZ} (together with Remark \ref{RemBMOCZO}) to conclude. Let us then check the assumptions in Lemma \ref{LEMSCCZ}. First, we note that the $L_2$-column condition in Lemma \ref{LEMSCCZ} means that the map $\Phi: L_\infty(\mathcal{B}(\ell_2(\G)); L_2^c(\R^n)) \to L_\infty(\mathcal{B}(\ell_2(\G)); L_2^c(\R^n))$ is cb. However, we have $$\Phi(\rho) = (\alpha_{g^{-1}} T_{\widetilde{m}_{gh^{-1}}} \alpha_g) \bullet \rho = (\alpha_{g^{-1}}) \bullet  (T_{\widetilde{m}_{gh^{-1}}}) \bullet (\alpha_g) \bullet \rho.$$ Using that $\beta$ is measure preserving, the generalized Schur product $$\summ_{g,h} a_{g,h} \otimes e_{g,h} \mapsto \summ_{g,h} \alpha_g(a_{g,h}) \otimes e_{g,h}$$ is a complete isometry on $L_\infty(\mathcal{B}(\ell_2(\G));L_2^c(\R^n))$. Hence, the $L_2$-column condition in Lemma \ref{LEMSCCZ} for $\Phi$ reduces to the complete boundedness condition in the statement. The smoothness condition matches exactly that of Lemma \ref{LEMSCCZ}. \fin

\begin{lemma} \label{LemNonequiv2}
Let $T_{\widetilde{m}}: \A_{\hatRd^n} \to \mathrm{BMO}_{\T'}$ given by $$\widehat{T_{\widetilde{m}}f'}(\xi) = \widetilde{m}(\xi) \widehat{f'}(\xi) \quad \mbox{for} \quad f' \in \A_{\hatRd^n}.$$ Then, $T_{\widetilde{m}} \rtimes id_\G: L_\infty(\hatRd^n) \rtimes \G \to \mathrm{BMO}_{\T_{\rtimes}'}^c$ is completely bounded provided
\begin{itemize}
\item[i)] $\widetilde{m}: \Rd^n \to \C$ bounded,

\vskip3pt

\item[ii)] $\displaystyle \esssup_{x_1, x_2} \int_{|x_1-y| > 2 |x_1-x_2|} \sup_{g \in \G} \big| \widehat{\widetilde{m}}(\beta_g (x_1-y)) - \widehat{\widetilde{m}}(\beta_g (x_2-y)) \big| \, dy \ < \ \infty$.
\end{itemize}
\end{lemma}

\dem This is a particular case of Lemma \ref{LemNonequiv} with $T_{\widetilde{m}_g} = T_{\widetilde{m}}$ for all $g \in \G$. The $L_2$-column condition clearly reduces to the $L_2$-boundedness of $T_{\widetilde{m}}$ ---since $L_2$ comes equipped with the OH operator space structure--- which in turn is equivalent to the boundedness of $\widetilde{m}$. On the other hand, the kernel in Lemma \ref{LemNonequiv} has the following form now $K(x,y) = \sum_{g,h} \widehat{\widetilde{m}}(\beta_g(x-y)) \otimes e_{g,h}$. Hence, we deduce
\begin{eqnarray*}
K(x,y) (f(y)) & = & \Big( \widehat{\widetilde{m}}(\beta_g (x-y)) \Big) \bullet \Big( f_{g,h}(y) \Big) \\ & = & \Big[ \summ_g \widehat{\widetilde{m}}(\beta_g (x-y)) \otimes e_{gg} \Big] \Big[ \summ_{g,h} f_{g,h}(y) \otimes e_{g,h} \Big].
\end{eqnarray*}
In particular, regarding $K(x,y)$ as a left multiplication map (not a Schur multiplier) it is a diagonal matrix in $\mathcal{B}(\ell_2(\G))$. Therefore, we may easily rewrite the H\"ormander smoothness condition for the kernel in Lemma \ref{LemNonequiv} as in the statement. \fin

\section{{\bf H\"ormander-Mihlin multipliers}}
\label{SectHM}

We now study Fourier multipliers over the group von Neumann algebra of an arbitrary discrete group $\G$.  In the language of quantum groups, these algebras are regarded as the compact dual of $\G$. Our main result is a cocycle form of H\"ormander-Mihlin multiplier theorem in this setting.

\subsection{Length functions and cocycles}
\label{SubLFC}

A left cocycle associated to a discrete group $\G$ is a triple $(\H, \alpha, b)$ formed by a Hilbert space $\H$, an isometric action $\alpha: \G \to \mathrm{Aut}(\H)$ and a map $b: \G \to \H$ so that $\alpha_g(b(h)) = b(gh) - b(g)$. A right cocycle satisfies the relation $\alpha_g(b(h)) = b(hg^{-1}) - b(g^{-1})$ instead. In this paper a (cocycle) \emph{length function} $\psi: \G \to \R_+$ is any symmetric conditionally negative function vanishing at the identity of $\G$, as defined in the Introduction. The fact that any length function takes values in $\R_+$ is easily justified. Any cocycle $(\H, \alpha,b)$ can be identified with an affine representation $$g \in \G \mapsto \big( b(g) \rtimes \alpha_g \big) \in \mathrm{Aff}(\H).$$ In what follows, we only consider cocycles with values in real Hilbert spaces. Note that $\mathrm{Aut}(\H)$ is the orthogonal group on $\H$ and $\mathrm{Aff}(\H) \simeq \H \rtimes O(\H)$. Any cocycle $(\H, \alpha,b)$ gives rise to an associated length function $\psi_b(g) = \langle b(g), b(g) \rangle_\H$, as it can be checked by the reader. Reciprocally, any length function $\psi$ gives rise to a left and a right cocycle. This is a standard application of the ideas around Schoenberg's theorem \cite{Sc}, which claims that $\psi: \G \to \R_+$ is a length function if and only if the mappings $S_{\psi,t}(\lambda(g)) = \exp(-t\psi(g)) \lambda(g)$ extend to a semigroup of unital completely positive maps on $\V$. Let us collect these well-known results.

\begin{lemma} \label{Lemmapsi}
If $\psi: \G \to \R_+$ is a length function$\, :$
\begin{itemize}
\item[i)] The forms
\begin{eqnarray*}
K_\psi^1(g,h) & = & \frac{\psi(g) + \psi(h) - \psi(g^{-1}h)}{2}, \\
K_\psi^2(g,h) & = & \frac{\psi(g) + \psi(h) - \psi(gh^{-1})}{2},
\end{eqnarray*}
define positive matrices on $\G \times \G$ and lead to $$\Big\langle \summ_g a_g \delta_g, \summ_h b_h
\delta_h \Big\rangle_{\psi,j} = \summ_{g,h} a_g K_\psi^j(g,h) b_h$$ on the group algebra $\R[\G]$ of finitely supported real functions on $\G$.

\vskip5pt

\item[ii)] Let $\H_\psi^j$ be the Hilbert space completion of $$(\R[\G]/N_\psi^j, \langle \cdot \hskip1pt, \cdot \rangle_{\psi,j}) \quad \mbox{with} \quad N_\psi^j = \mbox{null space of} \ \langle \cdot \hskip1pt , \cdot \rangle_{\psi,j}.$$ If we consider the mapping $b_\psi^j: g \in \G \mapsto \delta_g + N_\psi^j \in \H_\psi^j$
\begin{eqnarray*}
\alpha_{\psi,g}^1 \Big( \sum_{h \in \G} a_h b_\psi^1(h) \Big) & = & \sum_{h \in \G} a_h \big( b_\psi^1(gh) - b_\psi^1(g) \big), \\ \alpha_{\psi,g}^2 \Big( \sum_{h \in \G} a_h b_\psi^2(h) \Big) & = & \sum_{h \in \G} a_h \big( b_\psi^2(hg^{-1}) - b_\psi^2(g^{-1}) \big),
\end{eqnarray*}
determine isometric actions $\alpha_\psi^j: \G \to \mathrm{Aut}(\H_\psi^j)$ of $\G$ on $\H_\psi^j$.

\vskip5pt

\item[iii)] Imposing the discrete topology on $\H_\psi^j$, the semidirect product $\G_\psi^j = \H_\psi^j \rtimes \G$ becomes a discrete group and we find the following group homomorphisms
\begin{eqnarray*}
\pi_\psi^1: g \in \G & \mapsto & b_\psi^1(g) \rtimes g \in \G_\psi^1, \\ \pi_\psi^2: g \in \G & \mapsto & b_\psi^2(g^{-1}) \rtimes g \in \G_\psi^2.
\end{eqnarray*}
\end{itemize}
\end{lemma}

The previous lemma allows us to introduce two pseudo-metrics on our discrete group $\G$ in terms of the length function $\psi$. Indeed, a short calculation leads to the crucial identities $$\psi(g^{-1}h) \, = \, \big\langle b_\psi^1(g) - b_\psi^1(h), b_\psi^1(g) - b_\psi^1(h) \big\rangle_{\psi,1} \, = \, \big\| b_\psi^1(g) - b_\psi^1(h) \big\|_{\H_\psi^1}^2,$$ $$\psi(gh^{-1}) \, = \, \big\langle b_\psi^2(g) - b_\psi^2(h), b_\psi^2(g) - b_\psi^2(h) \big\rangle_{\psi,2} \, = \, \big\| b_\psi^2(g) - b_\psi^2(h) \big\|_{\H_\psi^2}^2.$$ In particular, $$\mathrm{dist}_1(g,h) = \sqrt{\psi(g^{-1}h)} = \|b_\psi^1(g) - b_\psi^1(h) \|_{\H_\psi^1}$$ defines a pseudo-metric on $\G$, which becomes a metric when the cocycle map is injective. Similarly, we may work with $\mathrm{dist}_2(g,h) = \sqrt{\psi(gh^{-1})}$. The following elementary observation will be crucial for what follows.

\begin{lemma} \label{LemLeftRight}
Let $(\H_1, \alpha_1, b_1)$ and $(\H_2, \alpha_2, b_2)$ be a left and a right cocycle on $\G$. Assume that the associated length functions $\psi_{b_1}$ and $\psi_{b_2}$ coincide, then the following map defines an isometric isomorphism $$\Lambda_{12}: b_1(g) \in \H_1 \mapsto b_2(g^{-1}) \in \H_2.$$ In particular, given a length function $\psi$ we have $\H_{\psi}^1 \simeq \H_\psi^2$ via $b_\psi^1(g) \mapsto b_\psi^2(g^{-1})$.
\end{lemma}

\dem By polarization, we see that $$\big\langle b_1(g), b_1(h) \big\rangle_{\H_1} = \frac{1}{2} \Big( \|b_1(g)\|_{\H_1}^2 + \|b_1(h)\|_{\H_1}^2- \|b_1(g) - b_1(h)\|_{\H_1}^2 \Big).$$ Since $b_1(g) - b_1(h) = \alpha_{1,h}(b_1(h^{-1}g))$, we obtain
\begin{eqnarray*}
\big\langle b_1(g), b_1(h) \big\rangle_{\H_1} & = & \frac{\psi_{b_1}(g) + \psi_{b_1}(h) - \psi_{b_1}(g^{-1}h)}{2} \\ & = & \frac{\psi_{b_2}(g) + \psi_{b_2}(h) - \psi_{b_2}(g^{-1}h)}{2} \ = \ \big\langle b_2(g^{-1}), b_2(h^{-1}) \big\rangle_{\H_2}.
\end{eqnarray*}
The last identity uses polarization and $b_2(g^{-1}) - b_2(h^{-1}) = \alpha_{2,h}(b_2(g^{-1}h))$.
\fin

\subsection{Smooth Fourier multipliers}
\label{SubHM}

We are now ready to prove our extension of H\"ormander/Mihlin's sufficient condition for Fourier multipliers to arbitrary discrete groups. The ideas leading to the next result probably go back to H\"ormander, but we could not find the specific statement given below in the literature. We provide a proof based on Stein's approach to these questions in his book \cite{St2}.

\begin{lemma} \label{LemStein}
Let $k_{\widetilde{m}}$ be a tempered distribution on $\R^n$ which coincides with a locally integrable function on $\R^n \setminus \{0\}$. Let $\widetilde{m}$ stand for its Fourier transform $\widetilde{m} = \widehat{k}_{\widetilde{m}}$ and consider an orthogonal action $\beta: \G \to O(n)$. Then we obtain$\, :$
\begin{itemize}
\item[i)] If $| \partial_\xi^\gamma \widetilde{m}(\xi) | \, \le \, c_n |\xi|^{-|\gamma|}$ for all $|\gamma| \le n+2$ $$\esssup_{x \in \R^n} \int_{|y| > 2 |x|} \, \sup_{g \in \G} \big| k_{\widetilde{m}}(\beta_g y - \beta_g x) - k_{\widetilde{m}}(\beta_g y) \big| \, dy \ < \ \infty.$$

\item[ii)] If $| \partial_\xi^\gamma \widetilde{m}(\xi) | \, \le \, c_n |\xi|^{-|\gamma|}$ for all $|\gamma| \le [\frac{n}{2}]+1$, the operator 
\begin{eqnarray*}
\lefteqn{\hskip-10pt \Phi: \sum_{g,h \in \G} f_{gh} \otimes e_{gh} \in L_2(\RR)} \\ \hskip25pt & \mapsto & \sum_{g,h \in \G} \int k_{\widetilde m}(\beta_g x - \beta_g y) f_{gh}(y) \, dy \otimes e_{gh} \in L_2(\RR)
\end{eqnarray*}
extends to a cb-map from $\RR$ to $\mathrm{BMO}_\RR^c$, where $\RR = L_\infty(\R^n) \bar\otimes \mathcal{B}(\ell_2(\G))$.
\end{itemize}
\end{lemma}

\dem For i), it suffices to show that $|\nabla k_{\widetilde{m}}(z)| \lesssim |z|^{-(n+1)}$. Let $\eta \in \mathcal{C}^\infty(\R^n)$ with $\chi_{\mathrm{B}_1(0)} \le \eta \le \chi_{\mathrm{B}_2(0)}$ and take $\delta(\xi) = \eta(\xi) - \eta(2\xi)$ so that $\sum_{j \in \Z} \delta(2^{-j}\xi) = 1$ for all $\xi \neq 0$. This gives rise to $\widetilde{m}(\xi) = \sum_j \widetilde{m}(\xi) \delta(2^{-j}\xi) = \sum_j \widetilde{m}_j(\xi)$ and we set $$k_{\widetilde{m}}^j(x) = \int_{\R^n} \widetilde{m}_j(\xi) e^{2 \pi i \langle x, \xi \rangle} \, d\xi.$$ We have $\sum_j k_{\widetilde{m}}^j \to k_{\widetilde{m}}$ as distributions, so that it suffices to estimate $\sum_{j} | \partial_x^\alpha k_{\widetilde{m}}^j(x) |$ for any $x
\neq 0$ and any multi-index $\alpha$ with $|\alpha| = 1$. Now we claim
that a) $\Rightarrow$ b) with
\begin{itemize}
\item[a)] $\displaystyle \big| \partial_\xi^\beta
\widetilde{m}(\xi) \big| \le c_M |\xi|^{- |\beta|}$ for all
multi-index $\beta$ s.t. $0 \le |\beta| \le M$.

\item[b)] $\displaystyle \big| \partial_x^\alpha
k_{\widetilde{m}}^j(x) \big| \le c_M |x|^{-M} 2^{j(n-M+1)}$ for
all multi-index $\alpha$ s.t. $|\alpha|=1$.
\end{itemize}
Let us first see how the assertion follows from the claim. Indeed,
we know from our hypotheses that a) holds for any $0 \le M \le
n+2$. If we apply our claim for $M=0$ on those $j$'s for which
$2^j \le |x|^{-1}$ and we apply it for $M=n+2$ on those $j$'s for
which $2^j > |x|^{-1}$, we find $$\sum_{j \in \Z}
\big|\partial_x^\alpha k_{\widetilde{m}}^j(x) \big| \ \lesssim \
\sum_{2^j \le |x|^{-1}} 2^{j(n+1)} + \frac{1}{|x|^{n+2}} \sum_{2^j
> |x|^{-1}} 2^{-j} \ \sim \ \frac{1}{|x|^{n+1}}.$$ To prove our
claim, we use the properties of the Fourier transform to get
$$(-2\pi i x)^\gamma \partial_x^\alpha k_{\widetilde{m}}^j(x) \, = \,
\int_{\R^n} \partial_\xi^\gamma \big[ (2 \pi i \xi)^\alpha
\widetilde{m}_j(\xi) \big] e^{2 \pi i \langle x, \xi \rangle} \,
d\xi.$$ On the other hand, using condition a) it is not difficult
to check that we have $$\Big| \partial_\xi^\gamma \big[ (2 \pi i
\xi)^\alpha \widetilde{m}_j(\xi) \big] \Big| \le \sum_{\gamma_1 +
\gamma_2 = \gamma} c_{\gamma_1 \gamma_2} \Big|
\partial_\xi^{\gamma_1} \big( (2\pi i \xi)^\alpha \big)
\partial_\xi^{\gamma_2} \widetilde{m}_j(\xi) \Big| \lesssim
|\xi|^{1-|\gamma|}.$$ Moreover, since $\widetilde{m}_j$ is
supported by an annulus of radius $\sim 2^j$, we conclude that
$$\Big| \int_{\R^n} \partial_\xi^\gamma \big[ (2 \pi i \xi)^\alpha
\widetilde{m}_j(\xi) \big] e^{2 \pi i \langle x, \xi \rangle} \,
d\xi \Big| \, \lesssim \, 2^{jn} 2^{j(1-|\gamma|)}.$$ Given $x
\in \R^n$ there exists a multi-index $\gamma$ such that $|\gamma|
= M$ and $|x^\gamma| \sim |x|^M$. Hence, taking such a multi-index $\gamma$ in
the identity above we deduce our claim. Let us now prove ii). If $f=\sum_g f_{gh}\ten e_{gh} \in L_2(\RR) \cap L_\infty(\RR)$, we have
\begin{eqnarray*}
\lefteqn{\hskip-20pt \Big\| \frac1{|Q|} \int_Q \big| \Phi(f)-(\Phi(f))_Q \big|^2 \Big\|_{\mathcal{B}(\ell_2(\G))}^\frac12} \\ & = & \sup_{\|\xi\|_{\ell_2(\G)}=1} \Big( \frac1{|Q|} \int_Q \big\| \Phi(f)\xi - (\Phi(f)\xi)_Q \big\|_{\ell_2(\G)}^2 \Big)^\frac12 \\ & = & \sup_{\|\xi\|_{\ell_2(\G)}=1} \Big( \frac1{|Q|} \int_Q \big\| \Phi(f \xi) - (\Phi(f\xi))_Q \big\|_{\ell_2(\G)}^2 \Big)^\frac12,
\end{eqnarray*}
where $f\xi = \sum_g(\sum_h f_{gh}\xi_h) \ten e_{ge}$ satisfies $\|f\xi\|_{L_\infty(\R^n; \ell_2(\G))} \le \|f\|_\RR$. The problem is then reduced to show that the restriction of $\Phi$ to column matrices extends to a bounded map $L_\infty(\R^n; \ell_2(\G)) \to \mathrm{BMO}(\R^n; \ell_2(\G))$. Let us decompose $f\xi$ in the usual way $f_{\xi,1} = f\xi \chi_{5Q}$ and $f_{\xi,2} = f\xi - f_{\xi,1}$. By the $L_2(\R^n; \ell_2(\G))$ boundedness of $\Phi$, we have
\begin{eqnarray*}
\lefteqn{\hskip-20pt \Big\| \frac1{|Q|} \int_Q \big| \Phi(f)-(\Phi(f))_Q \big|^2 \Big\|_{\mathcal{B}(\ell_2(\G))}^\frac12} \\ & \lesssim & \sup_{\|\xi\|_{\ell_2(\G)} = 1} \Big( \frac1{|Q|} \int_Q \|\Phi (f_{\xi,1})\|_{\ell_2(\G)}^2 \Big)^\frac12 \\ & + & \sup_{\|\xi\|_{\ell_2(\G)} = 1} \Big( \frac1{|Q|} \int_Q \big\| \Phi(f_{\xi,2}) - (\Phi(f_{\xi,2}))_Q \big\|^2_{\ell_2(\G)} \Big)^\frac12 \\ [7pt] & \lesssim & \|f\|_\RR \, + \, \sup_{\begin{subarray}{c} \|\xi\|_{\ell_2(\G)} = 1 \\ x,z \in Q \end{subarray}} \big\| \Phi(f_{\xi,2})(x) - \Phi(f_{\xi,2})(z) \big\|_{\ell_2(\G)}.
\end{eqnarray*}
The last term on the right hand side can be estimated by
\begin{eqnarray*}
\lefteqn{\sup_{\|\xi\|_2, \|\eta\|_2 = 1} \Big| \int_{|x-y| \ge 2|x-z|} \Big\langle \big( \eta_g [ k_{\widetilde{m}}(\beta_g x - \beta_g y) - k_{\widetilde{m}}(\beta_g z - \beta_g y) ] \big),f_{\xi,2}(y) \Big\rangle \, dy \Big|} \\ & \le & \sup_{\|\eta\|_2 = 1} \int_{|x-y| \ge 2|x-z|} \Big\| \big( \eta_g [ k_{\widetilde{m}}(\beta_g x - \beta_g y) - k_{\widetilde{m}}(\beta_g z - \beta_g y) ] \big) \Big\|_{\ell_2(\G)} \, dy \ \|f\|_\RR.
\end{eqnarray*}
Following a classical argument, it is easy to check that the last term in the inequality above is finite. In fact, arguing as for inequality (32) of  \cite[VI.4.4.2]{St2} ---see also our estimates for i)--- we may decompose $k_{\widetilde{m}} = \sum_j k_{\widetilde{m}}^j$ and conclude that 
\begin{eqnarray*}
\lefteqn{\hskip-20pt \int_{\R^n} |x|^{2M} \Big\| \sum_g \eta_g k_{\widetilde{m}}^j(\beta_g x) \delta_g \Big\|_{\ell_2(\G)}^2 \, dx} \\ & = & \summ_g |\eta_g|^2 \int_{\R^n} |\beta_{g^{-1}}x|^{2M}  |k_{\widetilde{m}}^j(x) |^2 \, dx \\ & = & \summ_g |\eta_g|^2 \int_{\R^n} |x|^{2M}  |k_{\widetilde{m}}^j(x) |^2 \, dx \ \lesssim \ 2^{j(n-2M)},
\end{eqnarray*}
for any $0 \le M \le [\frac{n}{2}]+1$ and any $j \in \Z$. The remaining part of the estimation is the same to that of \cite[VI.4.4.2 page 247]{St2}. Thus, taking the supremum over $Q$ we deduce the estimate for the norm of $T_{k_{\widetilde{m}}}$. The cb-norm is estimated similarly. This shows the $L_\infty \to \mathrm{BMO}$ boundedness for elements in $L_2 \cap L_\infty$. The extension to $L_\infty$ functions follows the classical argument, see e.g. \cite{Gr}. \fin

We are now ready for the main result of this paper. Let $\G$ be a discrete group and consider a bounded symbol $m: \G \to \C$. Then, the associated Fourier multiplier map is constructed as $$T_m: \sum_{g \in \G} \widehat{f}(g) \lambda(g) \mapsto \sum_{g \in \G}^{\null} m_g \widehat{f}(g) \lambda(g)$$ and becomes a completely bounded map on $L_2(\widehat{\mathbb{G}})$. In particular, being a finite von Neumann algebra this map is also well-defined on $\V$. If $1 < p < \infty$, the $L_p$-boundedness of such a map may be obtained by standard interpolation and duality arguments from a suitable $L_\infty \to \mathrm{BMO}$ inequality. In the following result we provide smoothness conditions on $m$ for this to happen.

\begin{theorem} \label{MainHorMult}
Let $\G$ be a discrete group equipped with a length $\psi: \G \to \R_+$ and set $(\H_j, \alpha_j, b_j)$ for the left and right cocycles associated to it $(j=1,2)$. Assume $\dim \H_j = n < \infty$ and let $\widetilde{m}_j: \H_j \to \C$ be lifting multipliers for $m$, so that $m = \widetilde{m}_j \circ b_j$. Assume $\widetilde{m}_j \in \mathcal{C}^{[\frac{n}{2}]+1}(\H_j \setminus\{0\})$ and $$\big|
\partial_\xi^\beta \widetilde{m}_j(\xi) \big| \, \le \, c_n |\xi|^{-
|\beta|} \quad \mbox{for all multi-index $\beta$ s.t.} \quad
|\beta| \le \mbox{$[\frac{n}{2}]+1$}.$$ Then, $T_m: L_p(\widehat{\mathbb{G}})
\to L_p(\widehat{\mathbb{G}})$ is bounded for $1 < p < \infty$ and $T_m:
\V \to \mathrm{BMO}_{\T_\psi}$.
\end{theorem}

\dem We divide it in several steps:

\noindent A. \emph{Reduction to $L_\infty \to \mathrm{BMO}$.} Assume that the hypotheses imply $L_\infty \to \mathrm{BMO}_{\T_\psi}$ boundedness. Since the condition for $\beta=0$ implies that $\widetilde{m}_1$ is bounded, the same holds for $m = \widetilde{m}_1 \circ b_1$ and we deduce the $L_2$ boundedness for $T_m$. The $L_p$ boundedness for $2 < p < \infty$ follows by interpolation from \cite{JM2}. Indeed, if we let $J_p$ the projection map onto the complemented subspace $$L_p^\circ(\widehat{\mathbb{G}}) = \Big\{ f \in L_p(\widehat{\mathbb{G}}) \, \big| \, \lim_{t \to \infty} S_{\psi,t}f = 0 \Big\},$$ we get from Theorem \ref{Interpolation} that $$J_p T_m: L_p(\widehat{\mathbb{G}}) \to L_p^\circ(\widehat{\mathbb{G}}).$$ However, $E_p = id_{L_p} - J_p$  projects onto the fixed point subspace, which is the closure of the span of $\lambda(g)$'s such that $\psi(g) = 0$. Since $\G_0 = \{g \in \G \, | \ \psi(g) = 0\}$ is a subgroup of $\G$, we deduce that $E_p$ is a conditional expectation. This implies that $T_m = m_e E_p + J_p T_m$ is also bounded. To prove the case $1 < p < 2$, we proceed by duality since $T_m^* = T_{\overline{m}}$ and the argument above also applies to $\overline{m}$.

\noindent B. \emph{Reduction to the column \textnormal{BMO} estimate.} According to the normal extension in Lemma \ref{NormalExt}, it suffices to see that $T_m: \A_\G \to \mathrm{BMO}_{\S_\psi}$. Let us assume now that $T_m: \A_\G \to \mathrm{BMO}_{\T_\psi}^c$ holds. Then $T_m$ is also a bounded map $\A_\G \to \mathrm{BMO}_{\T_\psi}$. Indeed the row BMO boundedness of $T_m$ is equivalent to the column BMO boundedness of $$T_m^\dag \Big( \sum_{g \in \G} \widehat{f}(g) \lambda(g) \Big) = T_m \Big(
\sum_{g \in \G} \overline{\widehat{f}(g)} \lambda(g^{-1}) \Big)^* = \sum_{g \in \G} \overline{m}_{g^{-1}} \widehat{f}(g) \lambda(g).$$ This shows that $$T_m^\dag = T_k \quad \mbox{with} \quad k_g = \overline{m}_{g^{-1}} = \overline{\widetilde{m}}_j \circ b_j(g^{-1}).$$ By Lemma \ref{LemLeftRight}, $k_g = \widetilde{k}_j \circ b_j$ where $\widetilde{k}_1 = \overline{\widetilde{m}}_2 \circ \Lambda_{12}$ and $\widetilde{k}_2 = \overline{\widetilde{m}}_1 \circ \Lambda_{12}^{-1}$. Since $\Lambda_{12}$ is an orthogonal transformation on $\R^n$ and the complex conjugation is harmless, it turns out that the $\widetilde{k}_j$'s satisfy one more time the same conditions as the $\widetilde{m}_j$'s.

\noindent C. \emph{Reduction to a cross product estimate.} We will only work here with the left cocycle $(\H_1, \alpha_1, b_1)$. According to the discrete topology imposed in Lemma \ref{Lemmapsi} and since $\dim \H_1 = n$, taking $\H_1 = \Rd^n$ is a suitable realization of $\H_1$. The algebra $\mathcal{L}(\H_\psi)$ is the $L_\infty$ space on the Bohr compactification. Let $\lambda_1$ and $\lambda_\rtimes$ denote the left regular representations on $\H_1$ and $\G_\rtimes = \H_1 \rtimes \G$ respectively, while $\exp b_1(g)$ will stand for $\lambda_1(b_1(g)) \simeq \exp(2\pi i \langle b_1(g), \cdot \rangle)$. Consider the trace preserving, normal homomorphism given by $\pi_1: \lambda(g) \in \V \mapsto \lambda_{\rtimes}(b_1(g) \rtimes g) \in \mathcal{L}(\G_\rtimes)$. It is very tempting and in fact very useful to use that $\mathcal{L}(\H_1)$ is commutative, by switching between the language of von Neumann algebras of discrete groups and semidirect products of von Neumann algebras. Indeed, it is a simple exercise to show that $\mathcal{L}(\G_\rtimes) \simeq \mathcal{L}(\H_1) \rtimes \G$. In particular, $\pi_1$ takes the form $$\pi_1: \lambda(g) \in \V \mapsto \exp b_1(g) \lambda(g) \in \mathcal{L}(\H_1) \rtimes \G.$$ Let $\T_\rtimes' = (S_{\rtimes,t}')_{t \ge 0}$ denote the crossed product extension $S_{\rtimes,t}' = S_t' \rtimes id_\G$ of the heat semigroup $\T' = (S_t')_{t \ge 0}$ on the Bohr compactification. It is evident that the heat semigroup is $\G$-equivariant with respect to any isometric action on $\H_1$. We now claim that it suffices to show that $T_\rtimes: \A_{\hatRd^n \rtimes \G} \to \mathrm{BMO}_{\S_\rtimes}^c$ is bounded where $$T_\rtimes \Big( \sum_{g \in \G} f_g \lambda(g) \Big) = \sum_{g \in \G} T_{\widetilde{m}_1}(f_g) \lambda(g) \quad \mbox{with} \quad T_{\widetilde{m}_1}(\exp b_1(h)) = \widetilde{m}_1 (b_1(h)) \exp b_1(h).$$ The key points are the intertwining identities $$\pi_1 \circ S_{\psi,t} \, = \, S_{\rtimes,t} \circ \pi_1 \quad \mbox{and} \quad \pi_1 \circ T_m \, = \, T_\rtimes \circ \pi_1.$$ Indeed, it is easily checked that the first one follows from $\psi(g) = \langle b_1(g), b_1(g) \rangle_{\H_1}$ while the second one from $m_g = \widetilde{m}_1(b_1(g))$. Our claim now holds for $f \in \A_\G$ as follows (note that $\pi_1(f) \in \A_{\hatRd^n \rtimes \G}$)
\begin{eqnarray*}
\|T_mf\|_{\mathrm{BMO}_{\T_\psi}^c} & = & \sup_{t > 0} \Big\| S_{\psi,t} |T_mf|^2 - |S_{\psi,t} T_m f|^2 \Big\|_{\V}^\frac12 \\ & = & \sup_{t > 0} \Big\| \pi_1 \Big( S_{\psi,t} |T_mf|^2 - |S_{\psi,t} T_m f|^2 \Big) \Big\|_{\mathcal{L}(\H_1) \rtimes \G}^\frac12 \\ & = & \sup_{t > 0} \Big\| S_{\rtimes,t} |T_\rtimes \pi_1 f|^2 - |S_{\rtimes,t} T_\rtimes \pi_1 f|^2 \Big\|_{\mathcal{L}(\H_1) \rtimes \G}^\frac12 \\ [6pt] & = & \big\| T_\rtimes (\pi_1 f) \big\|_{\mathrm{BMO}_{\T_{\rtimes}}^c} \ \lesssim \ \big\| \pi_1 f \big\|_{\mathcal{L}(\H_1) \rtimes \G} \ = \ \|f\|_{\V}.
\end{eqnarray*}
This proves the claim since we have reduced the assertion to $T_m: \A_\G \to \mathrm{BMO}_{\S_\psi}^c$.

\noindent D. \emph{Smoothness of the lifting multipliers.} Here the smoothness conditions come into play. Indeed, according to
Lemma \ref{LemStein} ii) we know that our assumptions on $\widetilde{m}_1$ imply that its Fourier inverse transform $k_{\widetilde{m}_1}$ defines a completely bounded map $\Phi$ from $L_\infty(\RR)$ to $\mathrm{BMO}_\RR^c$. However, arguing as in the proof of Lemma \ref{LemNonequiv}, we conclude that $T_\rtimes: \A_{\hatRd^n \rtimes \G} \to \mathrm{BMO}_{\T_\rtimes}^c$ is bounded. This completes the proof. \fin

\begin{remark} \label{MaincbRem}
\emph{The proof above is easily adapted to show that $T_m$ is completely bounded.}
\end{remark}

The drawback is that we need to find two lifting multipliers for the left and right cocycles. To simplify these conditions, we begin with  Theorem A ---stated in the Introduction only for left cocycles--- showing that for general discrete groups we may work with one lifting multiplier under stronger smoothness conditions.

\demA If we set $$\widetilde{m}^{\delta}(\xi) = \widetilde{m}(\xi) |\xi|^{\delta} \quad \mbox{for} \quad \xi \neq 0$$ and $\widetilde{m}^\delta(0) = 0$ with $\delta = \pm \varepsilon$, we find $| \partial_\xi^\beta \widetilde{m}^{\pm \varepsilon} (\xi) | \le c_n |\xi|^{-|\beta|}$ for all $0 \le |\beta| \le [\frac{n}{2}]+1$ by the chain rule and our hypotheses. In particular, letting $m^{\pm \varepsilon} = \widetilde{m}^{\pm \varepsilon} \circ b_\psi$ we may follow the proof of Theorem \ref{MainHorMult} to show that $T_{m^{\pm \varepsilon}}: \mathcal{L}(\G) \to \mathrm{BMO}_{\T_\psi}^c$ when $b_\psi$ is a left cocycle and $T_{m^{\pm \varepsilon}}: \mathcal{L}(\G) \to \mathrm{BMO}_{\T_\psi}^r$ when $b_\psi$ is a right cocycle. In fact these maps are completely bounded, as it follows from Remark \ref{MaincbRem}. On the other hand, we recall from \cite{JM2} that
\begin{eqnarray*}
\big[ \mathrm{BMO}_{\T_\psi}^r, L_2^\circ(\widehat{\mathbb{G}})
\big]_{2/p} & = & H_p^r(\T_\psi), \\ \big[
\mathrm{BMO}_{\T_\psi}^c, L_2^\circ(\widehat{\mathbb{G}}) \big]_{2/p} & =
& H_p^c(\T_\psi),
\end{eqnarray*}
see \cite{JM} for the definition of the Hardy spaces $H_p^r(\T_\psi)$ and $H_p^c(\T_\psi)$. Arguing as in  Theorem \ref{MainHorMult} [Point A], we get $T_{m^{\pm \varepsilon}} = \widetilde{m}^{\pm \varepsilon}(0) E_p + J_p T_{m^{\pm \varepsilon}} = J_p T_{m^{\pm \varepsilon}}$.  Therefore we conclude by interpolation that $$T_{m^{\pm \varepsilon}}: L_p(\widehat{\mathbb{G}}) \stackrel{cb}{\longrightarrow} H_p^c(\T_\psi)$$ for $2 < p < \infty$ whenever $b_\psi$ is a left cocycle and we must replace column by row if $b_\psi$ is a right cocycle. At any rate, if $A_\psi(\lambda(g)) = \psi(g) \lambda(g)$ stands for the infinitesimal generator of $\T_\psi$, we know from \cite{J22} that $$\|h\|_p \, \lesssim_{cb} \, \big\| A_\psi^{+ \gamma} h \big\|_{H_p^c}^\frac12 \big\| A_\psi^{- \gamma} h \big\|_{H_p^c}^\frac12$$ for all $\gamma > 0$ and $h \in L_p^\circ(\widehat{\mathbb{G}})$. Taking $\gamma = \varepsilon / 2$ and $h = J_p T_mf$, we see that
\begin{eqnarray*}
\|T_mf\|_p & \le_{cb} & |m_e| \|E_pf\|_p + \|J_p T_mf\|_p \\ [3pt] & \lesssim_{cb} & |m_e| \|f\|_p + \|A_\psi^{+\gamma}h\|_{H_p^c}^\frac12 \|A_\psi^{-\gamma}h\|_{H_p^c}^\frac12 \\ & = & |m_e| \|f\|_p + \|T_{m^{+\varepsilon}}f\|_{H_p^c}^\frac12 \|T_{m^{-\varepsilon}}f\|_{H_p^c}^\frac12 \ \lesssim_{cb} \ \|f\|_p.
\end{eqnarray*}
The complete boundedness for $1 < p < 2$ follows by duality as in Theorem \ref{MainHorMult}. \fin

\begin{remark} \label{RemThA}
\emph{If $\psi$ is bounded in $\G$ it suffices to know that $| \partial^\beta \widetilde{m}(\xi) \big| \le c_n |\xi|^{- |\beta|+\varepsilon}$ for all $|\beta| \le [\frac{n}{2}]+1$. If $\psi^{-1}$ is bounded in $\G \setminus \G_0 = \{g \in \G : \psi(g) \neq 0\}$, we just need to control by $|\xi|^{- |\beta|-\varepsilon}$ for the same $\beta$'s. The first condition holds for inner cocycles and the second for well-separated ones. The argument is very similar to the proof of Theorem A.}
\end{remark}

\demB As in Theorem \ref{MainHorMult}, the $L_p$-boundedness reduces to the $\A_\G \to \mathrm{BMO}_{\S_\psi}$ boundedness. Assume first that $(\H_\psi, \alpha_\psi, b_\psi)$ is a left cocycle, then the argument in Theorem \ref{MainHorMult} gives that $T_m: \A_\G \to \mathrm{BMO}_{\T_\psi}^c$ is bounded. Let us now consider the row case. One more time following our proof above, this is a matter of showing that $$T_\rtimes^\dag: \A_{\hatRd^n \rtimes \G} \to \mathrm{BMO}_{\T_\rtimes}^c$$ where $T_\rtimes = T_{\widetilde{m}} \rtimes id_\G$. As noticed in Paragraph \ref{NEQCZOSect}, we have $$T_\rtimes^\dag \big( \summ_g f_g \lambda(g) \big) = \summ_g \alpha_{\psi,g} T_{\widetilde{m}}^\dag \alpha_{\psi,g^{-1}} (f_g) \lambda(g) = \summ_g \Pi_g (f_g)
\lambda(g)$$ and $j ( \sum_g \Pi_g(f_g) \lambda(g) ) = \big( \alpha_{\psi,h^{-1}} T_{\widetilde{m}}^\dag \alpha_{\psi,h} \big) \bullet j ( \sum_g f_g \lambda(g) ) = \Phi \big( j ( \sum_g f_g \lambda(g) ) \big)$, where $$\alpha_{\psi,h^{-1}} T_{\widetilde{m}}^\dag \alpha_{\psi,h} f (x) \, = \, \Pi_{h^{-1}} f (x) \, = \, \int_{\R^n} \overline{k}_{\widetilde{m}}(\beta_hx - \beta_h y) f(y) \, dy$$ with $f(\beta_h x) = \alpha_{\psi,h^{-1}}f(x)$ and $\widehat{k}_{\widetilde{m}} = \widetilde{m}$. In particular, $T_\rtimes^\dag: \A_{\hatRd^n \rtimes \G} \to \mathrm{BMO}_{\T_\rtimes}^c$ will be bounded if the conditions in Lemma \ref{LemNonequiv} hold for $\Phi$. In fact, since the Schur product defining $\Phi$ is constant in rows, we may argue as for the proof of Lemma \ref{LemNonequiv2} ---$\Phi$ acts like a diagonal matrix by right multiplication--- and apply Lemma \ref{LemStein} i) to conclude that our smoothness condition is strong enough to imply that of Lemma \ref{LemNonequiv}. Thus, it remains to check that $$\Big\| \Big( \int_{\R^n} \Big| \Big( \Pi_{h^{-1}} f_{gh}(x) \Big) \Big|^2 \, dx \Big)^\frac12 \Big\|_{\mathcal{B}(\ell_2(\G))} \lesssim \Big\| \Big( \int_{\R^n} \Big| \Big( f_{gh}(x) \Big) \Big|^2 \, dx \Big)^\frac12 \Big\|_{\mathcal{B}(\ell_2(\G))}.$$ Indeed, the statement of Lemma
\ref{LemNonequiv} is written in terms of $\Pi_{gh^{-1}}$'s, but a quick look at the proof shows that we may replace them by $\Pi_{h^{-1}}$'s, since we have the identity $\Pi_{h^{-1}} = \alpha_{\psi,g^{-1}} \Pi_{gh^{-1}} \alpha_{\psi,g}$. On the other hand, the inequality in the operator space level follows from the argument below after matrix amplification. Thus, let us prove this inequality. Since $$\widehat{\Pi_{h^{-1}}f}(\xi) = \widehat{\alpha_{\psi,h^{-1}} \overline{k}_{\widetilde{m}}}(\xi) \widehat{f}(\xi) = \overline{\widetilde{m}(-\beta_h \xi)} \widehat{f}(\xi),$$ by Fubini and Plancherel theorems we may write the left hand side as $$\mathrm{LHS}^2 \ = \ \sup_{\|\gamma\|_{\ell_2(\G)} \le 1} \summ_g \int_{\R^n} \Big| \summ_h \widetilde{m}(- \beta_h \xi) \overline{\widehat{f}_{gh}(\xi) \gamma_h} \Big|^2 d\xi.$$ Since $\widetilde{m}(-\beta_h \xi) = \widetilde{m}(\alpha_{\psi,h^{-1}} (-\xi))$ and we are assuming that $\|\widetilde{m}\|_{schur} < \infty$, there exists a factorization $\widetilde{m}(\alpha_{\psi,h^{-1}} (-\xi)) = \langle A_{-\xi}, B_{h^{-1}} \rangle_\mathcal{K}$ and some positive constant $c$ for which $\sup_\xi \|A_\xi\|_\mathcal{K}, \ \sup_g \|B_g\|_{\mathcal{K}} \le \sqrt{c}.$ This yields
\begin{eqnarray*}
\mathrm{LHS}^2 & = & \sup_{\|\gamma\|_{\ell_2(\G)} \le 1} \summ_g \int_{\R^n} \Big| \Big\langle A_{-\xi}, \summ_h \overline{\widehat{f}_{gh}(\xi) \gamma_h} B_{h^{-1}} \Big\rangle_\mathcal{K} \Big|^2 d\xi \\ & \le & c \, \sup_{\|\gamma\|_{\ell_2(\G)} \le 1} \summ_g \int_{\R^n} \summ_j \, \Big| \summ_h \overline{\widehat{f}_{gh}(\xi) \gamma_h} B_{h^{-1}}^j \Big|^2 \, d\xi,
\end{eqnarray*}
where $B_{h^{-1}}^j$ denotes the $j$-th component of $B_{h^{-1}}$. Taking $\gamma^j = (\gamma_h \overline{B_{h^{-1}}^j})_{h \in \G}$
\begin{eqnarray*}
\mathrm{LHS}^2 & \le & c \sup_{\|\gamma\|_{\ell_2(\G)} \le 1} \summ_j \summ_g \int_{\R^n} \Big| \summ_h \widehat{f}_{gh}(\xi) \gamma_h^j \Big|^2 d\xi \\ & = & c \sup_{\|\gamma\|_{\ell_2(\G)} \le 1} \summ_j \summ_g \int_{\R^n} \Big| \summ_h f_{gh}(x) \gamma_h^j \Big|^2 dx \\ & = & c \sup_{\|\gamma\|_{\ell_2(\G)} \le 1} \summ_j \Big\langle \gamma^j, \int_{\R^n} \Big| \Big( f_{gh}(x) \Big) \Big|^2 dx \, \gamma^j \Big\rangle_{\ell_2(\G)} \ \le \ c^2 \ \mathrm{RHS}^2.
\end{eqnarray*}
This completes the proof for left cocycles. Alternatively, if we deal with a right cocycle $(\H_\psi, \alpha_\psi, b_\psi)$ everything is row/column switched. More concretely, this means that the row BMO estimate follows from our argument in Theorem \ref{MainHorMult} and the column BMO requires Lemma \ref{LemNonequiv}, details are left to the reader. \fin

\subsection{Noncommutative Riesz transforms}

Let $\G$ be a discrete group, let $\psi$ be a length function on it and construct $(\H_\psi, \alpha_\psi, b_\psi)$ to be either the left or right cocycle associated to $\psi$. The Riesz transform on $\V$ associated to an element $\eta \in \H_\psi$ is the multiplier $$R_\eta \Big( \sum_{g \in \G} \widehat{f}(g) \lambda(g) \Big) = - i \sum_{g \in \G} \frac{\langle b_\psi(g), \eta \rangle_{\H_\psi}}{\sqrt{\psi(g)}} \widehat{f}(g) \lambda(g).$$ Indeed, note that $b_\psi(g) / \sqrt{\psi(g)}$ is just the normalized vector in the direction of $b_\psi(g)$, so that the classical Riesz symbol $\widetilde{m}_\eta(\xi) = - i \langle \xi, \eta \rangle_{\H_\psi}/ \|\xi\|_{\H_\psi}$ is a lifting multiplier for $R_\eta$. The classical Mihlin condition clearly holds for $\widetilde{m}_\eta$, but it fails the more restrictive condition in Theorem A for $\varepsilon > 0$. Theorem \ref{MainHorMult} imposes alternatively to find another lifting multiplier $\widetilde{m}_\eta'$ so that $\widetilde{m}_\eta \circ b_\psi = \widetilde{m}_\eta' \circ b_\psi'$. We do not know how to find such function in general. The following result is on the contrary a simple consequence of Theorem B.

\begin{corollary} \label{RieszTransf1}
Given a discrete group $\G$, consider a length function $\psi: \G \to \R_+$ and set $(\H_\psi, \alpha_\psi, b_\psi)$ to be either the left or right cocycle associated to it. Assume that $\dim \H_\psi < \infty$, then any operator in the algebra   $\mathcal{R}$ generated by the Riesz transforms $$\mathcal{R} = \mathrm{span} \Big\{ \prod_{\eta \in
\Gamma} R_\eta \, \big| \ \Gamma \ \mbox{finite set in $\H_\psi$} \Big\}$$ defines a cb-map $\mathcal{L}(\G) \to
\mathrm{BMO}_{\T_\psi}$ and $L_p(\widehat{\mathbb{G}}) \to L_p(\widehat{\mathbb{G}})$ for all $1 < p < \infty$.
\end{corollary}

\dem Note that $\|\widetilde{m}_1 \widetilde{m}_2\|_{schur} \le \| \widetilde{m}_1\|_{schur} \|\widetilde{m}_2\|_{schur}$ by taking the Hilbertian tensor product $\mathcal{K} = \mathcal{K}_1 \otimes_2 \mathcal{K}_2$. Moreover, according to the chain rule the product $\widetilde{m}_1 \widetilde{m}_2$ satisfies the smoothness
conditions whenever $\widetilde{m}_1$ and $\widetilde{m}_2$ do. Therefore, the Fourier multipliers satisfying the hypotheses of Theorem B form an algebra. In particular, it suffices to check the conditions for a single Riesz transform $R_\eta$. We have
\begin{eqnarray*}
\widetilde{m}_\eta(\alpha_{\psi,g}(\xi)) & = & -i \frac{\langle
\alpha_{\psi,g}(\xi), \eta \rangle_{\H_\psi}}{\sqrt{\langle
\alpha_{\psi,g}(\xi), \alpha_{\psi,g}(\xi) \rangle_{\H_\psi}}} \\ & = &
-i \Big\langle \frac{\xi}{\sqrt{\langle \xi, \xi \rangle_{\H_\psi}}},
\alpha_{\psi,g^{-1}}(\eta) \Big\rangle_{\H_\psi} \ = \ \big\langle
A_\xi, B_g \big\rangle_{\R^n},
\end{eqnarray*}
with $A_\xi$ and $B_g$ satisfying the estimates $\sup_{\xi \in \R^n} |A_\xi| = 1$ and $\sup_{g \in \G}
|B_g| = |\eta|$. Hence, the assertion follows since the H\"ormander smoothness condition holds. \fin

\subsection{Mild algebraic/geometric assumptions}
\label{SubAG}

We continue our analysis just imposing the existence of one lifting multiplier. Let us prove our assertion ---in the Introduction--- that the additional $\varepsilon > 0$ in Theorem A can be removed under any of the following alternative assumptions:

\begin{itemize}
\item[i)] $\G$ is abelian,

\item[ii)] $b_\psi(\G)$ is a lattice in $\R^n$,

\item[iii)] $\alpha_\psi(\G)$ is a finite subgroup of $O(n)$,

\item[iv)]The multiplier is $\psi$-radial, i.e. $m_g = h(\psi(g))$.
\end{itemize}

\dem If $\G$ is abelian, the Hilbert space $\H_\psi$ and the inclusion map $b_\psi: \G \to \H_\psi$ coincide for both left and right cocycles, since both Gromov forms in Lemma \ref{Lemmapsi} i) coincide. Thus, the hypotheses of Theorem \ref{MainHorMult} are satisfied and we deduce the first assertion. For radial multipliers, we note that our smoothness condition implies the boundedness of $T_{\widetilde{m}}: L_\infty(\R^n) \to \mathrm{BMO}_{\R^n}$ where $\widetilde{m} = h \circ | \, |^2$. According to Lemma \ref{deLeeuwBMO}, we may replace $\R^n$ by its Bohr compactification and $\mathrm{BMO}_{\R^n}$ by the BMO space associated to the heat semigroup $\S'$ on $\hatRd^n$ in the subalgebra of trigonometric polynomials. Since radial multipliers on $\R^n$ are $\G$-equivariant with respect to any isometric action $\alpha: \G \to O(n)$, we may apply Lemma \ref{LemEquivCZO} with the cocycle action and deduce $$T_\rtimes: \A_{\hatRd^n \rtimes \G} \to \mathrm{BMO}_{\T_\rtimes'}$$ is a cb-map with $T_\rtimes = T_{\widetilde{m}} \rtimes id_\G$ and $\S_\rtimes' = \S' \rtimes id_\G$. However, as noticed in the proof of Theorem B, this is all that is really needed. When $\alpha_\psi(\G)$ is a finite subgroup of $O(n)$ we use Theorem B. By \cite{PiSim} the norm $$\big\| \widetilde{m} \big\|_{schur} \, = \, \inf_{\begin{subarray}{c}\widetilde{m}(\alpha_{\psi,g}(\xi)) = \langle A_\xi, B_g \rangle_\mathcal{K} \\ \mathcal{K} \ \mathrm{Hilbert} \end{subarray}} \Big( \sup_{\xi \in \R^n} \|A_\xi\|_\mathcal{K} \sup_{g \in \G} \|B_g\|_{\mathcal{K}} \Big)$$ coincides with the norm of the Schur multiplier $$\widetilde{m}: \sum_{\xi \in \R^n} \sum_{\gamma \in \alpha_\psi(\G)} a_{\xi, \gamma} e_{\xi, \gamma} \mapsto \sum_{\xi \in \R^n} \sum_{\gamma \in \alpha_\psi(\G)} \widetilde{m}(\gamma(\xi)) \, a_{\xi, \gamma} e_{\xi, \gamma}$$ on $L_2^c(\R^n) \otimes_h \ell_2^r(\alpha_\psi(\G))$, where $\otimes_h$ stands for the Haagerup tensor product, see \cite{P3} for precise definitions. If $\alpha_\psi(\G)$ is a finite set, we may factorize $\widetilde{m}$ as
\begin{eqnarray*}
L_2^c(\R^n) \otimes_h \ell_2^r(\alpha_\psi(\G)) &
\stackrel{id}{\longrightarrow} & L_2^c(\R^n) \otimes_h \ell_2^c
(\alpha_\psi(\G)) \\ & \stackrel{\widetilde{m}}{\longrightarrow} &
L_2^c(\R^n) \otimes_h \ell_2^c(\alpha_\psi(\G)) \\ &
\stackrel{id}{\longrightarrow} & L_2^c(\R^n) \otimes_h
\ell_2^r(\alpha_\psi(\G)),
\end{eqnarray*}
which immediately shows that $$\big\| \widetilde{m} \big\|_{schur} \, \le \, |\alpha_\psi(\G)| \sup_{(\xi,g) \in \R^n \times \G} \big| \widetilde{m}(\alpha_{\psi,g}(\xi)) \big| \, \le \, |\alpha_\psi(\G)| \, \big\| \widetilde{m} \big\|_\infty \, < \,
\infty.$$ On the other hand, the smoothness condition in Theorem B is used to ensure $$\Omega_{\widetilde{m}, \alpha_\psi} \, = \, \esssup_{x \in \R^n} \int_{|y| > 2 |x|} \, \sup_{g \in \G} \big|
k_{\widetilde{m}}(\alpha_{\psi,g} y - \alpha_{\psi,g} x) -
k_{\widetilde{m}}(\alpha_{\psi,g} y) \big| \, dy \, < \, \infty.$$
However, if $\displaystyle \Omega_{\widetilde{m}} = \esssup_{x \in \R^n}
\int_{|y| > 2 |x|} | k_{\widetilde{m}}(y - x) -
k_{\widetilde{m}}(y) | \, dy$, it is well-known that $$\big|
\partial_\xi^\beta \widetilde{m}(\xi) \big| \, \le \, c_n
|\xi|^{-|\beta|} \ \mbox{for} \ |\beta| \le \Big[ \frac{n}{2}
\Big] + 1 \quad \Rightarrow \quad \Omega_{\widetilde{m}} <
\infty.$$ In particular, since $\alpha_\psi(\G)$ is a finite set
we find that
$$\Omega_{\widetilde{m}, \alpha_\psi} \, \le \, \sum_{\alpha_{\psi,g}
\in \alpha_\psi(\G)} \esssup_{x \in \R^n} \int_{|y| > 2 |x|} \big|
k_{\widetilde{m}}(\alpha_g y - \alpha_g x) -
k_{\widetilde{m}}(\alpha_g y) \big| \, dy \, \le \, \big|
\alpha_\psi(\G) \big| \, \Omega_{\widetilde{m}}$$ is also finite, which proves assertion iii). It remains to study the case when the image $b_\psi(\G)$ lives in a lattice $\Lambda_\psi$ of the Hilbert space $\H_\psi$. If $\dim \H_\psi = n < \infty$, it is a simple observation that $\alpha_\psi(\G)$ must be a finite subgroup of $O(n)$, so that ii) follows from iii). Indeed, since there are finitely many orthogonal transformations leaving $\Lambda_\psi$ invariant, it suffices to see that $\alpha_{\psi,g}(\Lambda_\psi) \subset \Lambda_\psi$ for all $g \in \G$. We may clearly assume that $b_\psi(\G)$ generates $\H_\psi$, so that $\Lambda_\psi$ is the space of linear combinations $\sum_{h \in \G} \gamma_h b_\psi(h)$ with $\gamma_h \in \Z$. Since $$\alpha_{\psi,g}(b_\psi(h)) = b_\psi(gh) - b_\psi(g),$$ the $\Z$-linear combinations are stable under $\alpha_{\psi,g}$ for all $g$ and the claim follows. \fin

\subsection{Radial Fourier multipliers}

We now present a transference method between radial Fourier multipliers on discrete groups and their Euclidean counterparts. As usual, given a length function $\psi: \G \to \R_+$, we write $\T_\psi = (S_{\psi, t})_{t \ge 0}$ for the semigroup $\lambda(g) \mapsto \exp(-t \psi(g)) \lambda(g)$ and $\T = (S_t)_{t \ge 0}$ or $\T' = (S_t')_{t \ge 0}$ for the heat semigroup on $\R^n$ or its Bohr compactification respectively.

\demD The equivalence i) $\Leftrightarrow$ ii) follows easily from our de Leeuw compactification Lemma \ref{deLeeuwBMO}. Indeed, according to it we know the boundedness over trigonometric polynomials. The normal extension in ii) follows from Lemma \ref{NormalExt}. The normal extension in $\R^n$ is very similar. All we need to know to follow the same argument is that $H_1(\R^n)^* = \mathrm{BMO}_{\R^n}$ and that the class of Schwartz functions is dense in $H_1(\R^n)$, something which easily follows from the atomic description, see e.g. \cite{Gr}. The implication iii) $\Rightarrow$ ii) follows by taking $(\G,\psi) = (\Rd^n, | \hskip3pt |^2)$, while the argument for i) $\Rightarrow$ iii) is implicit in the proof of the result proved in Paragraph \ref{SubAG}, concluded taking the normal extension provided by Lemma \ref{NormalExt}. Indeed, if $b_\psi$ is the map associated to either the left or the right cocycle for $\psi$, we note that $$m = \widetilde{m} \circ b_\psi \quad \mbox{for} \quad m = h \circ \psi \quad \mbox{and} \quad \widetilde{m} = h \circ | \, |^2.$$ This proves the first statement. In fact, the argument for
i) $\Rightarrow$ iii) also applies assuming boundedness on the Bohr compactification instead. Moreover, a careful look at this argument shows that all that is needed is Lemma \ref{LemEquivCZO} for equivariant extension and the intertwining identities in Theorem \ref{MainHorMult}. Particularly, nothing is affected when we take $n = \infty$ as far as we remove condition i). This shows that ii) $\Leftrightarrow$ iii) even in the infinite-dimensional setting. \fin

\begin{remark}
\emph{Boundedness is equivalent to complete boundedness for all these maps. Indeed, let cb-$j$) denote the cb-version of $j$). Then, the assertion follows from the chain $\mathrm{i) \Leftrightarrow ii) \Leftrightarrow iii)} \Rightarrow \mbox{cb-iii)} \Rightarrow \mbox{cb-ii)} \Rightarrow \mbox{cb-i)}$. The implication cb-iii) $\Rightarrow$ cb-ii) is trivial, while the last implication follows again from Lemma \ref{deLeeuwBMO}. Therefore, it suffices to show that ii) $\Rightarrow$ cb-iii) which follows again the last statement in Lemma \ref{LemEquivCZO} and the intertwining identities. This completes the proof.}
\end{remark}

\begin{remark}
\emph{It is standard that $$\big| \partial_\xi^\beta
\widetilde{m}(\xi) \big| \, \le \, c_n |\xi|^{- |\beta|} \quad
\Rightarrow \quad \sup_{R>0} \Big( \frac{1}{R^{n-2|\beta|}}
\int_{R < |\xi| < 2R} \big| \partial_\xi^\beta \widetilde{m}(\xi)
\big|^2 \, d\xi \Big)^\frac12 \, \le \, c_n.$$ If the inequality
on the right holds for all $|\beta| \le [\frac{n}{2}] +1$, we say
that $\widetilde{m}$ satisfies H\"ormander's smoothness condition.
This condition also implies the $L_p$ as well as the $L_\infty \to
\mathrm{BMO}$ boundedness of the Fourier multiplier
$T_{\widetilde{m}}$ on $\R^n$. Thus, by Theorem D
we see that whenever $\widetilde{m}$ satisfies the (weaker)
H\"ormander smoothness condition and $\widetilde{m} = h \circ | \,
|^2$, the Fourier multipliers $T_{h \circ \psi}$ are $L_p$ and
$L_\infty \to \mathrm{BMO}$ bounded for any discrete group $\G$
with $\dim \H_\psi = n$.}
\end{remark}

\section{{\bf Littlewood-Paley theory}} \label{SectLP}

We now prove some square function estimates. The boundedness of new square functions for noncommutative martingale transforms and semicommutative CZO's was recently investigated in \cite{MP}. The smoothness assumptions there were needed for additional weak-type $(1,1)$ estimates, here we will find weaker conditions. Set $\RR_1 = L_\infty(\R^n) \bar\otimes \M_1$ and $\RR_{12} = L_\infty(\R^n) \bar\otimes \M_1 \bar\otimes \M_2$ with $\M_1$ and $\M_2$ semifinite algebras. Consider a CZO with kernel representation $$\Lambda f(x) = \int_{\R^n} k(x,y) \otimes f(y) \, dy = \int_{\R^n} \widetilde{k}(x,y) (f(y)) \, dy,$$ where $\widetilde{k}(x,y) (\cdot) = k(x,y) \otimes \cdot$ and $k$ takes values in $\M_2$. If
\begin{itemize}
\item $\Lambda: L_\infty(\M_1; L_2^c(\R^n)) \to L_\infty(\M_1 \bar\otimes \M_2; L_2^c(\R^n))$,

\vskip5pt

\item $\displaystyle \esssup_{x_1,x_2} \int_{|x_1-y| > 2 |x_1-x_2|} \big\| k(x_1,y) - k(x_2,y) \big\|_{\M_2} \ dy
\, < \, \infty$,
\end{itemize}
we deduce from Lemma \ref{LEMSCCZ} that $\Lambda: \RR_1 \to \mathrm{BMO}_{\RR_{12}}^c$. Take $\M_2 = \mathcal{B}(\ell_2)$ and $k(x,y) = \sum_j k_j(x,y) \otimes (e_{1j} \oplus_\infty e_{j1})$, where the $k_j$'s are scalar-valued and $e_{ij}$ stands for the $(i,j)$-th matrix unit. Consider the CZO $L_j$ associated to $k_j$ and such that $\Lambda = \sum_j L_j \otimes (e_{1j} \oplus_\infty e_{j1})$. The column part $\Lambda_c = \sum_j L_j \otimes e_{j1}$ satisfies the first condition if $\|\Lambda_cf\|_2^2 = \sum_j \|L_jf\|_2^2 \lesssim \|f\|_2^2$ since the kernel acts by left multiplication, see Remark \ref{Op-Vld model}. For the row part, we use some basic operator space theory \cite{P3}. Namely, the condition is equivalent to the complete boundedness of $\Lambda_r: L_2^c(\R^n) \to L_2^c(\R^n) \otimes_h R$.
In particular, such a map defines an element in $L_2^c(\R^n) \bar\otimes L_2^r(\R^n) \otimes_h R$ with norm $\|\sum_j L_j L_j^*\|^{1/2}$. This leads to the same condition with $L_j^*$ in place of $L_j$. Finally $$\esssup_{x_1,x_2} \int_{|x_1-y| > 2 |x_1-x_2|} \Big( \sum_{j=1}^\infty \big| k_j(x_1,y) - k_j(x_2,y) \big|^2 \Big)^\frac12 \, dy \, < \, \infty$$ is the form of the smoothness assumption. By the symmetry of the kernel, the map $\Lambda^\dag(f) = \Lambda(f^*)^*$ essentially equals $\Lambda$ and $\RR_1 \to \mathrm{BMO}_{\RR_{12}}$ boundedness follows with no extra assumptions. $L_p$-boundedness follows by interpolation and duality if the H\"ormander condition holds also on the second variable. In particular, if we let $\RR = L_\infty(\R^n) \bar\otimes \M$, we recover the main result in \cite{MP} in terms of the spaces $L_p(\RR; \ell_{rc}^2)$ defined in the Introduction.

\begin{lemma} \label{LPSemicommt}
Let $$\Lambda f(x) = \int_{\R^n} k(x,y) \otimes f(y) \, dy = \summ_j L_j f(x) \otimes (e_{1j} \oplus_\infty e_{j1})$$ be the $\mathrm{CZO}$ above. Assume that
\begin{itemize}
\item[i)] $\displaystyle \sum_{j=1}^\infty \|L_jf\|_2^2 +
\|L_j^*f\|_2^2 \ \lesssim \, \|f\|_2^2$ for $f \in L_2(\RR)$,

\vskip5pt

\item[ii)] $\displaystyle \esssup_{x_1,x_2} \int_{|x_1-y| > 2 |x_1-x_2|} \big\|
k(x_1,y) - k(x_2,y) \big\|_{\ell_2} \, dy \, < \, \infty$,

\vskip5pt

\item[iii)] $\displaystyle \esssup_{x_1,x_2} \int_{|x_1-y| > 2 |x_1-x_2|} \big\|
k(y,x_1) - k(y,x_2) \big\|_{\ell_2} \, dy \, < \, \infty$.
\end{itemize}
Then $\Lambda: \RR_1 \to \mathrm{BMO}_{\RR_2}$ is bounded and the inequality below holds for $1 < p < \infty$ $$\Big\|
\sum_{j=1}^\infty L_jf \otimes \delta_j \Big\|_{L_p(\RR_1; \ell_{rc}^2)} \lesssim \, \frac{p^2}{p-1} \, \|f\|_{L_p(\RR_1)}.$$
\end{lemma}
This gives the $L_p$-boundedness of operator-valued $g$-functions and Lusin square functions, see \cite{MP} for more applications. The conditions above hold for convolution
maps with kernels satisfying $( \sum_j | \partial_\xi^\beta \, \widehat{k}_j(\xi) |^2 )^\frac12 \le c_n \, |\xi|^{-|\beta|}$ for $|\beta| \le [\frac{n}{2}] + 1$, which is just a form of H\"ormander-Mihlin multiplier theorem for $\ell_2$-valued kernels.

\begin{lemma} \label{LemmaMT}
Let $\psi: \G \to \R_+$ be a length function with $\dim \H_\psi = n$. Let $\Gamma$ stand for the free group $\mathbb{F}_\infty$ with infinitely many generators
$\gamma_1, \gamma_2, \ldots$ and left regular representation $\lambda_\Gamma$. Consider a sequence of functions $(h_j)_{j \ge 1}$ in $\mathcal{C}^{k_n}(\R_+ \setminus \{0\})$ for $k_n = [\frac{n}{2}]+1$ such that $$\Big( \sum_{j=1}^\infty \Big| \frac{d^k}{d\xi^k} h_j(\xi) \Big|^2\Big)^\frac12 \, \le \, c_n \, |\xi|^{-k} \quad \mbox{for
all} \quad k \le \Big[ \frac{n}{2} \Big] +1.$$ Let $s_j: \R^n \setminus \{0\} \to \C$ given by $\widehat{s}_j(\xi) = h_j(|\xi|^2)$. Then, we have a cb-map $$\Phi:
L_\infty(\R^n) \bar\otimes \mathcal{L}(\Gamma) \ni \sum_{\gamma \in \Gamma} f_\gamma \otimes \lambda_\Gamma(\gamma) \mapsto \sum_{j=1}^\infty \overline{s}_j * f_{\gamma_j} \in \mathrm{BMO}_{\R^n}.$$
\end{lemma}

\dem According to the noncommutative Khintchine inequality for free generators \cite{P3}, the map $e_{1j} \oplus_\infty e_{j1} \mapsto \lambda_\Gamma(\gamma_j)$ is a complete isomorphism and the span of $\lambda_\Gamma(\gamma_j)$'s is completely complemented in $\mathcal{L}(\Gamma)$. In particular, it suffices to show that we
have a cb map $L_\infty(\R^n) \bar\otimes \mathcal{B}(\ell_2) \ni \sum_j f_j \otimes \big( e_{1j} \oplus_\infty e_{j1} \big) \mapsto \sum_j \overline{s}_j * f_j \in \mathrm{BMO}_{\R^n}$. Since we have an intersection of row and column at both sides, it is enough to prove the row-row and column-column complete boundedness. By symmetry, we just consider the column case, so we are reduced to show that $\sum_j f_j \otimes e_{j1} \mapsto \sum_j \overline{s}_j * f_j$ defines a cb-map $C(L_\infty(\R^n)) \to \mathrm{BMO}_{\R^n}^c$, where $C(L_\infty(\R^n)) = L_\infty(\R^n; \ell_2^c)$ and $C$ stands for the column subspace of $\mathcal{B}(\ell_2)$. Recall the following simple isometries for a sequence of matrix-valued functions $(f_j)_{j \ge 1}$ in $\RR = L_\infty(\R^n) \bar\otimes \mathcal{B}(\ell_2)$
\begin{itemize}
\item[i)] $\displaystyle \Big\| \sum_{j=1}^\infty f_j \otimes e_{j1} \Big\|_{C(\RR)} = \sup_{\|\xi\|_{\ell_2} \le 1} \Big\| \Big( \sum_{j=1}^\infty \|f_j \xi \|_{\ell_2}^2 \Big)^\frac12 \Big\|_{L_\infty(\R^n)}$,

\item[ii)] $\displaystyle \Big\| \underbrace{\sum_{j=1}^\infty \overline{s}_j * f_j}_f \Big\|_{\mathrm{BMO}_\RR^c} = \sup_{Q \ \mathrm{ball}} \sup_{\|\xi\|_{\ell_2} \le 1} \Big( \mean_Q \big\| (f \xi) - (f\xi)_Q \big\|_{\ell_2}^2 \Big)^\frac12$.
\end{itemize}
With this in mind, it suffices to show that $$\Big\| \Big( \sum_{j=1}^\infty \big\| s_j * \varphi \big\|_{\ell_2}^2 \Big)^\frac12 \Big\|_{L_1(\R^n)} \, \lesssim \, \|\varphi\|_{H_1(\ell_2)}$$ for $\ell_2$-valued functions, since this implies the predual inequality and we conclude taking adjoints. By the atomic characterization of $H_1(\ell_2)$, we may write its norm as $\inf \sum_k |\lambda_k|$ where the infimum runs over all possible decompositions $\varphi = \sum_k \lambda_k a_k$ as a linear combination of atoms $a_k$, which are mean zero functions $\R^n \to \ell_2$ supported by cubes and such that $\|a_k\|_{L_2(\ell_2)} \le | \mathrm{supp} \, a_k |^{-1/2}$. By the triangle inequality, it suffices to see that the left hand side is $\lesssim 1$ when $\varphi$ is an arbitrary atom $a$ supported by an arbitrary cube $Q$. We have $$\int_{\R^n} \Big( \sum_{j=1}^\infty \big\| s_j * a (x) \big\|_{\ell_2}^2 \Big)^\frac12 \, dx \, = \, \int_{5Q} + \int_{\R^n \setminus 5Q} \, = \, \mathrm{A} + \mathrm{B}.$$ Our hypotheses easily give $$\Big( \sum_{j=1}^\infty \big| \partial_{\xi}^\beta \widehat{s}_j(\xi) \big|^2 \Big)^\frac12 \, \le \, c_n \, |\xi|^{-|\beta|} \quad \mbox{for all $\beta$ such that} \quad
|\beta| \le \Big[ \frac{n}{2} \Big] + 1.$$ As we remarked before the statement of this result, this implies the hypotheses of Lemma \ref{LPSemicommt}. In particular, we have $\sum_j \|s_j * f\|_2^2 \lesssim \|f\|_2^2$. This, together with H\"older's inequality gives rise to
\begin{eqnarray*}
\mathrm{A} & \le & \sqrt{|5Q|} \, \Big( \int_{5Q} \sum_{j=1}^\infty \big\| s_j * a(x) \big\|_{\ell_2}^2 \, dx \Big)^\frac12 \\ & \le & \sqrt{|5Q|} \, \Big( \sum_{k=1}^\infty \sum_{j=1}^\infty \big\| s_j * a_k(x) \big\|_2^2 \Big)^\frac12 \ \lesssim \ \sqrt{|5Q|} \, \Big( \int_{\R^n} \big\| a(x) \big\|_{\ell_2}^2 \, dx \Big)^\frac12 \ \lesssim \ 1,
\end{eqnarray*}
for $a = (a_k)_{k \ge 1}$. On the other hand, using the mean-zero condition
\begin{eqnarray*}
\mathrm{B} & = & \int_{\R^n \setminus 5Q} \Big( \sum_{j=1}^\infty \Big\| \int_Q \big( s_j(x-y) - s_j(x-c_Q) \big) a(y) \, dy \Big\|_{\ell_2}^2 \Big)^\frac12 \, dx \\ & \le & \int_Q \Big[ \int_{\R^n \setminus 5Q} \Big( \sum_{j=1}^\infty \big| s_j(x-y) - s_j(x-c_Q) \big|^2 \Big)^\frac12 \, dx \Big] \, \big\| a(y) \big\|_{\ell_2} \, dy \\ & \lesssim & \int_Q \big\| a(y) \big\|_{\ell_2} \, dy \ \le \ \sqrt{|Q|} \, \Big( \int_{\R^n} \big\| a(y) \big\|_{\ell_2}^2 \, dy \Big)^\frac12 \ \le \ 1,
\end{eqnarray*}
according to condition iii) in Lemma \ref{LPSemicommt}, which holds as a consequence of the H\"ormander-Mihlin condition in the statement. The estimates for A and B show that the predual inequality holds and the proof is complete. \fin

\begin{theorem} \label{ThNCLP}
Given $f = \sum_{g} \widehat{f}(g) \lambda(g)$, let $$\Lambda_\psi f = \summ_j L_{\psi,j} f \otimes \lambda_\Gamma(\gamma_j) = \summ_{g,j} h_j(\psi(g)) \widehat{f}(g) \lambda(g) \otimes \lambda_\Gamma(\gamma_j),$$ where the $h_j$'s satisfy the same smoothness conditions as above. Then$\, :$
\begin{itemize}
\item[i)] If $\G_\Gamma = \G \times \Gamma$, we have cb-maps $$\hskip25pt \Lambda_\psi: \V \stackrel{cb}{\longrightarrow} \mathrm{BMO}_{\T_{\psi,
\otimes}}(\mathcal{L}(\G_\Gamma)) \quad \mbox{and} \quad \Lambda_\psi: L_p(\widehat{\mathbb{G}}) \stackrel{cb}{\longrightarrow} L_p(\widehat{\mathbb{G}}_\Gamma)$$ for $1 < p < \infty$, where $\T_{\psi, \otimes} = (S_{\psi,t} \otimes id_{\mathcal{L}(\Gamma)})_{t \ge 0}$. In particular, we obtain $$\Big\| \sum_{j=1}^\infty L_{\psi,j} f \otimes \delta_j \Big\|_{L_p(\widehat{\mathbb{G}}; \ell_{rc}^2)} \, \le_{cb} \, c_p \, \|f\|_{L_p(\widehat{\mathbb{G}})}.$$

\vskip5pt

\item[ii)] Additionally, we have $$\hskip35pt \sum_{j=1}^\infty |h_j(\xi)|^2 = 1 \quad \Rightarrow \quad \|f\|_{L_p(\widehat{\mathbb{G}})} \, \le_{cb} \, c_p \, \Big\| \sum_{j=1}^\infty L_{\psi,j} f \otimes \delta_j \Big\|_{L_p(\widehat{\mathbb{G}}; \ell_{rc}^2)}.$$
\end{itemize}
\end{theorem}

\dem By our smoothness assumption on the $h_j$'s, the map $$f \in L_\infty(\R^n) \mapsto \sum_{j=1}^\infty (s_j * f) \otimes (e_{1j} \oplus_\infty e_{j1}) \in \mathrm{BMO}_{\mathcal{R}}$$ is completely bounded from $L_\infty(\R^n)$ to $\mathrm{BMO}_{\mathcal{R}}$ with $\mathcal{R} = L_\infty(\R^n) \bar\otimes \mathcal{B}(\ell_2)$, since the smoothness of the $h_j$'s imply the hypotheses of Lemma \ref{LPSemicommt}. Arguing as in Lemmas \ref{HeatBMO} and \ref{deLeeuwBMO} ---like in the start of the proof of Lemma \ref{LemNonequiv}--- in conjunction with the Khintchine inequalities for the free group generators \cite{HP} implies $$f \in \A_{\hatRd^n} \mapsto \sum_{j=1}^\infty T_{\widetilde{m}_j} f \otimes \lambda_\Gamma(\gamma_j) \in \mathrm{BMO}_{\S_\otimes'}(L_\infty(\R^n) \bar\otimes \mathcal{L}(\Gamma))$$ for the multipliers $\widetilde{m}_j = h_j \circ | \ |^2$ and $\S_\otimes' = (S_t' \otimes id_{\mathcal{L}(\Gamma)})_{t \ge 0}$, the tensor product extension of the heat semigroup on the Bohr compactification of $\R^n$. Now, since the involved multipliers are radial, the whole map is $\G$-equivariant for any orthogonal action and Lemma \ref{LemEquivCZO} yields the crossed product extension $$T': \sum_{g \in \G} f_g \lambda(g) \in \A_{\hatRd^n \rtimes \G} \stackrel{cb}{\longmapsto} \sum_{j=1}^\infty \sum_{g \in \G} T_{\widetilde{m}_j} (f_g) \lambda(g) \otimes \lambda_\Gamma(\gamma_j) \in \mathrm{BMO}_{\S_\rtimes'}(\RR_{\rtimes \G_\Gamma}),$$ where $\RR_{\rtimes \G_\Gamma} = (L_\infty(\hatRd^n) \rtimes \G) \bar\otimes \mathcal{L}(\Gamma)$ and $\T_{\rtimes} = ( S_{t} \rtimes id_{\mathcal{L}(\G_\Gamma)} )_{t \ge 0}$. Note that $$T' \circ \pi_\psi = (\pi_\psi \otimes id_{\mathcal{L}(\Gamma)})  \circ \Lambda_\psi$$ using the embedding $\pi_\psi: \V \to \mathcal{L}(\H_\psi) \rtimes \G = L_\infty(\hatRd^n) \rtimes \G$ from the proof of Theorem \ref{MainHorMult}. Indeed, recall that $T_{\widetilde{m}_j}(\exp b_\psi(g)) = h_j(\psi(g)) \, \exp b_\psi(g)$. This immediately gives $$\Lambda_\psi: \A_\G \stackrel{cb}{\longrightarrow} \mathrm{BMO}_{\T_{\psi, \otimes}}(\mathcal{L}(\G_\Gamma)).$$ Finally, the same argument as in Lemma \ref{NormalExt} provides a unique normal extension $\Lambda_\psi: \V \to \mathrm{BMO}_{\S_{\psi,\otimes}}(\mathcal{L}(\G_\Gamma))$. The complete boundedness on $L_2$ follows immediately from the smoothness condition for $k=0$ on the $h_j$'s. Thus, by interpolation we obtain the complete boundedness for $2 < p < \infty$. The case $1 < p < 2$ is slightly different because $\Lambda_\psi$ is not self-dual $\Lambda_\psi^* (\sum_{\gamma \in \Gamma} f_\gamma \otimes \lambda_\Gamma(\gamma)) = \sum_j L_{\psi,j}^*(f_{\gamma_j})$. Since the $L_2$ boundedness is clear, we claim it suffices to show that $\Lambda_\psi^*: \mathcal{L}(\G_\Gamma) \to \mathrm{BMO}_{\T_\psi}$ is completely bounded. Indeed, arguing once more as in Theorem \ref{MainHorMult} we find $$\Lambda_\psi^* \Big( \underbrace{\sum_{\gamma \in \Gamma} f_\gamma \otimes \lambda_\Gamma(\gamma)}_{f} \Big) = J_p \Lambda_\psi^*(f) + \sum_{g \in \G_0} \sum_{j=1}^\infty \overline{h_j(0)} \, \widehat{f}_{\gamma_j}(g) \lambda(g),$$ where $J_p: L_p(\widehat{\mathbb{G}}) \to L_p^\circ (\widehat{\mathbb{G}})$ and $\G_0 = \{g \in \G \, | \, \psi(g) = 0\}$. The first term on the right is completely bounded on $L_p$ by interpolation. To estimate the $L_p$-norm of the second term we use Cauchy-Schwarz, the conditional expectation $\mathcal{E}_0$ onto the closure of $\mathrm{span} \, \lambda(\G_0)$, the noncommutative Khintchine inequality for free generators and the fact that the span of the $\lambda_\Gamma(\gamma_j)$'s is completely complemented in $L_p(\mathcal{L}(\Gamma))$, see e.g. \cite{PP,P3}. Altogether gives rise to the following estimate
\begin{eqnarray*}
\Big\| \sum_{g \in \G_0} \sum_{j=1}^\infty \overline{h_j(0)} \, \widehat{f}_{\gamma_j}(g) \lambda(g) \Big\|_p & \le_{cb} & \Big( \sum_{j=1}^\infty |h_j(0)|^2 \Big)^\frac12 \Big\| \Big( \sum_{j=1}^\infty |\mathcal{E}_0(f_{\gamma_j})|^2 \Big)^\frac12 \Big\|_p \\ & \le_{cb} & \Big( \sum_{j=1}^\infty |h_j(0)|^2 \Big)^\frac12 \Big\| \Big( \sum_{j=1}^\infty |f_{\gamma_j}|^2 \Big)^\frac12 \Big\|_p \\ & \lesssim_{cb} & \Big( \sum_{j=1}^\infty |h_j(0)|^2 \Big)^\frac12 \Big\| \sum_{j=1}^\infty f_{\gamma_j} \otimes \lambda_\Gamma(\gamma_j) \Big\|_p \\ & \lesssim_{cb} & \Big( \sum_{j=1}^\infty |h_j(0)|^2 \Big)^\frac12 \Big\| \sum_{\gamma \in \Gamma} f_{\gamma} \otimes \lambda_\Gamma(\gamma) \Big\|_p.
\end{eqnarray*}
This proves the claim. For the $L_\infty \to \mathrm{BMO}$ estimate, we recall that $$\Phi: L_\infty(\R^n) \bar\otimes \mathcal{L}(\Gamma) \ni \sum_{\gamma \in \Gamma} f_\gamma \otimes \lambda_\Gamma(\gamma) \mapsto \sum_{j=1}^\infty \overline{s}_j * f_{\gamma_j} = \sum_{j=1}^\infty T_{\overline{\widetilde{m}}_j} (f_{\gamma_j}) \in \mathrm{BMO}_{\R^n}$$ is completely bounded from Lemma \ref{LemmaMT}. Replacing again $\R^n$ by $\hatRd^n$ and $\mathrm{BMO}_{\R^n}$ by $\mathrm{BMO}_{\T'}$, we use that $\Phi$ is $\G$-equivariant ($s_j$ is radial) with respect to the natural action $\alpha_\psi$ and apply Lemma \ref{LemEquivCZO}. This shows that $\Phi \rtimes id_\G: \RR_{\rtimes \G_\Gamma} \to \mathrm{BMO}_{\T_{\rtimes}'}(\RR_{\rtimes \G})$ is completely bounded. Then observe that $\Phi \rtimes id_\G = T'^*$ and the intertwining identity $T'^* \circ (\pi_\psi \otimes id_{\mathcal{L}(\Gamma)}) = \pi_\psi \circ \Lambda_\psi^*$ still holds. Hence, $\Lambda_\psi^*: \A_{\G_\Gamma} \to \mathrm{BMO}_{\T_\psi}$ is a cb-map and (after normal extension) we deduce $\Lambda_\psi$ is completely bounded on $L_p$ for $1 < p < \infty$. Thus, we conclude $$\Big\| \sum_{j=1}^\infty T_j f \otimes \delta_j \Big\| _{L_p(\widehat{\mathbb{G}}; \ell_{rc}^2)} \, \le_{cb} \, c_p \, \|f\|_{L_p(\widehat{\mathbb{G}})}$$ according to the noncommutative Khintchine inequality for free generators. The proof of ii) is straightforward. Indeed, if $\sum_j |h_j(\xi)|^2 = 1$ it is clear that we find an isometry $\|\Lambda_\psi f\|_2 = \|f\|_2$. By polarization, $\langle f_1, f_2 \rangle_{L_2(\widehat{\mathbb{G}})} = \langle \Lambda_\psi f_1, \Lambda_\psi f_2 \rangle_{L_2(\widehat{\mathbb{G}}_\Gamma)}$. Therefore, if $f \in L_2(\widehat{\mathbb{G}}) \cap L_p(\widehat{\mathbb{G}})$ we see that $$\|f\|_p \, = \, \sup \Big\{ \big\langle \Lambda_\psi f, \Lambda_\psi g \big\rangle_{L_2(\widehat{\mathbb{G}}_\Gamma)} \, \big| \ g \in L_2(\widehat{\mathbb{G}}) \cap L_{p'}(\widehat{\mathbb{G}}), \ \|g\|_{p'} \le 1 \Big\} \, \lesssim \, \|\Lambda_\psi f\|_p.$$ By density, this inequality still holds in the whole $L_p(\widehat{\mathbb{G}})$. Moreover, the same estimate is valid after matrix amplification and we deduce the assertion once more by means of the noncommutative Khintchine inequality for free generators. \fin

\section{{\bf Examples and comments}}
\label{SectApp}

We finally illustrate our main results in a variety of scenarios. This includes the classical groups $\mathbb{T}^n$ and $\R^n$, noncommutative tori or the free group algebra. We also provide new examples of Rieffel's quantum metric spaces. 

\subsection{The $n$-torus}

Since $$\summ_{j,k} \overline{\beta}_j \beta_k \, e^{-t \|k-j\|^2} = \Big( \frac{\pi}{t} \Big)^{\frac{n}{2}} \int_{\R^n} e^{- \pi^2 \|x\|^2/t} \Big| \summ_j \beta_j \, e^{2 \pi i \langle j,x \rangle} \Big|^2 \, dx \ge 0,$$ Schoenberg's theorem gives that $\psi(k) = \|k\|^2$ is a length function on $\Z^n$. Being an abelian group, both Gromov products $K_\psi^1$ and $K_\psi^2$ coincide, so that there is just one Hilbert space $\H_\psi$ and one inclusion map $b_\psi: \G \to \H_\psi$. In the specific case considered, the inner product takes the form $$\Big\langle \sum_{j \in \Z^n} a_j \delta_j, \sum_{k \in \Z^n} a_k \delta_k \Big\rangle_{\H_\psi} = \sum_{j,k \in \Z^n} a_j a_k \, \big\langle j,k \big\rangle_{\R^n} = \Big\| \sum_{j \in \Z^n} a_j j \Big\|_{\R^n}^2.$$ According to Lemma \ref{Lemmapsi}, we have to quotient out the subspace of finitely supported sequences $(a_j)_{j \in \Z^n}$ for which $\sum_j a_j j = 0$. It is easily checked that the resulting quotient is $n$-dimensional ---so that $\H_\psi \simeq \R^n$--- and the map $b_\psi: \Z^n \to \R^n$ becomes the canonical inclusion. This cocycle is equipped with the trivial action $\alpha_{\psi,k} = id_{\R^n}$ for all $k \in \Z^n$. In particular, we find that $|\alpha_\psi(\Z^n)| = 1$ in this case and Theorem A $(\varepsilon=0)$ meets exactly the classical Mihlin condition, so that we recover the original formulation of H\"ormander-Mihlin theorem for $\mathbb{T}^n$.

\subsection{Euclidean multipliers}
\label{SubFMTN}

Taking $\G = \R^n_{\mathrm{disc}}$ the Euclidean space equipped with the discrete topology and recalling K. de Leeuw's compactification theorem \cite{dL} we obtain from Theorem A for $\G$ abelian ($\varepsilon = 0$) that $$T_mf (x) = \int_{\R^n} \widetilde{m}(b(\xi)) \hskip1pt \widehat{f}(\xi) \hskip1pt e^{2\pi i \langle x, \xi \rangle} \, d\xi$$ is $L_p(\R^n)$-bounded for any cocycle $b: \R^n \to \R^d$ with $\widetilde{m}$ Mihlin-smooth of degree $[\frac{d}{2}]+1$. By picking suitable cocycles, this is related to H\"ormander-Mihlin theorem \cite{Ho,Mi} and de Leeuw's restriction/periodization theorems \cite{dL}, see below. Also new $L_\infty \to \mathrm{BMO}$ estimates will be given. Before illustrating our point, let us analyze the variety of finite-dimensional cocycles of $\R^n$. To construct a generic $d$-dimensional cocycle for $\G= \R^n$, assume that we have $(n,d)  = (n_1,d_1) + (n_2,d_2)$ with $n_j = 0$ iff $d_j = 0$ for $j=1,2$. Consider a triple $\Sigma = (\eta, \pi, \gamma)$ composed by $\eta \in \R^{d_1}$, a representation $\pi: \R^{n_1} \to O(d_1)$ and a group homomorphism $\gamma: \R^{n_2} \to \R^{d_2}$. Then $b_{\Sigma}(\xi) = b_{\Sigma}(\xi_1 \oplus \xi_2) = ( \pi(\xi_1) \eta - \eta ) \oplus \gamma(\xi_2)$ is a cocycle of $\R^n$ with $\H_{\Sigma} = \R^d$ and $\alpha_{\Sigma, \xi} = \pi(\xi_1) \oplus id_{\R^{d_2}}$. In fact, all possible cocycles $\R^n \to \R^d$ break up into an orthogonal sum of an inner and a proper part (any of which may vanish) as above. The proper part is always associated to the trivial action. This characterization is not hard and it may be folklore. It was already noticed in \cite{dCTV1} and a proof can be easily reconstructed from \cite[Exercise 4.5]{Val}. 

\vskip3pt

\noindent \textbf{1.} \emph{Mihlin theorem.} Apply Theorem A $(\varepsilon = 0)$ with the trivial cocycle $\R^n \to \R^n$.

\vskip3pt

\noindent \textbf{2.} \emph{de Leeuw's restriction theorem.} K. de Leeuw proved in \cite{dL} that the restriction to $\R^k$ of any sufficiently smooth function $m: \R^n \to \C$ which defines an $L_p$-bounded Fourier multiplier, is also $L_p$-bounded. In our setting, this corresponds (under Mihlin regularity of the original multiplier) to take the standard cocycle $\R^k \to \R^n$ given by the inclusion map with the trivial action. The same argument works for restriction onto integer lattices $\Z^k$ or affine deformations of it.

\vskip3pt

\noindent \textbf{3.} \emph{de Leeuw's periodization theorem.} Another consequence of de Leeuw's approach is that $\Z^n$-periodizations of $L_p$-multipliers in $\R^n$ supported by the unit cube remain in the same class, see also Jodeit \cite{Jo}. Under Mihlin regularity of the original multiplier, this corresponds to the (inner) cocycle $\xi \in \R^n \mapsto \sum (e^{2\pi i \xi_j} - 1) e_j \in \R^{2n}$ with the corresponding action as described above. Here the lifting multiplier is taken to coincide with the original multiplier in an $n$-torus of $\R^{2n}$ and smoothly truncated outside it. Our $L_\infty \to \mathrm{BMO}$ estimate is apparently new.

\vskip3pt

\noindent \textbf{4.} \emph{Directional multipliers.} Taking $u_\gamma = (\gamma_1, \gamma_2, \ldots, \gamma_n) \in \R^n$ and the $1$-dimensional cocycle $b_\Sigma: \xi \in \R^n \mapsto \summ \xi_j \gamma_j \in \R$, we just need to control $1$ derivative of the lifting multiplier $\widetilde{m}$. Letting $\gamma_j = \delta_{j = j_0}$ we obtain multipliers depending only on the $j_0$-th coordinate. Taking $\gamma_1, \gamma_2, \ldots, \gamma_n$ to be $\Z$-independent, we obtain injective cocycles in $\Z^n$ (i.e. arbitrary smooth multipliers) and multipliers depending only on the direction $u_\gamma$ in $\R^n$. Our $L_\infty \to \mathrm{BMO}$ estimates are of particular interest when $\widetilde{m}(\xi) = - i \mathrm{sgn}(\xi)$ and $b_\Sigma$ as above. Namely, it turns out that $m_\xi = \widetilde{m}(b_\Sigma(\xi))$ induces the directional Hilbert transform $H_{u_\gamma}$ in the direction of $u_\gamma$. It is well-known that $H_{u_\gamma}$ is not $L_\infty \to \mathrm{BMO}$ bounded for the classical BMO space. However, Theorem \ref{MainHorMult} provides the alternative space $\mathrm{BMO}_{u_\gamma} = \mathrm{BMO}_{\T_\psi}$ for $\psi (\xi) = |\langle \xi, u_\gamma \rangle |^2$. Recall that this BMO space interpolates with $L_p$ and thus provides the right endpoint estimate for the directional Hilbert transform. Working with proper $d$-dimensional cocycles we obtain some generalizations for $1 \le d \le n$.

\subsection{Donut type multipliers.} 

In fact, other cocycles provide a large family of exotic $L_p$ multipliers in $\R^n$. The same construction applies in $\mathbb{T}^n$. As an illustration consider the cocycle $$b(\xi) = (\cos 2 \pi \alpha \xi - 1, \sin 2 \pi \alpha \xi, \cos 2 \pi \beta \xi - 1, \sin 2 \pi \beta \xi)$$ for some $\alpha, \beta \in \R_+$. Theorem A shows that the restriction of a Mihlin multiplier in $\R^4$ to this donut helix will be an $L_p$ multiplier on $\R$. It is useful to compare it with de Leeuw's periodization theorem for compactly supported multipliers. The main difference here is the irregularity obtained from choosing $\alpha/\beta$ irrational, leading to a geodesic flow with dense orbit. Hence, $m$ oscillates infinitely often with no periodic pattern. Taking for instance $0 < \gamma < \frac12$ and $\widetilde{m}(\gamma) = |\gamma|^{2\gamma}$ for $\gamma$ small and $\widetilde{m}$ smoothly truncated outside $\mathrm{B}_3(0)$, Theorem A shows that $$\widetilde{m}(b(\xi)) = \big( 2 - \cos(2 \pi \alpha \xi) - \cos(2 \pi \beta \xi) \big)^\gamma$$ is an $L_p$ multiplier in $\R$. These examples are certainly less standard. With some hindsight, they can be obtained via a clever combination of classical results, we invite the reader to try! However, it seems fair to say that such a general statement follows naturally from our approach, see \cite{PPR,PRo} for related results.

Geometrically, we embed $\R$ in a $2$-dimensional torus as an infinite non-periodic helix.
This geodesic flow clearly generalizes by taking cocycles $\R^n \to \R^{2d}$ of the form $$b_\Sigma(\xi) = \bigoplus_{s=1}^d \big( e^{2\pi i \sum_j \xi_j \gamma_j^s} - 1 \big).$$ Further examples arise from mixed ---neither inner nor proper--- cocycles.

\subsection{The noncommutative tori}
\label{SectNCTori}

We now generalize for noncommutative tori the H\"ormander-Mihlin conditions. Given $n \ge 1$ and an $n \times n$ antisymmetric matrix $\Theta$ with entries $0 \le \theta_{ij} < 1$, we define the noncommutative torus with $n$ generators associated to the angle $\Theta$ as the von Neumann algebra $\mathcal{A}_\Theta$ generated by $n$ unitaries $u_1, u_2, \ldots, u_n$ satisfying the relations $u_j u_k = e^{2 \pi i \theta_{jk}} u_k u_j$. Every element of $\mathcal{A}_\Theta$ can be written as an element in the closure of the span of words of the form $w_k = u_1^{k_1} u_2^{k_2} \cdots u_n^{k_n}$ with $k = (k_1, k_2, \ldots, k_n) \in \Z^n$. Moreover, we equip $\mathcal{A}_\Theta$ with the normalized trace $$\tau(f) = \tau \Big( \sum_{k \in \Z^n} \widehat{f}(k) \, w_k \Big) = \widehat{f}(0).$$ The classical $n$-dimensional torus corresponds to $\Theta=0$, so that $\mathcal{A}_0 = L_\infty(\mathbb{T}^n)$. On the other hand, once we have defined $\mathcal{A}_\Theta$, it is clear what should be the aspect of the heat semigroup for noncommutative tori. Namely $$S_{\Theta,t}(f) = S_{\Theta,t} \Big( \sum_{k \in \Z^n} \widehat{f}(k) \, w_k \Big) = \sum_{k \in \Z^n} \widehat{f}(k) \, e^{-t |k|^2} \, w_k.$$ We may not apply directly any of our results in Section \ref{SectHM} since $\mathcal{A}_\Theta$ is not the group von Neumann algebra of a discrete group. We will use instead that $\mathcal{A}_\Theta$ embeds in the von Neumann algebra of a discretized Heisenberg group. Given an antisymmetric $n \times n$ matrix $\Theta$ with entries $0 \le \theta_{jk} < 1$, consider the bilinear form $B_\Theta: \Z^n \times \Z^n \to \R$ given by $B_\Theta(\xi, \zeta) = \frac12 \sum_{j,k=1}^n \theta_{jk} \xi_j \zeta_k = \frac12 \langle \xi, \Theta \zeta \rangle$. Define the discretized Heisenberg group $\mathrm{H}_\Theta = \R \times \Z^n$ with the product $$(x,\xi) \cdot (z,\zeta) = \big( x + z + B_\Theta(\xi,\zeta), \xi + \zeta \big).$$

\begin{lemma} \label{LemNCTn}
We have $$\mathcal{L}(\mathrm{H}_\Theta) = \int_{\R}^{\oplus} \mathcal{A}_{x \Theta} \, dx.$$
\end{lemma}

\dem Let $\lambda$ denote the left regular representation of $\mathrm{H}_\Theta$. Since $(x,0)$ commutes in $\mathrm{H}_\Theta$ with every $(z, \zeta)$, it turns out that $\lambda(\R,0)$ lives in the center of the algebra $\mathcal{L}(\mathrm{H}_\Theta)$. Using von Neumann's decomposition theorem for subalgebras of the center $$\mathcal{L}(\mathrm{H}_\Theta) = \int_{\mathrm{sp}(\lambda(\R))}^\oplus \mathcal{M}_x \, dx = \int_\R^\oplus \mathcal{M}_x \, dx.$$ Given $\xi \in \Z^n$, we set $w_\xi = \lambda(0,\xi)$ and observe that $w_\xi w_\zeta = \lambda(B_\Theta(\xi,\zeta),\xi+\zeta)$ implies $w_\xi w_\zeta = \lambda \big( B_\Theta(\xi,\zeta) - B_\Theta(\zeta,\xi),0 \big) w_\zeta w_\xi$. The $w_\xi$'s are generated by the unitaries $u_j = \lambda(0,e_j)$ which satisfy $u_j u_k = e^{2 \pi i \theta_{jk} \cdot} u_ku_j$. Moreover, since $\lambda(\R)$ acts on $\mathcal{M}_x$ by scalar multiplication we see that $\mathcal{M}_x = \langle u_j(x) \, | \ 1 \le j \le n \rangle$ where the $u_j(x)$'s arise from $$u_j = \int_\R^\oplus u_j(x) \, dx \quad \mbox{and satisfy} \quad u_j(x) u_k(x) = e^{2 \pi i \theta_{jk} x} u_k(x) u_j(x).$$ Therefore, we have proved that $\mathcal{M}_x = \mathcal{A}_{x \Theta}$ as expected. \fin

\begin{corollary} \label{CorNCTori}
Given an angle $\Theta$ with $n$ generators, let $$T_m: \sum_{k
\in \Z^n} \widehat{f}(k) \, w_k \mapsto \sum_{k \in \Z^n} m_k \,
\widehat{f}(k) \, w_k$$ be the Fourier multiplier on
$\mathcal{A}_\Theta$ associated to $m: \Z^n \to \C$. Let
$\widetilde{m}: \R^n \to \C$ be a lifting multiplier for $m$, so
that $\widetilde{m}_{|_{\Z^n}} = m$. Assume that $\widetilde{m} \in \mathcal{C}^{[\frac{n}{2}]+1}(\R^n
\setminus\{0\})$ and $$\big| \partial_\xi^\beta \widetilde{m}(\xi)
\big| \, \le \, c_n |\xi|^{- |\beta|} \quad \mbox{for all
multi-index $\beta$ s.t.} \quad |\beta| \le
{\mbox{$[\frac{n}{2}]$}}+1.$$ Then, $T_m: L_p(\mathcal{A}_\Theta)
\stackrel{cb}{\to} L_p(\mathcal{A}_\Theta)$ for all $1 < p <
\infty$ and $T_m: L_\infty(\mathcal{A}_\Theta) \stackrel{cb}{\to}
\mathrm{BMO}_{\T_\Theta}$.
\end{corollary}

\dem Let us consider the heat semigroup
$S_{\Theta,t}(\lambda(x,\xi)) = e^{-t|\xi|^2} \lambda(x,\xi)$ and
also the length function $\psi(x,\xi) = |\xi|^2$ in
$\mathrm{H}_\Theta$. Note that $\T_\Theta = \T_\psi$ in the
terminology of Section \ref{SectHM}. The length function yields to
the non-injective cocycle $\mathrm{H}_\Theta \to \R^n$ given by
$b_\psi(x,\xi) = \xi$. The associated action is trivial since
$$\alpha_{\psi,(z,\zeta)}(\xi) =
\alpha_{\psi,(z,\zeta)}(b_\psi(0,\xi)) = b_\psi((z,\zeta) \cdot
(0,\xi)) - b_\psi(z,\zeta) = \xi.$$ In particular,
$|\alpha_\psi(\mathrm{H}_\Theta)| < \infty$ and we know that $$T_M: \sum_{h \in \mathrm{H}_\Theta}
\widehat{f}(h) \lambda(h) \mapsto \sum_{h \in \mathrm{H}_\Theta}
M_h \widehat{f}(h) \lambda(h)$$ will be completely bounded
$\mathcal{L}(\mathrm{H}_\Theta) \to \mathrm{BMO}_{\T_\Theta}$ as
far as we can find a lifting multiplier $\widetilde{m} \circ
b_\psi(h) = M_h$ satisfying the smoothness condition in the
statement. Note also that the non-injectivity of the cocycle
imposes $M_{(x,\xi)} = M_{(z,\xi)}$ for $x,z \in \R$. However,
this is not a restriction for the multiplier $m_k$ in the
statement since $b_\psi(x,\cdot)$ is injective for any $x$. In
other words, we use $m_k = M_{(0,k)} = \widetilde{m} \circ
b_\psi(0,k) = \widetilde{m}(k)$ as expected. Therefore, since
$$\mbox{$M_{(x,\xi)}$ is $x$-independent} \ \Rightarrow \ T_M
= \int_\R^\oplus {T_m}_{|_{\mathcal{A}_{x\Theta}}} \, dx,$$ we
conclude that $$\esssup_{x \in \R}
\big\| T_m: \mathcal{A}_{x \Theta} \to \mathrm{BMO}_\Theta
\big\|_{cb} < \infty.$$ To show complete boundedness for $x=1$, we restrict
the above inequality to the $\mathrm{C}^*$-algebra generated by
the $u_j$'s, where the $\mathcal{A}_{x\Theta}$-norm is
$x$-continuous in the sense of continuous fields \cite{Rie}. This proves the $L_\infty \to
\mathrm{BMO}$ complete boundedness for $x=1$ by weak-$*$ density. The
$L_p(\mathcal{A}_\Theta) \to L_p(\mathcal{A}_\Theta)$
complete boundedness is proved as usual by interpolation and duality
since $\T_\Theta$ is regular as in Theorem \ref{Interpolation}. \fin

\begin{remark}
\emph{There is an alternative proof of Corollary \ref{CorNCTori} by transference. The authors in \cite{CXY} have recently extended to $\mathcal{A}_\Theta$ several results from classical harmonic analysis on $\mathbb{T}^n$ using this idea. In particular, they prove that cb-multipliers on the quantum $n$-torus are exactly those on the usual $n$-torus with equal cb-norms.}
\end{remark}

\subsection{The free group algebra $\mathcal{L}(\mathbb{F}_n)$}
\label{SubFree}

All that is needed to apply Theorems A and B for the free group is to know the more we can about finite-dimensional cocycles on $\mathbb{F}_n$. These cocycles are easy to classify. It suffices to know $b(g_k)$ and $\alpha_{g_k}$ for the generators $g_k$, but any choice of points and unitaries in $\R^d$ is admissible by freeness. Thus the family of finite-dimensional cocycles of $\mathbb{F}_n$ is too rich. We will concentrate on describing low dimensional injective cocycles since they can be regarded as basic building blocks of our family. Since $\mathbb{F}_n$ embeds isomorphically into $\mathbb{F}_2$ for all $n \ge 2$ let us just consider the free group $\mathbb{F}_2$ with two generators $a_1, a_2$. The construction below is well-known to group/measure theorists. Our first observation goes back to the proof of the Banach-Tarski paradox. Namely, if $\theta \in \R \setminus 2 \pi \mathbb{Q}$ the subgroup of $SO(3)$ generated by $$A_1 = \left( \begin{array}{ccc} \cos \theta & - \sin \theta & 0 \\ \sin \theta & \hskip10pt \cos \theta & 0 \\ 0 & 0 & 1 \end{array} \right) \quad \mbox{and} \quad A_2 = \left( \begin{array}{ccc} 1 & 0 & 0 \\ 0 & \cos \theta & - \sin \theta \\ 0 & \sin \theta & \hskip10pt \cos \theta \end{array} \right)$$ is isomorphic to $\mathbb{F}_2$ under the mapping $$\mathbb{F}_2 \ni \underbrace{a_{k_1}^{n_1} a_{k_2}^{n_2} \cdots a_{k_r}^{n_r}}_w \mapsto \underbrace{A_{k_1}^{n_1} A_{k_2}^{n_2} \cdots A_{k_r}^{n_r}}_{W_\theta} \in SO(3)$$ with $k_1, k_2,
\ldots, k_r \in \{1,2\}$, $k_j \neq k_{j+1}$ and $n_1, n_2, \ldots, n_r \in \Z$. On the other hand $SO(3)$ acts naturally on $\R^3$ and $\alpha_{\theta}(w) = W_\theta$ defines an isometric action $\mathbb{F}_2 \curvearrowright \R^3$ with associated cocycle map $b_{\theta\xi}(w) = W_\theta(\xi) - \xi$ for some $\xi
\in \H_\theta = \R^3$. Therefore we find a $3$-dimensional cocycle $(\H_\theta, b_{\theta\xi}, \alpha_\theta)$ for any $\xi \in \R^3$. In order to pick $\xi$ so that $b_{\theta\xi}$ is injective we must show that $$A_\theta = \bigcap_{w \in \mathbb{F}_2 \setminus \{e\}} \big\{ \gamma \in \R^3 \, | \, W_\theta(\gamma) \neq \gamma \big\}$$ is nonempty. However, given $w \in \mathbb{F}_2 \setminus \{e\}$, the orthogonal map $W_\theta$ is a nonidentity linear map on $\R^n$. In particular, the Lebesgue measure of $\R^3 \setminus A_\theta$ is zero since it is a countable union of linear subspaces with codimension at least $1$. This proves that the action $\alpha_\theta$ is weakly free with respect to almost every $\xi \in \R^3$ and for all such $\xi$'s we find an injective $b_{\theta\xi}: \mathbb{F}_2 \to \R^3$. Our construction above is not completely constructive since we have not provided a criterium to pick the right $\xi$'s. If $e_1, e_2, e_3$ denotes the standard basis of $\R^3$, this can be fixed taking $\H_\theta = \R^9$ and
\begin{eqnarray*}
\alpha_\theta(w) & = & W_\theta \oplus W_\theta \oplus W_\theta,
\\ b_\theta(w) & = & \big( W_\theta(e_1) - e_1 \big) \oplus \big(
W_\theta(e_2) - e_2 \big) \oplus \big( W_\theta(e_3) - e_3 \big).
\end{eqnarray*}

\begin{corollary} \label{CorFree}
Given $\theta \in \R \setminus 2 \pi \mathbb{Q}$, consider the free group algebra $\mathcal{L}(\mathbb{F}_2)$ equipped with the cocycle $(\H_\theta, b_\theta, \alpha_\theta)$ above. Let $\psi_\theta$ denote the associated length function and fix a function $\widetilde{m} \in \mathcal{C}^{5}(\R^9 \setminus \{0\})$ with $$\big| \partial_\xi^\beta \widetilde{m}(\xi) \big| \, \lesssim \, |\xi|^{-|\beta|+\varepsilon} \quad \mbox{for all} \quad |\beta| \le 5.$$ Then, if $m: \mathbb{F}_2 \to \C$ is of the form $m_w = \widetilde{m} \circ b_\theta(w)$ $$T_m: \, \summ_w \widehat{f}(w) \lambda(w) \mapsto \summ_w m_w \widehat{f}(w) \lambda(w)$$ defines a cb-map on $L_p(\mathcal{L}(\mathbb{F}_2),\tau)$ for $1 < p < \infty$ and also from $\mathcal{L}(\mathbb{F}_2)$ to $\mathrm{BMO}_{\T_{\psi_\theta}}$.
\end{corollary}

\dem This is a direct application of Theorem A and Remark \ref{RemThA}. \fin

\vskip-5pt

\noindent Further results follow inspecting the conditions of the other results from Section \ref{SectHM}.

\begin{remark} \label{RemFreeLattice}
\emph{Given the free group $\mathbb{F}_n = \langle g_1, g_2, \ldots, g_n \rangle$, $$g_{k_1}^{r_1} g_{k_2}^{r_2} \cdots g_{k_m}^{r_m} \mapsto \sum_{s=1}^n \Big( \sum_{k_j = s} r_j \Big) e_s$$ defines a (non-injective) $\Z^n$-valued cocycle with respect to the trivial action. It vanishes on a normal subgroup $\mathbb{H}_n$ with $\mathbb{F}_n / \mathbb{H}_n \simeq \Z^n$. Hence, the corresponding semigroup $\mathrm{BMO}_{\T_\psi}$ lives in $\mathcal{L}(\mathbb{F}_n/\mathbb{H}_n) \simeq L_\infty(\mathbb{T}^n)$. Since the associated action is trivial, Theorem A with $\varepsilon=0$ shows that Fourier multipliers on $\mathbb{F}_n$ constant in the cosets of $\mathbb{H}_n$ can be analyzed in terms of the corresponding multiplier in the $n$-dimensional torus. We may also compose the given cocycle with other cocycles of $\Z^n$ to obtain cocycles of $\mathbb{F}_n$. That way, our exotic examples for the $n$-torus can be transferred to produce examples in the free group. In fact, this observation applies for many finitely-generated groups. Indeed, by the Grushko-Neumann theorem any finitely-generated $\G$ factorizes as a finite free product of finitely-generated freely indecomposable groups $\G_1 * \G_2 * \cdots * \G_n$. Thus, the same construction applies as long as all the factors $\G_j$ have independent generators $g_s$ satisfying $$g_{k_1}^{r_1} g_{k_2}^{r_2} \cdots g_{k_m}^{r_m} = e \hskip10pt \Rightarrow \hskip10pt \sum_{k_j = s} r_j = 0 \quad \mbox{for all} \quad 1 \le s \le n.$$}
\end{remark}

\vskip-13pt

\null

\subsection{New quantum metric spaces}

The notion of compact quantum metric space was originally introduced by Rieffel \cite{Ri1,Ri2}. Let $\mathcal{A}$ be a $\mathrm{C}^*$-algebra and $\mathcal{B}$ a unital, dense $*$-subalgebra of $\mathcal{A}$. Let $\| \cdot \|_{\mathrm{lip}}$ be a seminorm on ${\mathcal B}$ vanishing exactly on $\C \1_{\A}$. The triple $(\A,\mathcal{B}, \|\cdot\|_{\mathrm{lip}})$ is called a compact quantum metric space if the metric
$\rho(\phi_1,\phi_2) = \sup \{ | \phi_1(x) - \phi_2(x) | \ | \ x \in \mathcal{B} \ \mbox{and} \ \|x\|_{\mathrm{lip}} \le 1 \}$ coincides with the weak-$*$ topology on the state space $S(\A)$. This crucial property is hard to verify in general. Ozawa and Rieffel have found an equivalent condition to this property \cite[Proposition 1.3]{OR}, we rewrite it as a lemma.

\begin{lemma} \label{LemOR}
If $\sigma$ is a state on $\A$ and $$\Big\{ x \in \mathcal{B} \ \mbox{such that} \ \|x\|_{\mathrm{lip}} \le 1 \ \mbox{and} \ \sigma(x)=0 \Big\}$$ is relatively compact in $\A$, then $(\A, \mathcal{B}, \|\cdot\|_{\mathrm{lip}})$ is a compact quantum metric space.
\end{lemma}

Given a length function $\psi: \G \to \R_+$, let $(\H_\psi, \alpha_\psi, b_\psi)$ be either the left or right cocycle associated to it. We will say that $\psi$ yields a well-separated metric if we have $$\Delta_\psi \, = \, \inf_{b_\psi(g) \neq 0} \psi(g) \, = \, \inf_{b_\psi(g) \neq b_\psi(h)} \big\| b_\psi(g) - b_\psi(h) \big\|_{\H_\psi}^2 \, > \, 0.$$

\begin{lemma} \label{CountingLemma}
If $\dim \H_\psi = n$ and $\Delta_\psi > 0$, we have
$$\Big| \Big\{ b_\psi(g) \, \big| \ k \Delta_\psi \le \ |b_\psi(g)| \le (k+1) \Delta_\psi \Big\} \Big| \, \le \, 5^n \, k^{n-1} \,.$$
\end{lemma}

\dem If $\xi_1 \neq \xi_2$ belong to $b_\psi(\G)$, we have $$\Big( \xi_1 + \frac{\Delta_\psi}{2} \mathrm{B}_n \Big) \cap \Big( \xi_2 + \frac{\Delta_\psi}{2} \mathrm{B}_n \Big) = \emptyset,$$ where $\mathrm{B}_n$ denotes the Euclidean unit ball in $\H_\psi$. This shows that 
\begin{eqnarray*}
\lefteqn{\hskip-5pt \Big| \Big\{ b_\psi(g) \, \big| \ k \Delta_\psi \le \ |b_\psi(g)| \le (k+1) \Delta_\psi \Big\} \Big|} \\ & \le & \Big| \frac{\Delta_\psi}{2} \mathrm{B}_n \Big|^{-1} \Big[ \Big| \Big( (k+1) \Delta_\psi + \frac{\Delta_\psi}{2} \Big) \mathrm{B}_n \Big| - \Big| \Big( k \Delta_\psi - \frac{\Delta_\psi}{2} \Big) \mathrm{B}_n \Big| \Big] \\ & = & (2k+3)^n - (2k-1)^n \ = \ \sum_{j=0}^n {{n}\choose{j}} (2k)^{n-j} \big( 3^j - (-1)^j \big) \ \le \ 5^n k^{n-1} \, . \ \ \square 
\end{eqnarray*}

Consider now a discrete group $\G$ equipped with a length function $\psi$. We have noticed above that $\G_0 = \{g \in \G \, | \ \psi(g) = 0\}$ is a subgroup of $\G$. If we consider the usual semigroup $\T_\psi$ given by $S_{\psi,t}(\lambda(g)) = e^{-t \psi(g)} \lambda(g)$, recall that $$L_p^\circ (\widehat{\mathbb{G}}) = \Big\{ f \in L_p(\widehat{\mathbb{G}}) \, \big| \, \lim_{t \to \infty} S_{\psi,t} f = 0 \Big\} = \Big\{ f \in L_p(\widehat{\mathbb{G}}) \, \big| \, \widehat{f}(g) = 0 \ \mbox{for all} \ g \in \G_0 \Big\}.$$

\begin{lemma} \label{LemWS}
If $\dim \H_\psi = n$, $\Delta_\psi > 0$, $|\G_0| < \infty$ and $$\big| \widetilde{m}(\xi) \big| \, \le \, c_n |\xi|^{-(n+\varepsilon)} \ \ \mbox{for some} \ \ \varepsilon > 0,$$ then $T_{\widetilde{m} \circ b_\psi}: L_1(\widehat{\mathbb{G}}) \to L_\infty(\widehat{\mathbb{G}})$ is $cb$-bounded. In particular $$\big\| S_{\psi,t}:
L_1^{\circ}(\widehat{\mathbb{G}}) \to L_\infty(\widehat{\mathbb{G}}) \big\|_{cb} \, \le \, \frac{c_n(\Delta_\psi)}{t^{n/2}} \, .$$
\end{lemma}

\dem Given $f = \sum_g \widehat{f}(g) \otimes \lambda(g) \in S_1^r(L_1(\widehat{\mathbb{G}}))$ with $\widehat{f}(g) \in M_r$, we have
\begin{eqnarray*}
\lefteqn{\hskip-5pt \|f\|_{S_1^r(L_1(\widehat{\mathbb{G}}))} \ = \ \sup_{\|f'\|_{S_\infty^r(L_\infty(\widehat{\mathbb{G}}))} \le 1} \mathrm{tr} \otimes \tau \big( f^*f' \big)} \\ & \ge & \sup_{\|a_g \otimes \lambda(g) \|_{S_\infty^r(L_\infty(\widehat{\mathbb{G}}))} \le 1} \mathrm{tr} \big( \widehat{f}(g)^* a_g \big) \, = \, \| \widehat{f}(g) \|_{S_1^r} \, \ge \, \big\| \widehat{f}(g) \otimes \lambda(g) \big\|_{S_1^r(L_{\infty}(\widehat{\mathbb{G}}))}.
\end{eqnarray*}
This, together with the fact $|\{ g \in \G \, | \, b_\psi(g) = \xi \}| = |\G_0|$ for all $\xi \in b_\psi(\G)$, yield
\begin{eqnarray*}
\lefteqn{\Big\| \sum_{g \in \G} m_g \widehat{f}(g) \otimes \lambda(g) \Big\|_{S_1^r(L_{\infty}(\widehat{\mathbb{G}}))}} \\ [3pt] \!\!\! & \le &
\!\!\! |m_e| \, \Big\| \sum_{\psi(g) = 0} \widehat{f}(g) \otimes \lambda(g) \Big\|_{S_1^r(L_{\infty}(\widehat{\mathbb{G}}))} \\ \!\!\! & + & \!\!\! \sum_{k \ge 1}
\sum_{\begin{subarray}{c} \xi \in b_\psi(\G) \\ k \Delta_\psi
\le |\xi| < (k+1) \Delta_\psi \end{subarray}}
|\widetilde{m}(\xi)| \, \Big\| \sum_{b_\psi(g) = \xi}
\widehat{f}(g) \otimes \lambda(g) \Big\|_{S_1^r(L_{\infty}(\widehat{\mathbb{G}}))} \\ \!\!\! & \le & \!\!\!
\Big( |m_e| + \sum_{k \ge 1} \sum_{\begin{subarray}{c} \xi \in
b_\psi(\G) \\ k \Delta_\psi \le |\xi| < (k+1) \Delta_\psi
\end{subarray}} |\widetilde{m}(\xi)| \Big) \, |\G_0| \, \|f\|_{S_1^r(L_1(\widehat{\mathbb{G}}))} \\ \!\!\! & \le & \!\!\! 5^n \Big( |m_e| + \sum_{k \ge 1} k^{n-1} (k \Delta_\psi)^{-n+\varepsilon} \Big) \, |\G_0| \, \|f\|_{S_1^r(L_1(\widehat{\mathbb{G}}))} \, = \, c_{n,\varepsilon}(\Delta_\psi) \, |\G_0| \, \|f\|_1.
\end{eqnarray*}
The third inequality follows from Lemma \ref{CountingLemma} and our growth assumption on $\widetilde{m}$. The second assertion follows similarly. Indeed, since $f \in L_1^\circ (\widehat{\mathbb{G}})$ we may ignore the term $|m_e|$ above and the result follows from the inequality \vskip-10pt \null \hfill $\displaystyle \summ_k k^{n-1} e^{-tk^2} \le  C(n) \Gamma(n/2) t^{-n/2}$. \hfill $\square$

To state our next result, we need to consider the gradient form associated to the infinitesimal generator $A_\psi(\lambda(g)) = \psi(g) \lambda(g)$ of our semigroup $\T_\psi$. Namely, if $\C[\G]$ stands for the algebra of trigonometric polynomials (whose norm closure is the reduced $\mathrm{C}^*$-algebra of $\G$), we set for $f_1, f_2 \in \C[\G]$ $$2\Gamma(f_1, f_2) = A_\psi(f_1^*) f_2 + f_1^* A_\psi(f_2) - A_\psi(f_1^* f_2).$$ Consider the seminorm $$\|f\|_\Gamma = \max \Big\{ \big\| \Gamma(f,f) \big\|_\infty^\frac12, \big\| \Gamma(f^*,f^*) \big\|_\infty^\frac12 \Big\}$$ and the pseudo-metric $\mathrm{dist}_\psi(g,h) = \sqrt{\psi(g^{-1}h)}$. We find the following result.

\begin{corollary} \label{EXQMS}
If $\dim \H_\psi < \infty$, the following implication holds $$\mathrm{dist}_\psi \ \mbox{well-separated metric} \ \Rightarrow \ (\mathrm{C}_{\mathrm{red}}^*(\G), \C[\G], \| \cdot \|_\Gamma) \ \mbox{quantum metric space}.$$
\end{corollary}

\dem Since $\mathrm{dist}_\psi(g,h) = \| b_\psi(g) - b_\psi(h) \|_{\H_\psi}$, it defines a metric iff $b_\psi: \G \to \H_\psi$ is injective iff $\G_0 = \{e\}$. On the other hand, recalling that $\Gamma(f,f) \ge 0$, we see that $\Gamma(f,f) = 0$ iff $\tau (\Gamma(f,f)) = 0$. It is easily checked that $\tau(\Gamma(f,f)) = \sum_g |\widehat{f}(g)|^2 \psi(g)$. Hence we deduce that $\| \cdot \|_\Gamma$ vanishes in $\C \mathbf{1}$ iff $\G_0 = \{e\}$ iff $\mathrm{dist}_\psi$ is a metric. It is also clear that $\mathrm{dist}_\psi$ is well-separated iff $\Delta_\psi > 0$. In particular, we can not have infinitely many points of $b_\psi(\G)$ inside any ball of the finite-dimensional Hilbert space $\H_\psi$. This means that the set $\{\psi(g)^{-1} \, | g \neq e \}$ can not have a cluster point different from $0$, so that $$A_\psi^{-1}: L_2^\circ(\widehat{\mathbb{G}}) \to L_2^\circ(\widehat{\mathbb{G}})$$ is a compact operator. According to \cite[Theorem 1.1.7]{JM}, $A_\psi^{-1/2}: L_p^\circ (\widehat{\mathbb{G}}) \to L_\infty^\circ (\widehat{\mathbb{G}})$ is also compact for any $p > n+ \varepsilon$. Lemma \ref{LemWS} has been essential at this point, see \cite{JM}. This means that $$\Big\{ f \in L_\infty^\circ (\widehat{\mathbb{G}}) \, \big| \, \big\| A_\psi^{1/2} f \big\|_p \le 1 \Big\} = \Big\{ f \in L_\infty(\widehat{\mathbb{G}}) \, \big| \, \big\| A_\psi^{1/2} f \big\|_p \le 1 \ \mbox{and} \ \tau(f) = 0 \Big\}$$ is relatively compact in $L_\infty(\widehat{\mathbb{G}})$. According to the main result in \cite{JM}, we see that
$$\big\| A_\psi^{1/2} f \big\|_p \ \le \ c_p \max \Big\{ \big\| \Gamma(f,f)^\frac12 \big\|_p, \big\| \Gamma(f^*,f^*)^\frac12 \big\|_p \Big\} \ \le \ c_p \, \|f\|_\Gamma.$$ We deduce from this inequality that $$\Big\{ f \in L_\infty(\widehat{\mathbb{G}}) \, \big| \, \|f\|_\Gamma \le 1 \ \mbox{and} \ \tau(f) = 0 \Big\}$$ is relatively compact in $L_\infty(\widehat{\mathbb{G}})$. The desired result follows from Lemma \ref{LemOR}. \fin

\begin {remark}
\emph{Let $\G = \Z$ and $\psi(k) = |k|^2$. Consider the commutator $[A_\psi^\alpha,f]$ of $A_\psi^\alpha$ and $f\in \C[\Z]$. Rieffel \cite{Ri2} showed that the triple $(\mathrm{C}_{\mathrm{red}}^*(\Z), \C[\Z], \| [A_\psi^\alpha, \cdot \, ] \|)$ is a compact quantum metric space for all $0 < \alpha \le \frac12$. The same argument of the previous corollary shows that this is true for $\frac12 < \alpha \le 1$ too. Indeed, in this case $n=1$ and applying Lemma \ref{LemWS} together with \cite[Theorem 1.1.7]{JM}, we have that $A_\psi^{-\alpha}$ is compact from $L_2^\circ(\mathbb{T})$ to $L_\infty^\circ(\mathbb{T})$ since we may choose $\varepsilon > 0$ such that
$2 > \frac{1+\varepsilon}{2\alpha}$ for any $\alpha > \frac14$. In particular, $\{ x \in L_\infty^\circ(\mathbb{T}) \, | \, \| A_\psi^{\alpha} f \|_2 \le 1 \}$ is relatively compact in $L_\infty (\mathbb{T}).$ Note that
$\| [A_\psi^\alpha, f] \| \, \ge \, \| [A_\psi^\alpha,f] \mathbf{1} \|_2 \, = \, \| A_\psi^\alpha(f) \|_2$. We conclude that $$\Big\{ f \in L_\infty(\mathbb{T}) \, \big| \, \big\| [A_\psi^\alpha, f] \big\| \le 1 \ \mbox{and} \ \int_\mathbb{T} f d\mu = 0 \Big\}$$ is relatively compact. Again, we deduce the assertion from Lemma \ref{LemOR}. Moreover, the exact same argument applies on $\Z^2$ with $\frac12 < \alpha \le 1$ and on $\Z^3$ with $\frac34 < \alpha \le 1$.}
\end{remark}

\section{{\bf Conclusions}}

The classical form of H\"ormander-Mihlin theorem on the compact dual of $\Z^n$ is applied either for \emph{testing} boundedness of a given multiplier or for \emph{constructing} multipliers out of smooth lifting functions. In both situations the standard length $\psi(k) = |k|^2$ with its associated cocycle $b_\psi: \Z^n \hookrightarrow \R^n$ are used in conjunction with transference. The properties which characterize this cocycle are $\Delta_\psi >0$ and the injectivity of the cocycle map $b_\psi$. The injectivity avoids additional restrictions on the multiplier $m$ under the lifting $m = \widetilde{m} \circ b_\psi$, while the well-separatedness $\Delta_\psi > 0$ preserves the discrete topology of $\Z^n$ in its image on $\R^n$. Here is a description of those finite-dimensional cocycles of $\Z^n$.

\begin{lemma} \label{LemStandZn}
Let $\psi$ be a length function on $\Z^n$ giving rise to a
finite-dimensional cocycle $(\H_\psi, \alpha_\psi, b_\psi)$.
Assume that $\dim \H_\psi = d$, $b_\psi$ is injective and
$\Delta_\psi > 0$, then we have
\begin{itemize}
\item $\alpha_\psi: \Z^n \to \mathrm{Aut}(\H_\psi)$ is the
trivial action,

\item $d=n$, $\H_\psi \simeq \R^n$ and $b_\psi: \Z^n \to \H_\psi$
is a group homomorphism.
\end{itemize}
\end{lemma}

\dem We know from Paragraph \ref{SubFMTN} that $$b_\psi(k) \, = \, b_\psi(k_1 \oplus k_2) \, = \, \big( \pi(k_1) \eta - \eta \big) \oplus \gamma(k_2),$$ where $(n,d) = (n_1,d_1) + (n_2,d_2)$ and $n_1=0$ iff $d_1=0$, the map $\pi: \Z^{n_1} \to O(d_1)$ is an orthogonal representation, $\gamma: \Z^{n_2} \to \R^{d_2}$ is a group homomorphism, $\eta \in \R^{d_1}$ and the action has the form $\alpha_\psi(k) = \pi(k_1) \oplus id_{\R^{d_2}}$. Moreover, we claim that $(n_1, d_1) =
(0,0)$ from the hypotheses. Indeed, assume that $n_1, d_1 > 0$. Then, the injectivity of $b_\psi$ and the condition $\Delta_\psi > 0$ imply that $\{ \pi(k_1) \eta - \eta \, | \, k_1 \in \Z^{n_1} \}$ is an infinite set of points in $\R^{d_1}$ mutually separated by a distance greater than or equal to $\sqrt{\Delta_\psi} > 0$. This means that the set must be unbounded, which is a contradiction since $\|\pi(k_1) \eta - \eta \|_{\H_\psi} \le 2 \|\eta\|_{\H_\psi}$ for all $k_1 \in \Z^{n_1}$. 
Once we know $(n_1, d_1) = (0,0)$, the action $\alpha_\psi$ must be trivial if $b_\psi$ has no inner part and we get that $b_\psi = \gamma$ is a group homomorphism. We also know that $d \le n$ since $b_\psi(e_1), b_\psi(e_2), \ldots, b_\psi(e_n)$ linearly generate $\H_\psi$. To prove that $d \ge n$, we use the injectivity of $b_\psi$ and simultaneous Dirichlet's Diophantine approximation. Namely, if $d < n$ we may find $\beta = (\beta_j) \in \R^n \setminus \{0\}$ such that $\sum_j \beta_j b_\psi(e_j) = 0$. This does not contradict the injectivity of $b_\psi$ since the $\beta_j$'s are not necessarily integers. However, given any $\mathrm{N} > 1$ and by Dirichlet approximation, we may find $p_1(\mathrm{N}), p_2(\mathrm{N}), \ldots, p_n(\mathrm{N}), q(\mathrm{N}) \in \Z$ such that $$\Big| \beta_j - \frac{p_j(\mathrm{N})}{q(\mathrm{N})} \Big| \, \le \, \frac{1}{q(\mathrm{N})^\mathrm{N} \mathrm{N}^{\frac{1}{n}}}.$$ This implies that $$\Big| b_\psi \Big( \sum_{j=1}^n p_j(\mathrm{N}) e_j \Big) \Big| \, \le \, \frac{1}{q(\mathrm{N})^{\mathrm{N}-1} \mathrm{N}^{\frac{1}{n}}} \sum_{j=1}^n |b_\psi(e_j)| \, \longrightarrow \, 0 \quad \mbox{as} \ \mathrm{N} \to \infty.$$ Note however that this contradicts the condition $\Delta_\psi > 0$ of well-separatedness. \fin

In particular, given a length function $\psi: \G \to \R_+$, it is natural to call $\psi$ a \emph{standard length function} if $\Delta_\psi > 0$ and the associated cocycle $b_\psi: \G \to \H_\psi$ is injective. Similarly, we will say that a cocycle is standard when so is the length function it gives rise to. Although standard cocycles are an important piece of the theory, they are definitely not the whole of it! We have already illustrated this with our \lq\lq donut multipliers" above. Let us analyze what new information can be extracted from our results so far.

\subsection{Small dimension vs smooth interpolation}
\label{Ending}

If we are given a fixed Fourier multiplier on $\mathcal{L}(\G)$, the problem of finding the optimal cocycle and lifting multiplier to
study its $L_p$ boundedness might be quite hard.

\begin{problem} \label{ProbInj}
Given a Fourier multiplier $$\summ_g \widehat{f}(g) \lambda(g) \mapsto \summ_g m_g \widehat{f}(g) \lambda(g)$$
\begin{itemize}
\item[a)] Find low dimensional injective cocycles $b_\psi: \G \to \H_\psi$.

\item[b)] Given such $(\H_\psi, b_\psi, \alpha_\psi)$, find $\widetilde{m} \in \mathcal{C}^{[\frac{n}{2}]+1}(\R^n \setminus \{0\})$ with $\widetilde{m} (b_\psi(g)) = m_g$ \\ and minimizing $$\hskip10pt \sup_{\xi \in \R^n} \sup_{\begin{subarray}{c} |\beta| \le [\frac{n}{2}]+1 \end{subarray}} \max \big\{ |\xi|^{|\beta|+\varepsilon}, |\xi|^{|\beta|-\varepsilon} \big\} \, \big| \partial_\xi^\beta \widetilde{m}(\xi) \big|,$$ where $\varepsilon$ may be $0$ under any of the situations considered in Paragraph \emph{\ref{SubAG}}.
\end{itemize}
\end{problem}

Once fixed a cocycle, we must
find a lifting multiplier for $m: \G \to \C$ optimizing the
constants. This means that we have to control a number of
derivatives of a smooth
function $\widetilde{m}$ taking certain preassigned values on a cloud of points $b_\psi(g)$ in $\R^n$. In particular, this fits in Fefferman's approach to the smooth interpolation of data carried out in \cite{Fe1,Fe2,Fe3}. There is no canonical answer for questions a) and b) and in general we find certain incompatibility. Indeed, if we pick $\gamma_1, \gamma_2, \ldots, \gamma_n \in \R$ linearly independent over $\Z$, the cocycle $k \mapsto \summ_j \gamma_j k_j$ associated to the trivial action is a $1$-dimensional injective cocycle for $\Z^n$. This minimizes the number of derivatives to estimate for the lifting multiplier. Note however that $\{ \sum_j \gamma_j k_j \, | \, k \in \Z^n \}$ is a dense cloud of points in $\R$ and the $\psi$-metric is far from being well-separated. In general this makes harder to solve b), and the lifting multiplier will be highly oscillating in many cases. On the other hand, as we have
seen for $\R^n$, certain multipliers can only be treated with alternative cocycles like this one. In summary our notion of \lq\lq smooth multiplier" is very much affected by the cocycle we use.

\begin{problem}
Solve Problem $\ref{ProbInj}$ using standard cocycles, not just injective ones.
\end{problem}

This is more restrictive and we will not always find finite-dimensional standard cocycles, see the next paragraph. On the other hand, if we content ourselves with not necessarily injective well-separated cocycles, we may apply our construction in Remark \ref{RemFreeLattice} for finitely generated groups.

\subsection{Infinite-dimensional standard cocycles}

There exists two distinguished length functions on $\Z$, the absolute value $\psi(k) = |k|$ and its square respectively related to the Poisson and heat semigroups. Both yield standard cocycles, but one of them is infinite-dimensional while the other has dimension $1$. If we take free products of $\Z$ only the Poisson like cocycle survives. Hence we wonder if there exist finite-dimensional standard cocycles for the free group. A negative answer follows from a classical theorem of Bieberbach \cite{Bie}. Let us recall that a group $\G$ is called virtually abelian whenever it has an abelian subgroup $\mathrm{H}$ of finite index, so that $\G$ has finitely many left/right $\mathrm{H}$-cosets. Bieberbach's theorem claims that every discrete subgroup of $\R^n \rtimes O(n)$ is virtually abelian.

\begin{theorem} \label{IDSC}
If $\G$ has a finite-dimensional standard cocycle, $\G$ is virtually abelian.
\end{theorem}

\dem Note that $g \mapsto (b_\psi(g), \alpha_{\psi,g}) \in \H_\psi \rtimes O(\dim \H_\psi)$ defines an injective group isomorphism for any standard cocycle. Moreover, the well-separatedness property shows that it is an homeomorphism. Thus, $\G$ can be regarded as a discrete subgroup of $\H_\psi \rtimes O(\dim \H_\psi)$ \hskip-2pt with $\dim \H_\psi < \infty$. \hskip-3pt We conclude from Bieberbach's theorem. \fin

According to this result, we see in particular that nonabelian free groups do not admit finite-dimensional standard cocycles. A unitary representation of a locally compact group $\G$ is called primary if the center of its intertwining algebra $\mathcal{C}(\pi)$ is trivial. The group $\G$ is said to be of type I whenever the von Neumann algebra $\A_\pi$ generated by every primary representation $\pi$ is a type I factor. This condition turns out to be crucial to admit Plancherel type theorems in terms of irreducible unitary representations, see \cite[Chapter 7]{Fol} for explicit results.

\begin{corollary} \label{IDSC2}
A discrete group is virtually abelian if and only if it is of type \emph{I}.
\end{corollary}

\dem By Thoma's theorem \cite{Th}, a discrete group is type I iff it has a normal abelian subgroup of finite index, hence virtually abelian. On the contrary, if $\G$ is virtually abelian it admits an abelian subgroup $\mathrm{H}$ of finite index. Let us show that we can pick another such $\mathrm{H}$ being a normal subgroup. The map $\gamma: g \mapsto \Lambda_g$ with $\Lambda_g(g' \mathrm{H}) = gg' \mathrm{H}$ defines a group homomorphism between $\G$ and the symmetric group of permutations $\mathcal{S}_{\G/\mathrm{H}}$ on the space of left $\mathrm{H}$-cosets. Its kernel is clearly a normal subgroup of $\G$, which is abelian since it is contained in $\mathrm{H}$ and of finite index since $\G / \ker \gamma \simeq \mathrm{Img} \, \gamma$ is a subgroup of a finite group, hence finite. \fin

A locally compact group $\G$ satisfies Kazhdan's property $(\mathrm{T})$ when the trivial representation is an isolated point in the dual object with the Fell topology. A discrete group $\G$ satisfies this property iff all its cocycles are inner. Moreover, a cocycle is inner iff it is bounded. Hence, infinite groups satisfying Kazhdan property $(\mathrm{T})$ do not admit finite-dimensional standard cocycles. In summary, many interesting discrete groups do not admit a \emph{finite-dimensional \lq\lq standard" H\"ormander-Mihlin theory} as it happens with the integer lattice $\Z^n$. Our results establish a more general theory which includes these cases.

\subsection{Bohr compactification}

We have $$\int_{\hatRd^n} \lambda(\xi) \, d\mu = \int_{\hatRd^n} e^{2 \pi i \langle \xi, x \rangle} \, d\mu(x) =
\delta_{\xi,0}$$ for the Haar measure $\mu$ on $\hatRd^n$. Being a Haar measure on a compact group, it is a
translation invariant probability measure on the Bohr compactification. Therefore it vanishes on every measurable
bounded set of $\R^n$ and $\mu$ is singular to the Lebesgue measure. In fact $$\delta_{\xi,0} = \lim_{t \to \infty} \exp (-t |\xi|^2) = \lim_{t \to \infty} \Big( \frac{\pi}{t} \Big)^{\frac{n}{2}} \int_{\R^n} e^{2 \pi i \langle \xi,x \rangle}
\, \exp \Big( - \frac{\pi^2 |x|^2}{t} \Big) \, dx,$$ so that we find $$\int_{\hatRd^n} f \, d\mu = \lim_{t \to \infty} \Big(
\frac{\pi}{t} \Big)^{\frac{n}{2}} \int_{\R^n} f(x) \, \exp \Big( - \frac{\pi^2 |x|^2}{t} \Big) \, dx.$$ In other words, the measure $\mu$ can be understood as a limit of averages along large balls. By subordination, the same holds for Poisson kernels. As it follows from Theorem D, a dimension-free Calder\'on-Zygmund theory for Fourier multipliers on arbitrary discrete groups would follow from a dimension-free CZ theory on the Bohr compactification. This leads us to the following very natural problem.

\begin{problem} \label{ProbBohr}
Develop a CZ theory for the heat$/$Poisson semigroups on $\widehat{\R}_{\mathrm{disc}}^{\infty}$.
\end{problem}

In order to bring some hope to the problem suggested above, we can construct non-trivial radial Fourier multipliers in the Bohr compactification of $\R^\infty$. In terms of Theorem D, we may equivalently say that the class of radial Fourier multipliers which are bounded $T_{h \circ \psi}: \V \to \mathrm{BMO}_{\T_\psi}(\V)$ for any discrete group $\G$ with $\dim \H_\psi = \infty$ is not trivial. Indeed, as it follows from \cite{JM2}, imaginary powers of length functions are bounded with dimension free constants. More concretely given any discrete group $\G$ and any length function $\psi: \G \to \R_+$, the family of functions of the form $$m_g = \sqrt{\psi(g)} \int_{\R_+} e^{- s \sqrt{\psi(g)}} \, f(s) \, ds$$ with $f: \R_+ \to \C$ bounded, define radial multipliers for which $$T_m: \V \to \mathrm{BMO}_{\T_\psi}$$ is bounded and its norm does not depend on $\dim \H_\psi$. Recall that $$f(s) = \frac{s^{-2i\gamma}}{\Gamma(1 - i\gamma)} \quad \mbox{with} \quad \gamma \in \R \quad \Rightarrow \quad m_g = \psi(g)^{i\gamma}.$$

\subsection{H\"ormander-Mihlin dimension}
\label{SubHMDim}

According to our geometrical analysis, low dimensional injective cocycles are basic building blocks to study the $L_p$ boundedness of Fourier multipliers on compact duals of discrete groups. We may for instance reconstruct (up to an orthogonal change of basis) the standard cocycle $\Z^n \to \R^n$ as the direct sum of $n$ one-dimensional injective cocycles. Thus, given a discrete group $\mathrm{G}$ we define its \emph{H\"ormander-Mihlin dimension} $\mbox{HM-dim} (\G)$ as $$\inf \Big\{ \dim \mathcal{H}_\psi \, \big| \ \psi: \mathrm{G} \to \R_+ \ \mbox{length function with $b_\psi: \mathrm{G} \to \mathcal{H}_\psi$ injective} \Big\}.$$ In other words, ignoring degenerate multipliers which are constant in the cosets of certain subgroup, the H\"ormander-Mihlin dimension gives a lower scale to construct meaningful cocycles. We have already proved that
\begin{itemize}
\item $\mbox{HM-dim}(\Z^n) = 1$ for all $n \ge 1$.

\item $\mbox{HM-dim}(\mathbb{F}_n) \le 3$ for all $2 \le n \le \infty$.
\end{itemize}

\begin{lemma} We have
\begin{itemize}
\item[i)] \emph{HM-dim} is defined for every discrete $\G$.

\vskip3pt

\item[ii)] Given a discrete group $\G$, we have that $$\mbox{\rm HM-dim}(\G) = \infty$$ if $\G$ is finitely generated, non-amenable and does not contain $\mathbb{F}_2$.
\end{itemize}
\end{lemma}

\dem The first assertion claims that every discrete $\G$ has an injective cocycle into a (possibly infinite-dimensional) Hilbert space. Indeed, just take $\H_\psi = \ell_2(\G,\R)$ the space of $\R$-valued square integrable functions on $\G$ with its usual inner product and $b_\psi(g) = \delta_g - \delta_e$, which is naturally implemented by the action of the left regular representation. The second assertion is a little more subtle. Assume that such a $\G$ admits an injective cocycle $(\H, \alpha, b)$ with $\dim \H = n < \infty$. This means that we have an injective group homomorphism $$\pi: g \in \G \mapsto \Big( \begin{array}{cc} \alpha_g & b(g) \\ 0 & 1 \end{array} \Big) \in \mathrm{Aff}(\R^n) \subset \mathrm{GL}_{n+1}(\R),$$ so that $\G$ is a finitely generated subgroup of $\mathrm{GL}_{n+1}(\R)$. By Tits alternative, $\G$ must be either amenable or contain $\mathbb{F}_2$ as a subgroup, a contradiction. Examples of infinite-HM-dimensional groups are therefore the Burnside groups $B(m,n)$ for $m \ge 2$ and $n \ge 665$ odd, see \cite{P4} for more on this topic. \fin

\begin{remark}
\emph{Note also that:}
\begin{itemize}
\item \emph{If $\mbox{{\rm HM}-dim}(\G_k) < \infty$ for $k \ge 1$ we find} $$\mbox{\rm HM-dim} \Big( \prod_{k=1}^n \mathrm{G}_k \Big) \le \sum_{k=1}^n \mbox{\rm HM-dim}(\mathrm{G}_k).$$

\item \emph{Assume $\mbox{\rm HM-dim}(\G), \mbox{\rm HM-dim}(\mathrm{H}) < \infty$ and}

\vskip3pt

\begin{itemize}
\item[$\circ$] \emph{$\psi$ is a length function of $\G$ with} $$\psi(g) \neq 0 \ \mbox{for} \ g \neq e \quad \mbox{and} \quad \dim \H_\psi = \mbox{\rm HM-dim}(\G),$$
\item[$\circ$] \emph{There exists an action $\beta: \mathrm{H} \curvearrowright \G$ such that $\psi \circ \beta_h = \psi$.}
\end{itemize}

\vskip3pt

\noindent \emph{Then we may estimate the H\"ormander dimension of $\G \rtimes_\beta \mathrm{H}$ as follows} $$\mbox{\rm HM-dim} \big( \G \rtimes_\beta \mathrm{H} \big) \le \mbox{\rm HM-dim}(\G) + \mbox{\rm HM-dim}(\mathrm{H}).$$
\end{itemize}
\emph{The proofs of these estimates are straightforward and we leave them to the reader.}
\end{remark}

\begin{remark}
\emph{We have already seen there exist certain discrete groups with no finite-dimensional injective cocycles. In fact, the Tarski monster group $\Gamma$ has no finite-dimensional cocycles at all. It is a simple, finitely generated, non-amenable group which does not contain an isomorphic copy of $\mathbb{F}_2$. Let us assume $\Gamma$ admits an $n$-dimensional cocycle $(\mathcal{H}, \alpha, b)$. As before, we construct a group homomorphism $\pi: \Gamma \to \mathrm{Aff}(\R^n) \subset \mathrm{GL}_{n+1}(\R)$. Since $\Gamma$ is simple, we find that $\mathrm{ker}(\pi)$ is either $\{e\}$ or $\Gamma$, so that $\pi$ must be injective because $b$ is assumed to be non-trivial. Therefore, we conclude as above by using that $\Gamma$ is finitely generated and Tits alternative.}
\end{remark}

Our goal now is to present a brief analysis of the H\"ormander-Mihlin dimension of finite groups. Given a finite group $\G$ equipped with an isometric action $\alpha: \G \to \mathrm{Aut}(\mathcal{H})$, we will say that $\alpha$ is weakly free if there exists $\eta \in \H$ such that $\alpha_g(\eta) \neq \eta$ for all $g \neq e$. Under these circumstances $b(g) = \alpha_g(\eta) - \eta$ defines an injective inner cocycle $\G \to \H$. Taking $\mathcal{H} = \ell_2(\G)$, the left action $\alpha_g(\delta_h) = \delta_{gh}$ is weakly free. This shows in particular that $$\mbox{HM-dim}(\G) \le |\G|$$ for any finite group. This estimate seems quite rough for many groups, but we do not know whether there exist finite groups of arbitrary large cardinality satisfying $\mbox{HM-dim}(\G) \gtrsim |\G|$ up to some universal constant. Note in passing that this property is opposite to factoriality since $$\mbox{HM-dim} \Big( \prod_{k \ge 1} \G_k \Big) \le \sum_{k \ge 1} \mbox{HM-dim}(\G_k) << \Big| \prod_{k \ge 1} \G_k \Big|.$$

\begin{problem}
Find a family of finite groups $(\G_k)_{k \ge 1}$ with $$|\G_k| < |\G_{k+1}| \quad \mbox{and} \quad \lim_{k \to \infty} \frac{\mbox{HM-dim}(\G_k)}{|\G_k|^\gamma}  > 0 \quad \mbox{for some} \quad 0 < \gamma \le 1.$$
\end{problem}

\noindent Let us give two more examples with $\mbox{HM-dim}(\G) << |\G|$:

\begin{itemize}
\item[\textbf{a)}] The cyclic abelian groups $\Z_n$ satisfy
$$\mbox{HM-dim}(\Z_2) = 1 \quad  \mbox{and} \quad
\mbox{HM-dim}(\Z_n) = 2$$ for any $n \ge 3$. Indeed, regarding
$\Z_n$ as the multiplicative group of $n$-th roots of unity in
$\mathbb{T}$, $g \in \Z_n \mapsto g-1 \in \C$ defines an injective
cocycle implemented by the action $\alpha_g(z) = gz$. This shows
that $\mbox{HM-dim}(\Z_2) = 1$ and $\mbox{HM-dim}(\Z_n) \le 2$. If
$\mbox{HM-dim}(\Z_n) = 1$ for some $n \ge 3$ there must exists an
isometric action $\Z_n \curvearrowright \R$. If $n$ is odd the action
must be trivial, so that $b(gh) = b(g) +
\alpha_g(b(h)) = b(g) + b(h)$ which gives $n b(1) = b(n) = b(0) =
0$ and thus we get a non-injective cocycle. When $n$ is even, we
may also consider the action $\alpha_g(z) = (-1)^g z$ which
implies $b(g^2) = b(g) + (-1)^g b(g)$ so that $b(2) = 0$. Thus, we get a non-injective cocycle unless
$n=2$.

\vskip7pt

\item[\textbf{b)}] Given $\Omega = \{1,2,\ldots,n\}$, let $\mathcal{S}_n$ be its symmetric permutation group. If we set $\mathcal{H} = \ell_2(\Omega,\C)$ the Hilbert space of functions
$\Omega \to \C$, we have a natural action $\alpha_w(\delta_k) =
\delta_{w(k)}$ for $w \in \mathcal{S}_n$ and $1 \le k \le n$. This
action is weakly free, since taking $\eta = \sum_{k=1}^n e^{2 \pi
i k/n} \delta_k \ \Rightarrow \ \alpha_w(\eta) \neq \eta$ for $w \neq id_\Omega$. Thus the map $w \mapsto
\alpha_w(\eta) - \eta$ defines an injective cocycle and we have
$$\mbox{HM-dim} (\mathcal{S}_n) \le \dim_\R(\ell_2(\Omega,\C)) =
2n.$$ Moreover, since $\ell_2(\Omega,\C)$ is equipped with its
$\R$-valued inner product $$\psi(w) \ = \ \sum_{k=1}^n \Big| e^{2\pi i w(k)/n} - e^{2\pi i k/n} \Big|^2.$$ In fact, taking $\eta = \sum_k k \delta_k$ we might work with $\H = \ell_2(\Omega,\R)$ giving rise to the estimate $\mbox{HM-dim}(\mathcal{S}_n) \le n$. However, this choice leads to a less natural length function $\psi$. On the other hand, there is another standard length function given by the number of crossings of the permutation $$\psi(w) = \Big| \Big\{ (i,j) \, \big| \quad 1 \le i < j \le n \quad \mbox{s.t.} \quad w(i) > w(j)\Big\} \Big|.$$ This coincides with the minimal number of transpositions which are needed to factorize $w$. As explained to us by Marek Bo\.zejko, this however leads to a cocycle of dimension $\binom{n}{2} = \frac12 n(n-1) >> n$ for $n$ large.

\vskip7pt

\item[\textbf{c)}] What about Thompson, Coxeter, Dihedral... groups? Performing a similar analysis for these other families of groups will provide explicit estimates for the $L_p$-norm of Fourier multipliers in terms of our H\"ormander/Mihlin smoothness conditions.
\end{itemize}

\section*{{\bf Appendix. An $H_1 - \mathrm{BMO}$ duality}}

Let $\psi: \G \to \R_+$ be any length on some given discrete group $\G$ and denote by $\S_\psi = (S_{\psi,t})_{t \ge 0}$ the Markov semigroup $\lambda(g) \mapsto \exp(-t \psi(g)) \lambda(g)$ on $\V$. Define $H_1(\S_\psi)$ as the completion of $L_2^\circ(\widehat{\mathbb{G}})$ with respect to the norm $$\|f\|_{H_1(\S_\psi)} = \sup_{\|h\|_{\mathrm{BMO}_{\S_\psi}} \le 1} |\tau_\G(fh^*)|,$$ with operator space structure determined by $\|f\|_{M_k(H_1(\S_\psi))}  = \|f\|_{\mathcal{CB}(\mathrm{BMO}_{\S_\psi},M_k)}$. The row/column analogues are defined similarly using the row/column forms of $\mathrm{BMO}_{\S_\psi}$ instead. We will prove the following result.

\begin{Atheorem} \label{H1BMODual}
The map $$f \in \mathrm{BMO}_{\S_\psi} \mapsto \phi_f \in H_1(\S_\psi)^*$$ with $\phi_f$ densely defined by $\phi_f(h) = \tau_\G(f^*h)$ for $h \in L_2^\circ(\widehat{\mathbb{G}})$ is a complete isomorphism.
\end{Atheorem}

\dem We will first consider the column case $$H_1^c(\S_\psi)^* = \mathrm{BMO}_{\S_\psi}^c,$$ the row case follows by a similar argument. The inclusion $\mathrm{BMO}_{\S_\psi}^c \subset H_1^c(\S_\psi)^*$ is trivially contractive from the definition of $H_1^c(\S_\psi)$. On the other hand, note that for any $f \in \mathrm{BMO}^c_{\S_\psi} \subset L_2^\circ(\widehat{\mathbb{G}})$, we have $\lim_{t \to \infty} \|S_tf\|_2=0$. Hence, given any $\delta > 0$ we may pick a large $t > 0$ satisfying $$\tau_\G (|f|^2) = \tau_\G(S_{\psi,t}|f|^2) \le (1+\delta) \tau_\G \big( S_{\psi,t} |f|^2 - |S_{\psi,t}f|^2 \big) \le (1+\delta) \|f\|_{\mathrm{BMO}_{\S_\psi}^c}^2.$$ This proves that $\|f\|_2 \le \|f\| _{\mathrm{BMO}_{\S_\psi}^c}$, which implies $$L_2^\circ(\widehat{\mathbb{G}}) \subset H_1^c(\S_\psi)$$ contractively. In particular, any $\phi \in H_1^c(\S_\psi)^*$ is a linear functional on $L_2^\circ(\widehat{\mathbb{G}})$ and $$\phi(h) = \tau_\G (f^*h) = \phi_f(h)$$ for some $f \in L_2^\circ(\widehat{\mathbb{G}})$ and all $h \in L_2^\circ(\widehat{\mathbb{G}})$. We need to show that $f \in \mathrm{BMO}_{\S_\psi}^c$. For this purpose, we will require a minimal amount of $L_p$-module theory. Given $1 \le p < \infty$ recall that $L_p(\V \bar{\otimes}_{S_{\psi,t}} \V)$ is the closure of algebraic tensors $z = \sum_{j=1}^m a_j \otimes b_j$ with respect to the norm $$\|z\|_{p,\psi,t} = \Big( \tau_\G \big( \langle z, z \rangle_{\psi,t}^{\frac{p}{2}} \big) \Big)^\frac1p,$$ where the Hilbert module bracket is given by $$\Big\langle \summ_j a_j \otimes b_j , \summ_k a_k' \otimes b_k' \Big\rangle_{\psi,t} \, = \, \summ_{j,k} b_j^* S_{\psi,t}(a_j^*a_k') b_k'.$$ When $p=\infty$, $L_p(\V \bar{\otimes}_{S_{\psi,t}} \V)$ is the  closure of algebraic tensors with respect to the strong operator topology determined by the seminorms $\varphi(\langle z,z \rangle_{\psi,t})$ with $\varphi \in \V_*$. We refer to \cite{JS,Ps} for the following facts 
\begin{itemize}
\item[i)] $L_\infty(\V \bar{\otimes}_{S_{\psi,t}} \V)$ is the dual space of $L_1(\V \bar{\otimes}_{S_{\psi,t}} \V)$.

\item[ii)] Given $1 \le p \le \infty$, the $L_p$-module $L_p(\V \bar{\otimes}_{S_{\psi,t}} \V)$ is isomorphic to a complemented subspace of the column space $L_p(\V; \ell_2^c(\mathcal{I}))$ for some index set $\mathcal{I}$.
\end{itemize}
In particular, given $x \in L_2(\V \bar{\otimes}_{S_{\psi,t}} \V)$ it will belong to the unit ball of $ L_\infty(\V \bar{\otimes}_{S_{\psi,t}} \V)$ provided $|\tau_\G ( \langle x,z \rangle_{\psi,t})| \le 1$ for all $z = \sum_{j=1}^m a_j \otimes b_j$ with $\|z\|_{1,\psi,t} \le 1$. Let $$u_t(f)=f\otimes \1 - \1 \otimes S_{\psi,t}(f)$$ for $f \in L_p(\widehat{\mathbb{G}})$ so that $$\big\langle u_t(f),u_t(f) \big\rangle_{\psi,t} \, = \, S_{\psi,t}(f^*f) - S_{\psi,t}(f)^*S_{\psi,t}(f).$$ It is clear that $u_t: L_p^\circ(\widehat{\mathbb{G}})\to L_p(\V \bar{\otimes}_{S_{\psi,t}} \V)$ for $2 \le p \le \infty$. Moreover
\begin{eqnarray*}
\tau_\G \big( \langle u_t(f),z \rangle_{\psi,t} \big) \!\! & = & \!\! \summ_j \tau_\G \big( S_{\psi,t}(f^*a_j)b_j \big) - \tau_\G \big( S_{\psi,t}(f^*) S_{\psi,t}(a_j)b_j \big) \\ \!\! & = & \!\! \summ_j \tau_\G \big( f^*(a_jS_{\psi,t}(b_j) - S_{\psi,t}(S_{\psi,t}(a_j)b_j)) \big) \ = \ \tau_\G (f^*u_t^*(z)),
\end{eqnarray*} 
for any $z=\sum_{j=1}^m a_j\otimes b_j$, with $u_t^*(z)=\sum_j a_j S_{\psi,t}(b_j) - S_{\psi,t}(S_{\psi,t}(a_j)b_j))$. In fact if $\Pi$ denotes the $L_2$-projection onto $\mathrm{ker} A_{\psi,2}$, we have $\Pi(f)=0$ since $f$ is a mean-zero element. This implies that 
$$\tau_\G \big( \langle u_t(f),z \rangle_{\psi,t} \big) = \tau_\G \big( f^*(1-\Pi) u^*_t(z) \big).$$ Applying this identity to all $f \in \mathrm{BMO}_{\S_\psi}^c \subset L_2^\circ(\widehat{\mathbb{G}})$, we see that $$(1-\Pi)u_t^*(z) \in H_1^c(\S_\psi) \quad \mbox{and} \quad \big\| (1-\Pi)u_t^*(z) \big\|_{H_1^c(\S_\psi)} \le \|z\|_{1,\psi,t}$$ for all $z\in \sum_{j=1}^m a_j\otimes b_j$ and all $m\in{\Bbb N}$. Now we are ready to go back to the proof. Namely, we have $\phi = \phi_f \in H_1^c(\S_\psi)^*$ and we are interested in showing that $f \in \mathrm{BMO}_{\S_\psi}^c$. Since $f \in L_2(\widehat{\mathbb{G}})$ we know that $u_t(f) \in L_2(\V \bar{\otimes}_{S_{\psi,t}} \V)$ and applying our identity above to $f$ and $z=\sum_{j=1}^m a_j\otimes b_j$ in the unit ball of $L_1(\V \bar{\otimes}_{S_{\psi,t}} \V)$, we get
\begin{eqnarray*}
\big| \tau_\G \big( \langle u_t(f),z \rangle_{\psi,t} \big) \big| \, = \,\big| \tau_\G \big( f^*(1-\Pi)u_t^*(z) \big) \big| \, \le \, \sup_{\begin{subarray}{c} h \in L_2^\circ(\widehat{\mathbb{G}}) \\ \|h\|_{H_1^c(\S_\psi)} \le 1 \end{subarray}} |\tau_\G(f^*h)|.
\end{eqnarray*} 
By the density of $L_2^\circ(\widehat{\mathbb{G}})$ in $H_1^c(\S_\psi)$ and of algebraic tensors in $L_1(\V \bar{\otimes}_{S_{\psi,t}} \V)$ $$\|f\|_{\mathrm{BMO}_{\S_\psi}^c}=\sup_{t > 0} \|u_t(f)\|_{\infty,\psi,t} \le \|\phi_f\|_{H_1^c(\S_\psi)^*}.$$ The same proof works in the operator space setting after matrix amplification. We may similarly define $H_1^r(\S_\psi)$ taking adjoints. Let us now consider the direct sum $\mathrm{X} = L_2^{\circ}(\widehat{\mathbb{G}}) \oplus L_2^{\circ}(\widehat{\mathbb{G}})$ as a subspace of $H_1^r(\S_\psi) \oplus H_1^c(\S_\psi)$. Since $L_2^{\circ}(\widehat{\mathbb{G}})$ is contractively contained and dense in both spaces, we may define the sum $$H_1^r(\S_\psi) + H_1^c(\S_\psi)$$ as the completion of $\mathrm{X}/\Delta$ where $\Delta = \big\{ (f,-f): f \in L_2^{\circ}(\widehat{\mathbb{G}}) \big\}$. It is clear that $L_2^\circ(\widehat{\mathbb{G}})$ is dense in the sum. Moreover, as a consequence of our row/column duality results, it turns out that the dual of this sum is $\mathrm{BMO}_{\S_\psi} = \mathrm{BMO}_{\S_\psi}^r \cap \mathrm{BMO}_{\S_\psi}^c$. Finally, these properties imply that $H_1(\S_\psi) = H_1^r(\S_\psi) + H_1^c(\S_\psi)$. \fin

\begin{Aremark}
\emph{The duality result above holds in the general setting of semigroup type BMO spaces over finite von Neumann algebras, although we have preferred to adapt the terminology to the one used in this paper.}
\end{Aremark}

\begin{Aremark}
\emph{The space $\1\otimes \V \subset \V \bar\otimes_{S_{\psi,t}} \V$ is a $\V$-right module and there exists a completely bounded projection $\mathrm{P}: \V \bar\otimes_{S_{\psi,t}} \V \to \1 \otimes \V$ given by $a \otimes b \mapsto \1 \otimes S_t(a)b$. Consider the set $\mathrm{A}_{c,t}$ of elements $$\xi = \sum_{j=1}^m a_j S_{\psi,t}(b_j)-S_{\psi,t}(S_{\psi,t}(a_j)b_j) \in L_2^{\circ}(\widehat{\mathbb{G}})$$ with $a_j,b_j \in L_2(\widehat{\mathbb{G}})$ and such that $$\Big\| \sum_{j=1}^m a_j \ten b_j \Big\|_{1,\psi,t} \le 1.$$ Note that $\mathrm{A}_{c,t}$ is the image of the unit ball in $L_1(\V \bar\otimes_{S_{\psi,t}} \V)$ via $u_t^*$. Let $$\mathrm{A}_c = \bigcup_{t>0} \mathrm{A}_{c,t}.$$ As in Theorem \ref{H1BMODual}, we then see that $\mathrm{A}_c$ is norming for $\mathrm{BMO}_{\S_\psi}^c$. Given $h \in L_2^\circ(\widehat{\mathbb{G}})$, define $$\|h\|_{H_{1,\mathrm{at}}^c(\S_\psi)} = \inf \Big\{ \summ_k |\lambda_k| : h = \summ_k \lambda_k \xi_k, \ \xi_k \in \mathrm{A}_c \Big\},$$ and let $H_{1,\mathrm{at}}^c(\S_\psi)$ denote the completion of $L_2^\circ(\widehat{\mathbb{G}})$ with respect to the atomic norm given above. The convergence in the identity $h = \summ_k \lambda_k \xi_k$ is understood in the $L_2$ norm. Then we have an isometric isomorphism $$f \in \mathrm{BMO}_{\S_\psi}^c \mapsto \phi_f \in H_{1,\mathrm{at}}^c(\S_\psi)^*.$$ Indeed, since $\mathrm{A}_c$ is norming for $\mathrm{BMO}_{\S_\psi}^c$ we easily obtain an isometric embedding $\mathrm{BMO}_{\S_\psi}^c \hookrightarrow H_{1,\mathrm{at}}^c(\S_\psi)^*$. To prove that it is surjective, it suffices to show that every $\phi \in H_{1,\mathrm{at}}^c(\S_\psi)^*$ arises as $\phi = \phi_f$ for some $f \in \mathrm{BMO}_{\S_\psi}^c$. Arguing as in the proof of Theorem \ref{H1BMODual}, it reduces to show that $$L_2^\circ(\widehat{\mathbb{G}}) \subset H_{1,\mathrm{at}}^c(\S_\psi)$$ contractively. For $f \in L_2^{\circ}(\mathbb{G})$ we have $$\|f \ten \1\|_{1,\psi,t} \le \big\| S_{\psi,t}(f^*f)^{\frac12} \big\|_1 \le \big\| S_{\psi,t}(f^*f)^{\frac12} \big\|_2 \le \|f\|_2,$$ and $$u_t^*(f \ten \1) = f - S_{\psi,t}(S_{\psi,t}(f)) = f - S_{\psi,2t}(f).$$ Applying this to $t/2$ we get $$\|f - S_{\psi,t}(f)\|_{H_{1,\mathrm{at}}^c(\S_\psi)} \le \|f\|_2.$$ Now, since $f \in L_2^\circ(\widehat{\mathbb{G}})$, we may fix $t_k \to \infty$ such that $\|S_{\psi,t_k}(f)\|_2 \le 2^{-k}\varepsilon$. Thus 
\begin{eqnarray*}
\lefteqn{\hskip-10pt \big\| S_{\psi,t_k}(f) - S_{\psi,t_{k+1}}(f) \big\|_{H_{1,\mathrm{at}}^c(\S_\psi)}} \\ & = & \big\| S_{\psi,t_k}(f) - S_{\psi,t_{k+1}-t_k}(S_{\psi,t_k}(f)) \big\|_{H_{1,\mathrm{at}}^c(\S_\psi)} \ \le \ \|S_{\psi,t_k}(f)\|_2 \ \le \ \varepsilon 2^{-k}.
\end{eqnarray*}
This implies 
\begin{eqnarray*}
\lefteqn{\|f\|_{H_{1,\mathrm{at}}^c(\S_\psi)}} \\ \!\!\! & \le & \!\!\! \big\| f-S_{\psi,t_1}(f) \big\|_{H_{1,\mathrm{at}}^c(\S_\psi)} + \sum_{k \ge 1} \big\| S_{\psi,t_k}(f)-S_{\psi,t_{k+1}}(f) \big\|_{H_{1,\mathrm{at}}^c(\S_\psi)} \le (1+\varepsilon) \|f\|_{2}.
\end{eqnarray*}
Letting $\varepsilon \to 0^+$ we complete the proof. This gives the atomic description for $H_1^c(\S_\psi)$ and combining row and column atoms also for $H_1(\S_\psi) = H_1^r(\S_\psi) + H_1^c(\S_\psi)$.}
\end{Aremark} 

\bibliographystyle{amsplain}

\vskip30pt

\hfill \noindent \textbf{Marius Junge} \\
\null \hfill Department of Mathematics
\\ \null \hfill University of Illinois at Urbana-Champaign \\
\null \hfill 1409 W. Green St. Urbana, IL 61891. USA \\
\null \hfill\texttt{junge@math.uiuc.edu}

\

\hfill \noindent \textbf{Tao Mei} \\
\null \hfill Department of Mathematics
\\ \null \hfill Wayne State University \\
\null \hfill 656 W. Kirby Detroit, MI 48202. USA \\
\null \hfill\texttt{mei@wayne.edu}

\

\hfill \noindent \textbf{Javier Parcet} \\
\null \hfill Instituto de Ciencias Matem{\'a}ticas \\ \null \hfill
CSIC-UAM-UC3M-UCM \\ \null \hfill Consejo Superior de
Investigaciones Cient{\'\i}ficas \\ \null \hfill C/ Nicol\'as Cabrera 13-15.
28049, Madrid. Spain \\ \null \hfill\texttt{javier.parcet@icmat.es}
\end{document}